\documentclass[12pt,a4paper,cd]{article}

\usepackage{a4wide,amssymb,amsmath,amsthm,amscd, xspace,epsfig,graphics}
\usepackage{graphicx}

\begin{document}

\newcommand{\RR}{\mathsf{R}}
\newcommand{\LL}{\mathbb{L}}
\newcommand{\OO}{\mathbb{O}}

\def\co{\colon}

\newcommand{\R}{\mathbb{R}}
\newcommand{\Z}{\mathbb{Z}}
\newcommand{\Q}{\mathbb{Q}}
\newcommand{\C}{\mathbb{C}}

\newcommand{\esp}{\vskip .3cm \noindent}

\newcommand\Cfty{${\cal C}^{\infty}$}

\mathchardef\flat="115B

\newcommand\adm{{\subset}^a}
\newcommand\cfp{cofinal${}^+$}
\newcommand\cfm{cofinal${}^-$}
\newcommand\bdp{bounded${}^+$}
\newcommand\bdm{bounded${}^-$}


\def\uptriangle{{\blacktriangle\mkern-3mu\blacktriangle}}


\def\CC#1{${\cal C}^{#1}$}
\def\h#1{\hat #1}
\def\t#1{\tilde #1}
\def\wt#1{\widetilde{#1}}
\def\wh#1{\widehat{#1}}
\def\wb#1{\overline{#1}}

\def\we#1{\text{Up}({#1})}
\def\ve#1{\text{Down}({#1})}
\def\de#1{{\Delta}({#1})}

\def\restrict#1{\bigr|_{#1}}

\newtheorem*{lemma*}{Lemma}
\newtheorem*{thmbis}{Theorem \ref{theorem1}'}

\newtheorem{thmalph}{Theorem}
\renewcommand\thethmalph{\Alph{thmalph}}




\newtheorem{prop}{Proposition}[section]
\newtheorem{ex}[prop]{Example}
\newtheorem{exs}[prop]{Examples}
\newtheorem{thm}[prop]{Theorem}
\newtheorem{thm*}[prop]{$\ast$Theorem}
\newtheorem{cor*}[prop]{$\ast$Corollary}
\newtheorem{lemma}[prop]{Lemma}
\newtheorem{cor}[prop]{Corollary}
\newtheorem{rem}[prop]{Remark}
\newtheorem{defi}[prop]{Definition}

\newtheorem*{defi*}{Definition}
\newtheorem*{exs*}{Examples}
\newtheorem*{conj*}{Conjecture}

\newtheorem{prob}{Problem}
\newtheorem{q}[prob]{Question}

\newtheorem*{prop*}{Proposition}
\newtheorem*{claim}{Claim}
\newtheorem*{lemma0}{Lemma 0}

\renewcommand\theprop{\thesection .\arabic{prop}}

\noindent

\title{Directions in Type I spaces}
\date{\empty}
\author{Mathieu Baillif}
\maketitle

\abstract{A direction in a Type I space $X=\cup_{\alpha<\omega_1}X_\alpha$ is a closed and unbounded subset $D$
of $X$ such that given any continuous $f:X\to\LL_{\ge 0}$ (the closed long ray), if $f$ is unbounded
on $D$ then $f$ is unbounded on each unbounded subset of $D$. A closed copy of $\omega_1$ is 
a direction in any Type I space.
We study various aspects of directions and show some independence results.
A sample: There is an $\omega$-bounded Type I space without direction; 
{\bf PFA} implies that a locally compact countably tight $\omega_1$-compact Type I space contains a direction;
if there is a Suslin tree then there is an $\omega_1$-compact Type I manifold without direction;
there are Type I first countable spaces which contain directions and whose closed unbounded subsets contain each 
a closed unbounded discrete subset. We also study a naturel order on the directions of a given space
and show that we may obtain various classical ordered types with the space a manifold (often $\omega$-bounded).}

\tableofcontents

\section{Introduction}

In this paper `space' means `topological space'.
A space $X$ is
{\em Type I} if $X=\cup_{\alpha<\omega_1}X_\alpha$
with $X_\alpha$ open, $\wb{X_\alpha}$ Lindel\"of and $\wb{X_\alpha}\subset X_\beta$ whenever $\alpha<\beta$.
We will assume that none of our Type I spaces are Lindel\"of, so that $X\not= X_\alpha$ for any $\alpha$,
although this is not a part of the usual definition.
A subset of a Type I space $X$ is called {\em bounded} if it is contained in some $X_\alpha$, otherwise
it is called {\em unbounded} (this definition does not depend on the particular sequence 
$\langle X_\alpha\, : \,\alpha\in\omega_1\rangle$). Likewise, a function into a Type I space is unbounded iff its range is unbounded.
We abbreviate `closed and unbounded' by {\em club}.
We use brackets $\langle \,\,\, , \,\, \rangle$ for ordered pairs to avoid confusion with open intervals.
$\LL_{\ge 0}$ denotes the closed long ray (with the point $0$), i.e. the space $\omega_1\times [ 0,1)$ endowed with
the lexicographic order topology. We shall often view $\omega_1$ as a subspace of $\LL_{\ge 0}$, identifying
the point $\langle\alpha,0\rangle\in\LL_{\ge 0}$ with $\alpha\in\omega_1$, and writing directly $\alpha\in\LL_{\ge 0}$.
Every function is assumed continuous unless
specified, and the restriction of a function $f$ to a subset $A$ of its domain will be denoted by $f\restrict{A}$.

This paper is about the notion (introduced in \cite{mesziguessurf}) of directions in Type I topological spaces.
There are many equivalent definitions
(see Section \ref{secdir}), one relies on the following partial order.

\begin{defi}
 The unbounded function order (ufo) $\preceq$ on the subsets of a Type I space $X$ is defined
 as $A\preceq B$ whenever for all $f:X\to\LL_{\ge 0}$, $f\restrict{A}$ unbounded implies $f\restrict{B}$ unbounded,
 or equivalently if $f\restrict{B}$ bounded implies $f\restrict{A}$ bounded.
 If $A\preceq B\preceq A$, we say that $A$ and $B$ are ufo-equivalent.
 We denote the strict order by $\prec$. 
\end{defi}

The ufo is actually a preorder and 
extends the order given by the 
inclusion (and its variant `inclusion outside of a bounded set'). Notice that $X$ itself is always ufo-maximal.
While $\preceq$ is defined on all subsets of $X$, it is of 
interest only for unbounded subsets, since $f(A)$ is bounded whenever $A$ is. It is moreover easy to check that 
$A$ is ufo-equivalent to $\wb{A}$ (see Section \ref{sec:gendir}). We will
thus consider the restriction of $\preceq$ to the set of club subsets of $X$.

\begin{defi}\label{lemmadirpreceq}
  A club subset $D$ of a Type I space $X=\cup_{\alpha\in\omega_1}X_\alpha$ is a direction in $X$ iff it is 
  ufo-equivalent to any of its club subsets. If $D=X$, we say that $X$ is a direction in itself.
  We denote by 
  $\langle\mathfrak{D}_X,\preceq\rangle$ the poset of ufo-equivalent classes of directions in $X$.
\end{defi}

In words: $D$ is a direction in $X$ iff for any function 
$f:X\to\LL_{\ge 0}$, if an unbounded subset of $D$ has a bounded image then $f$ is bounded
on all of $D$.
The definition of direction is relative: the space $X$ consisting of two disjoint copies of  $\LL_{\ge 0}$ is not a direction in itself, 
but any copy of $X$ in
$Y=\LL_{\ge 0}\times\mathbb{S}^1$ is a direction in it since $Y$ is a direction in itself by
Propositions \ref{prop1} and 
\ref{lemmadoublestar} below. \\
The concept of direction grew out of while investigating the properties of club subsets of finite powers of $\LL_{\ge 0}$.
The following elementary proposition (due to the author and D. Cimasoni)
serves as a motivation for the questions we investigate in this paper.
The first octant $\OO$ is the subspace of $\LL_{\ge 0}^2$ given by $\{\langle x,y\rangle\in\LL_{\ge 0}^2\,:\,y\le x\}$.

\begin{prop}[\cite{mesziguessurf, meszigues+Dave}]
  \label{prop1} 
   \ \\
   (a) $\omega_1$ endowed with the order topology 
      is a direction in itself, and thus any closed copy of $\omega_1$ in a Type I space $X$ is a direction in $X$. 
      The same holds for $\LL_{\ge 0}$.\\
   (b) $\langle{\mathfrak{D}}_{\OO},\preceq\rangle$ is isomorphic to the ordered set 
      $\langle\{0,1\},\le\rangle$.
      Indeed, any horizontal line in $\OO$ is $\prec$ the diagonal, any two horizontal lines are
      ufo-equivalent, the diagonal is ufo-maximal in $\langle{\mathfrak{D}}_{\OO},\preceq\rangle$,
      and any direction in $\OO$ is ufo-equivalent to either an horizontal line or the diagonal.
      \\
   (c) If $f:\LL_{\ge 0}^2\to\LL_{\ge 0}$ is unbounded, then $f$ restricted to the diagonal is unbounded.
\end{prop}

Notice that a closed copy of $\omega_1$ in a Type I space has to be unbounded.
Point (b) is the reason of the choice of terminology: 
in $\OO$ the (ufo equivalent classes of) directions correspond to the two possibilities of escaping:
horizontally or diagonally.

In view of this proposition, the following (intentionally vaguely stated) questions felt natural:
\begin{enumerate}
 \item 
     Which properties do ensure that a Type I space possesses a direction~? 
 \item
     What does a direction `look like'~? Can it be discrete~? Can there be many discrete directions in a space~?
 \item 
     Which properties do ensure that a space contains a direction in itself~?
     How do directions in themselves behave in products and subspaces~?
 \item
     Given a Type I space $X$ containing many directions and an unbounded map $f:X\to\LL_{\ge 0}$, is there some direction on which $f$ is unbounded~?
 \item
     What does $\langle{\mathfrak{D}}_X,\preceq\rangle$ look like~? When does it have a maximal/minimal element (if non-empty)~? 
     Can it be `complicated'~?
 \item
     Do we gain or lose something if one replaces $\LL_{\ge 0}$ by another Type I space $Y$ and define $Y$-directions in a similar manner~?
\end{enumerate}

By `many directions' we mean the following:

\begin{defi}\label{defispanned}
  We call an open subset of a Type I space large if it contains a club set. A Type I space is spanned by directions 
  if any large open subset contains a direction.
\end{defi}

This paper contains our findings about these questions in various classes of spaces.
We have a bias towards manifolds, and we tried to find manifold examples of spaces with such or such property when possible. 
(By a manifold we mean a Hausdorff connected space each of whose points has a neighborhood 
homeomorphic to $\R^n$ for some $n$.
In a manifold with boundary the boundary points have neighborhoods homeomorphic to $\R^{n-1}\times\R_{\ge 0}$.)
Many of our results have the following format:

\newtheorem*{proto}{Theorem prototype}
\begin{proto}
  Let $\mathcal{P}$ be a property claiming
that a  Type I space contains some object related to directions.
Then the following holds.
\begin{itemize}
  \item[(A)]
   There is an $\omega$-bounded space without the property $\mathcal{P}$.
  \item[(B)]
  {\bf$\neg$SH} implies that there is an $\omega_1$-compact first countable space (often locally compact, often a manifold) 
   which does not have property $\mathcal{P}$.
  \item[(C)]
   {\bf PFA} implies that
   any $\omega_1$-compact countably tight locally compact space has property $\mathcal{P}$.
\end{itemize}
\end{proto}

{\bf PFA} denotes the proper forcing axiom and {\bf$\neg$SH} the negation of the Suslin hypothesis, thus the existence of a Suslin tree.
Recall that a space is {\em $\omega_1$-compact} iff it does not contain an uncountable closed discrete subset, and 
{\em $\omega$-bounded} iff the closure of any countable subset is compact. An $\omega$-bounded space is countably compact,
and a Type I countably compact space $X$ is $\omega$-bounded.
A space $X$ is {\em countably tight} if whenever $x$ is in the closure of some $A\subset X$,
then there is a countable $B\subset A$ having $x$ in its closure.

It happens that for properties related to Questions 2 to 5, the way of getting a counterexample in (A) and (B) 
is to start with a space in the relevant class which possesses no direction
and to feed it to a constructing procedure (adapted to the property) which gives a space 
that fails to satisfy $\mathcal{P}$ in the same class (or almost). 
Using a claim by Nyikos (whose proof has not been supplied to this date, see Section \ref{sec:omega_1cpct}), we can even 
add the item: 
\begin{itemize}
  \item[$\ast$(D)]
  {\bf$\diamondsuit^+$} implies that there is an $\omega$-bounded first countable space (often a manifold) which does not have property $\mathcal{P}$.
\end{itemize}
Since this point comes almost for free if Nyikos' claim is indeed correct, we decided to include it but denoted
each occurrence with an asterisk `$\ast$' until a proof of the claim is given.

The scheme (A)--(C) (and $\ast$(D)) is not completely accurate since some of the construction procedures
either push us outside the class of spaces we are looking at or simply do not work in some cases, but it summarizes a good 
part of this paper and explains the organization of its sections which we now detail.

\vskip .4cm
The generalities about Type I spaces and directions are gathered in Section \ref{sec:gen}.

Section \ref{secdir} is related to Question 1: 
We show that there is a Type I $\omega$-bounded space containing no direction
and obtain independence results in the format (B)--(C) above. We also show 
in Sections \ref{sec:closedpreimages} and \ref{sec:subsets2omega_1}   
that it is not possible to obtain $\omega$-bounded spaces containing no direction by looking at subspaces of the tree $2^{<\omega_1}$ 
(with a natural topology) or
closed preimages of $\omega_1$ with `small' fibers.

Section \ref{sec:constructions} contains the constructions, divided in two parts,
and is the longest section of the paper.
In \ref{secconstr} we build surfaces by gluing copies of the octant in many ways, which yield
various examples of ufo-equivalent classes of directions.
We show in particular that we can find surfaces $S$ such that $\langle{\mathfrak{D}}_S,\preceq\rangle$
is order-isomorphic to the following: the Cantor set, an antichain of size the continuum, 
any ordinal $<\omega_2$, a compact Aronszajn or Suslin line. 
Some of these constructions can be generalized, yielding the
construction procedures mentioned above. This is done in \ref{sec:generalization}.
We have chosen to gather all the constructions (and the main lemmas about them) in this section, but 
the reader may skip it until they are needed later.

In Section \ref{secgenprop3} we investigate Question 4: given a space spanned with directions 
and an unbounded map $f:X\to\LL_{\ge 0}$, is $f$ 
unbounded on some direction~? We obtain results of the (B)--(C) type,
though (A) is not achieved. We use 
the spaces $\we{X,s}$ of Section \ref{sec:generalization}.

Section \ref{secdirthem} deals with direction in themselves and Question 3.
We first show in \ref{sec:monolithic} that they have a close relationship with monoclinic and monolithic spaces (as defined by Nyikos),
especially for normal spaces. 
In \ref{sec:inproducts} we investigate what happens when we take the product of a direction in itself 
with some space. 
Then in \ref{sec:nodirinitself} we look inside the directions to see whether 
they always contain a direction in itself. Our result takes the (A)--(C) (and $\ast$(D)) prototypical form above.
We use the spaces $\Sigma(X,s)$ of Section \ref{sec:generalization}.

In Section \ref{sec:directiondiscrete} we analyse Question 2.
We produce an example of
a Type I locally metrizable locally separable space possessing directions and which has the property that any unbounded subset 
contains an uncountable closed discrete set. In particular, any direction contains a discrete one.
If $\diamondsuit^*$ holds, there is a collectionwise Hausdorff such space.
Other similar spaces are given as well, all using spaces $\Omega(X,A,s)$ of Section \ref{sec:generalization}.
On the other hand,
we show that a regular strongly collectionwise Hausdorff Type I space cannot contain a discrete direction.

In Section \ref{sec:ufo} we investigate Question 5, in particular the existence of a ufo-maximal or minimal direction
in a space.
Using spaces $S^{\alpha,h,r}$ and $\Xi(X,A,s)$ of Sections \ref{secconstr} and \ref{sec:generalization}
we obtain a theorem of the form (A)--(C) (and $\ast$(D)), that is,
the existence of spaces in various classes containing no ufo-minimal or maximal direction.
The corresponding item (C) for maximal directions is missing, though.
In \ref{sec:partialordermaps} we introduce a partial order on maps $X\to\LL_{\ge 0}$ and 
show some relations with the ufo.

We briefly investigate in Section \ref{sec:dominant} 
a strengthening of the concept of direction called {\em locally dominant subsets}, and show in {\bf ZFC}
(following ideas of P. Nyikos exposed in Section \ref{secconstr}) 
the existence of an $\omega$-bounded manifold containing no locally dominant subsets and similar results.

Finally, we look at Question 6 in Section \ref{sec:Hausdorff}.
It is easy to see that if the space $Y$ is regular, then a direction in some space 
$X$ is also a $Y$-direction. But this may fail
for Hausdorff spaces, and we shall provide simple counter-examples.

\vskip .4cm
All sections use the small toolbox of lemmas in Section \ref{sec:gen}. 
Other than that,
Sections \ref{secdir}, \ref{sec:constructions}, \ref{sec:monolithic}, \ref{sec:inproducts}, \ref{sec:partialordermaps} and \ref{sec:Hausdorff}
are (reasonably) self contained, while the others depend on Sections \ref{sec:betaomega}--\ref{sec:omega_1cpct} and \ref{sec:constructions}.

It might be important to note that all of our uses of {\bf PFA} are indirect, in the sense that 
we use some of its consequences obtained by various authors. Likewise, 
our use of $\diamondsuit^*$ is through a result of Devlin and Shelah. Actually, no 
purely `set-theoretic' construction is done in this paper, and we
shall mainly manipulate Type I surfaces and $\omega_1$-trees, building more complicated spaces either
by gluing simple parts together, or through inductive processes of length $\omega_1$ 
(or both). A good knowledge of $\omega_1$-trees and ordinal manipulation should be mostly sufficient to read this paper.

Before ending this introduction, a word about separation axioms.
{\em Every space is assumed to be Hausdorff}, but
most results assume regularity (which implies Hausdorff in our definition) or some weaker form of it. 
The reason is that regularity of a Type I space
$X$ ensures that it is Tychonoff, and
that there is {\em at least} one unbounded map $X\to\LL_{\ge 0}$, see Lemma \ref{lemma0}.
This is not true for Hausdorff spaces (as seen in Section \ref{sec:Hausdorff}), which might render the
definition of direction somewhat vacuous.

\vskip .4cm
{\em Disclaimer.}
This paper was mainly undertaken as a solo free-time project after our time at the university.
As a consequence, 
we might have a very partial knowledge of both the present state of the subject and
the literature, which yields a tendency to
over-cite our own work, and maybe to ignore other's
containing similar (or stronger) results.
We apologize if it is indeed the case, and would appreciate
any comments.

\vskip .4cm
{\em Acknowledgements.} This work started through discussions with P. Nyikos 
and S. Greenwood in Auckland in june 2006.
Some of the results presented here are either due to P. Nyikos (and quoted such) 
or were obtained during these discussions
(and might be misquoted).
We thank both of them as well as D. Gauld and the Auckland University for the kind invitation.
We also thank D. Cimasoni.



\section{Generalities}\label{sec:gen}

\subsection{Type I spaces}\label{secgenTypeI}

Recall that every space is assumed to be Hausdorff and every function continuous, otherwise specified.
We identify an ordinal with the set of its predecessors, 
$\omega$ is thus the set of natural numbers and $\omega_1$ the set of countable ordinals.

Given a Type I space $X=\cup_{\alpha\in\omega_1}X_\alpha$, 
the sequence $\langle X_\alpha\, : \,\alpha\in\omega_1\rangle$ is {\em canonical} iff $X_\alpha=\cup_{\beta<\alpha}X_\beta$ for all limit
$\alpha$. Any sequence can be made canonical by adding the missing $X_\alpha$'s.
When a Type I space $X$ is given, we shall implicitly assume that it comes with a canonical sequence $\langle X_\alpha\,:\,\alpha\in\omega_1\rangle$.
Two canonical sequences agree on a club subset of $\omega_1$.
The sets $B_\alpha(X)=\wb{X_\alpha} -  X_\alpha$ are called the {\em bones} of $X$, and the union of
the bones is the {\em skeleton} of $X$, denoted by $Sk(X)$. 
For notational convenience, we assume that $X_0$ is empty and define $B_0(X)$ to be any point in $X_1$.
(This makes the definition of a {\em slicer} below more straightforward.)
The next three lemmas are well known, and we will often use them (sometimes implicitly). 
We give them without proof.

\begin{lemma}\label{convenientlemma}
   Let $C$ be club in a Type I space $X=\cup_{\alpha <\omega_1}X_\alpha$.
   Then the following holds.\\
   (a) If $X$ is $\omega_1$-compact, then $C\cap B_\alpha(X)\not=\varnothing$ for a stationary set of $\alpha$.\\
   (b) If $X$ is $\omega$-bounded, then $C\cap B_\alpha(X)\not=\varnothing$ for a club set of $\alpha$.
\end{lemma}

\begin{lemma}\label{ccob}
  Let $X$ be a Type I space. Then $X$ is countably compact iff $X$ is $\omega$-bounded.
\end{lemma}

\begin{lemma}\label{lemma:newsmall}
   Let $X$ by Type I. Then $A\subset X$ is open (resp. closed) iff $A\cap X_\alpha$ is open 
   (resp. $A\cap\wb{X_\alpha}$ is closed) for each $\alpha$,
   and $f:X\to\LL_{\ge 0}$ is continuous iff $f\restrict{X_\alpha}$ is continuous for each $\alpha$.
\end{lemma}

It will be convenient to slice Type I spaces with a continuous function.

\begin{defi} \label{def:slicer}
  Let $X$ be Type I. A slicer is a map $s:X\to\LL_{\ge 0}$ such that 
  $s(\wb{X_{\alpha+1}} -  X_\alpha)\subset[\alpha,\alpha+1]$.
\end{defi}

Notice that $s(B_\alpha(X))=\{\alpha\}$ if $B_\alpha(X)$ is non-empty.

\begin{lemma}\label{lemma0}
  Let $X=\cup_{\alpha\in\omega_1}X_\alpha$ be a Type I regular space.
  Then $X$ is Tychonoff, each $\wb{X_\alpha}$ is normal,
  and there is a slicer $s:X\to\LL_{\ge 0}$.
\end{lemma}
\proof
  Recall that a regular Lindel\"of space is normal, so $\wb{X_{\alpha}}$ is normal, and $X$ is Tychonoff.
  By induction, suppose that $f$ is defined on $\wb{X_\beta}$  for all $\beta<\alpha$ such that 
  $f(B_\beta(X))=\{\beta\}$ if $B_\alpha(X)\not=\varnothing$. If $\alpha$ is limit define $f$ to be constant 
  of value $\alpha$ on $B_\alpha(X)$,
  otherwise define $f$ on $\wb{X_{\alpha+1}}$ that agrees with $f$ on $\wb{\cup_{\beta<\alpha}X_\beta}$
  and takes the value $\alpha+1$ on $B_{\alpha+1}(X)$ using the
  normality of $\wb{X_{\alpha+1}}$, Tietze's extension theorem, and the fact that $B_{\alpha+1}(X)\cap\wb{\cup_{\beta<\alpha}X_\beta}=\varnothing$.
\endproof

The next lemma will be useful in various places.

\begin{lemma}\label{lemmastar}
   Let $X$ be Type I and $\omega_1$-compact, and $f:X\to\LL_{\ge 0}$ be unbounded on a club $E\subset X$.
   Then the set 
   $$
      C(f,E)=\{x\in E\,:\, \text{ there is }\alpha\in\omega_1\text{ with }x\in B_\alpha(X)\cap f^{-1}([\alpha,\omega_1))\}
   $$ 
   is club in $X$.
\end{lemma}
\proof
   It is clear that $C(f,E)$ is closed, we show that it is unbounded.
   Let $\alpha_0\in\omega_1$, we define $\alpha_\gamma$ by induction.
   If $\gamma$ is limit, set $\alpha_\gamma=\sup_{\beta<\gamma}\alpha_\beta$.
   Given $\alpha_\gamma$, take some $x_{\gamma+1}\in (X - X_{\alpha_\gamma + 1})\cap E$ such that
   $f(x_{\gamma+1})\ge\alpha_\gamma$. Set $\alpha_{\gamma + 1}$ be minimal such that $x_{\gamma + 1}\in X_{\alpha_{\gamma + 1}}$.
   Set $F\subset E$ to be the closure of $\{x_{\gamma+1}\,:\,\gamma\in\omega_1\}$, by Lemma \ref{convenientlemma}
   $F$ intersects a stationary subset of bones with index $\alpha_\gamma$ such that $\gamma$ is limit.
   But for such $\gamma$, a point $x$ in $F\cap B_{\alpha_\gamma}(X)$ 
   is either equal to $x_{\gamma+1}$ or is an accumulation point of some
   $x_\beta$'s for smaller $\beta$'s converging to $\gamma$, and thus
   $f(x)\ge \alpha_\gamma$.
\endproof


\subsection{Directions: Equivalent definitions and basic properties}\label{sec:gendir}

\begin{lemma}\label{pulpita} 
  Let $X$ be Type I and
  $D\subset X$ be club. The following are equivalent:\\
  (a) $D$ is a direction in $X$,\\
  (b) for all maps $f:X\to\LL_{\ge 0}$, either $f$ is bounded on $D$ or for all $\alpha\in\omega_1$,  
  $f^{-1}([0,\alpha])\cap D$ is bounded,\\
  (c) there is no map $f:X\to \LL_{\ge 0}$ with $f\restrict{D}$ unbounded and 
  $f^{-1}(\{0\})\cap D$ unbounded,\\
  (d) for every $f:X\to\LL_{\ge 0}$ with $f\restrict{D}$ unbounded and all 
    $\alpha\in\omega_1$ there is $\gamma(\alpha)\in\omega_1$ such that 
    $f(D -  X_{\gamma(\alpha)})\subset[\alpha,\omega_1)$,\\
  (e) for every $f:X\to\LL_{\ge 0}$ with $f\restrict{D}$ unbounded there is a club $C\subset\omega_1$ such that 
    $f(D -  X_{\alpha})\subset[\alpha,\omega_1)$ for each $\alpha\in C$,\\
  (f) for every $f:X\to\LL_{\ge 0}$ and $\alpha\in\omega_1$, either $f^{-1}([0,\alpha])\cap D$ 
     is bounded, or $f^{-1}([\alpha,\omega_1))\cap D$ is bounded.
\end{lemma}
\proof

(a) is easily seen equivalent to (b) and to (d), while (b) $\rightarrow$ (c) and (e) $\rightarrow$ (d) are immediate.
\\
(c) $\rightarrow$ (b) 
If (b) does not hold for $D$, then there is a
$g:X\to\LL_{\ge 0}$ unbounded on $D$ and an $\alpha$ such that $g^{-1}([0,\alpha])\cap D$ is unbounded. Let $p:\LL_{\ge 0}\to\LL_{\ge 0}$
be continuous such that $p([0,\alpha])={0}$ and $p\restrict{[\alpha+1,\omega_1)}=id$, then $f=p\circ g$ yields a contradiction.
\\
(d) $\rightarrow$ (e) By a leapfrog argument and continuity of $f$. \\
(b) $\rightarrow$ (f) If $f$ is unbounded on $D$, then $f^{-1}([\alpha,\omega_1))\cap D$ is unbounded.\\
(f) $\rightarrow$ (d) If $f$ is unbounded on $D$, then $f^{-1}([0,\alpha])\cap D\subset X_{\gamma(\alpha)}$ for some $\gamma(\alpha)$. 
\endproof

An unbounded subset $A$ of a Type I space is a {\em predirection} if it satisfies the definition of direction except
the requirement that $A$ is closed.

\begin{lemma} $A$ is a predirection in $X$ iff $\wb{A}$ is a direction in $X$.
\end{lemma}
\proof
  Let $f:X\to\LL_{\ge 0}$, then $f$ is unbounded on $\wb{A}$ iff it is unbounded on $A$,
  $f^{-1}(\{0\})\cap A$ is bounded if $f^{-1}(\{0\})\cap\wb{A}$ is bounded, so $A$ is a predirection whenever 
  $\wb{A}$ is a direction. 
  If $A$ is a predirection and $f$ is unbounded on $A$ then $f^{-1}([0,1))\cap A$ is bounded,
  say by $\alpha$. If $x\in(\wb{A}-A)\cap f^{-1}(\{0\})$
  then $f^{-1}([0,1))$, being a neighborhood of $x$, has a non-empty intersection
  with $A\cap f^{-1}([0,1))$. Thus $f^{-1}(\{0\})\cap\wb{A}\subset\wb{X_\alpha}$,
  and $\wb{A}$ is a direction.
\endproof

It might be worth mentioning also the following simple lemma.

\begin{lemma}\label{insideomega_1}
  Let $A\subset\omega_1$ be endowed with the subspace topology. 
  Then $A$ is a direction in itself iff $A$ is stationary.
\end{lemma}
\proof
If $A$ is not stationary, it avoids a club subset of $\omega_1$, and it is easy to 
define a function witnessing that $A$ is not a direction in itself.
If $A$ is stationary and $f:A\to\LL_{\ge 0}$, then the closures of $f^{-1}(\{0\})$ and $f^{-1}([1,\omega_1))$ in $\omega_1$ cannot be both 
unbounded, otherwise their common intersection with $A$ is non-empty, which contradicts their definition. Thus if $f$ is unbounded,
$f^{-1}(\{0\})$ must be bounded.
\endproof



\section{Spaces with and without directions}\label{secdir}

\subsection{A type I $\omega$-bounded subspace of $\beta\omega$ containing no direction}\label{sec:betaomega}

Recall that $\beta\omega$ is the \v{C}ech-Stone compactification of $\omega$, and that $\omega^\ast = \beta\omega -\omega$. Both
are compact Hausdorff non-countably tight spaces.

\begin{thm}\label{thm:betaomega}
  There is a Type I $\omega$-bounded space $Q$ (which is a subspace of $\beta\omega$) containing no direction.
\end{thm}

We start with a lemma.

\begin{lemma}\label{lemma:betaomega}
  There is an $\omega$-bounded Type I subspace $X$ of $\omega^\ast$ such that $X_\alpha$ is Lindel\"of and $X_{\alpha+1}$ compact
  for each $\alpha$.
\end{lemma}
\proof
  Let $p\in\omega^\ast$. 
  It is shown in \cite[Corollary 1.5.4]{vanMill:1984} 
  that the cover $\mathcal{U}=\{U\subset\omega^\ast\,:\,U\text{ is clopen and }p\not\in U\}$ of $\omega^\ast -\{p\}$
  has no countable subfamily whose union is dense in $\omega^*-\{p\}$.  
  We define $X_\alpha$ by induction on $\alpha$, each will be an at most countable union of members of $\mathcal{U}$.
  Set $X_0=\varnothing$. 
  Given $X_\alpha$, take a finite family $\mathcal{F}\subset\mathcal{U}$ that covers the compact set $\wb{X_\alpha}$
  and some $U\in\mathcal{U}$ not included in $\wb{\cup\mathcal{F}}$. Set $X_{\alpha+1}=U\bigcup\cup\mathcal{F}$. 
  Then $X_{\alpha+1}$ is clopen (and thus compact). 
  For limit $\alpha$ define $X_\alpha=\cup_{\beta<\alpha}X_\beta$, then $X_\alpha$ is Lindel\"of and not dense in $\omega^\ast-\{p\}$.
  Then $X=\cup_{\alpha<\omega_1}X_\alpha$ is a Type I $\omega$-bounded space. 
\endproof

We will need the following classical result (see \cite[Theorem 1.5.2]{vanMill:1984} for the proof of a more general fact). 

\begin{lemma}
  Let $Y\subset\beta\omega$ be Lindel\"of. Then every function $g:Y\to[0,1]$ can be continuously extended 
  to a function $\beta\omega\to[0,1]$.
\end{lemma}

\proof[Proof of Theorem \ref{thm:betaomega}]
  Take $Q=X=\cup_{\alpha<\omega_1}X_\alpha$ given by Lemma \ref{lemma:betaomega}. 
  Let $C$ be club in $X$. Up to changing the canonical sequence, we may assume that 
  $C\cap(X_{\alpha+1}-X_\alpha)$ contains at least two points.
  We shall define closed subsets $D_\alpha,E_\alpha\subset C\cap \wb{X_\alpha}$ and $f_\alpha:\beta\omega\to[0,\alpha]\subset\LL_{\ge 0}$
  such that the following holds:
  \\
  (1) $f_\alpha(D_\alpha)=\{0\}$, $f_\alpha(E_\alpha)=\{0,1,2,\dots,\alpha\}\subset\LL_{\ge 0}$),
  \\
  (2) $f_\alpha\restrict{X_\beta}=f_\beta\restrict{X_\beta}$ whenever $\beta\le\alpha$,
  \\
  (3) $D_\alpha\cap X_\beta= D_\beta$ and $E_\alpha\cap X_\beta= E_\beta$ whenever $\beta\le\alpha$.\\
  Once these are defined, set $D=\cup_{\alpha<\omega_1}D_\alpha$, $E=\cup_{\alpha<\omega_1}E_\alpha$.
  Then the map $f:X\to\LL_{\ge 0}$ defined by $f(x)=f_{\alpha}(x)$ for $\alpha$ such that $x\in X_\alpha$
  is continuous, unbounded on $E$ and constant on $0$ on $D$. Since $D,E\subset C$, this shows that $C$ is not a direction in $X$.
  \\
  So, start with $E_0=F_0=\varnothing$. 
  Let now $\alpha$ be successor or $0$, then $X_\alpha=\wb{X_\alpha}$. 
  Take $d,e\in C\cap(X_{\alpha+1}-X_\alpha)$. The map $g:X_\alpha\cup \{d,e\}\to[0,\alpha+1]\subset\LL_{\ge 0}$ 
  defined as 
  $$
      g(x)=\left\{
          \begin{array}{cl}
            f_\alpha(x) &\text{if }x\in X_\alpha \\
            0 & \text{if }x=d \\ 
            \alpha+1 &  \text{if }x=e
          \end{array}\right.
  $$
  is continuous. Since $X_\alpha\cup \{d,e\}$ is compact,
  there is some $f_\alpha:\beta\omega\to[0,\alpha+1]$ extending $g$. Set $D_{\alpha+1}=D_\alpha\cup\{d\}$,
  $E_{\alpha+1}=E_\alpha\cup\{e\}$, by construction the above conditions are fulfilled.
  \\
  Let now $\alpha$ be limit. The map $g:X_\alpha\to[0,\alpha]$ defined by $g\restrict{X_\beta}=f_\beta\restrict{X_\beta}$ 
  for $\beta<\alpha$ is continuous. Since $X_\alpha$ is Lindel\"of, $g$ extends to a function
  $f_\alpha:\beta\omega\to[0,\alpha]$. Set 
  \begin{align*}
      D_\alpha&=(\cup_{\beta<\alpha}D_\beta)\cup(f_\alpha^{-1}(\{0\})\cap C \cap B_\alpha(X)), \\
      E_\alpha&=(\cup_{\beta<\alpha}E_\beta)\cup(f_\alpha^{-1}(\{\alpha\})\cap C \cap B_\alpha(X)).
  \end{align*}
  By construction $D_\alpha$ and $E_\alpha$ are closed in $X$ and satisfy the above conditions (1)--(3).
  This finishes the proof.
\endproof


\subsection{$\omega_1$-compact and $\omega$-bounded countably tight spaces}\label{sec:omega_1cpct}

The aim of this section is to prove the following pair of theorems. 

\begin{thm}[{\bf PFA}]\label{thmomega_1compact}
  Let $X$ be countably tight, locally compact, $\omega_1$-compact and Type I, and $C\subset X$ be club. 
  Then $X$ contains a closed copy of $\omega_1$ and thus a direction.
\end{thm}

\begin{thm}[$\neg\mathbf{SH}$]\label{surfaceomega_1compact}
  There is an $\omega_1$-compact Type I surface $S_T$ containing no direction.
\end{thm}

We thank P. Nyikos for the ideas leading to these theorems.
Theorem \ref{thmomega_1compact} is a simple consequence of two results by other authors. 

\begin{thm}\label{AxB} 
   (Balogh, Eisworth  and al.) 
   {\rm({\bf PFA})} Any countably tight perfect preimage of $\omega_1$ contains a closed copy of $\omega_1$.
\end{thm}

Recall that a map is perfect if it is closed and any point has a compact preimage.
Balogh showed in \cite{Balogh:1989} the
result for spaces of character $\le\aleph_1$ and
Eisworth proved it for countably tight spaces in \cite{Eisworth:2002}.
For first countable spaces, the conclusion holds in a model of
{\bf ZFC} + {\bf CH} without inaccessible cardinals 
by a result of Eisworth and Nyikos \cite{EisworthNyikos:2005}.

\begin{lemma}\label{diromegabounded}
  {\rm({\bf PFA})}
  Let $X$ be Type I, countably tight and $\omega$-bounded. Then $X$ contains a closed copy of $\omega_1$.
\end{lemma}
\proof
  $X$ is $\omega$-bounded and thus countably compact,
  so the bones $B_\alpha(X)$ are compact, and $Sk(X)$ is a perfect preimage of $\omega_1$. 
\endproof

\proof[Proof of Theorem \ref{thmomega_1compact}]
The First Trichotomy Theorem of Eisworth and Nyikos \cite{EisworthNyikos} shows that
{\bf PFA} implies that 
a locally compact space is either
a countable union of $\omega$-bounded subspaces, 
or has a closed uncountable discrete space, or has a Lindel\"of subset with a non-Lindel\"of closure.
The latter two are impossible if $X$ is Type I and $\omega_1$-compact, thus
$X$ is a countable union of $\omega$-bounded subspaces, one
of which must be unbounded.
We may thus assume that $X$ is $\omega$-bounded and conclude with Lemma \ref{diromegabounded}.
\endproof

\begin{prob}\label{prob:directionsPFA}
   Is there a model of {\bf ZFC} where any countably tight (or first countable) $\omega$-bounded Type I
   space contains a direction but not necessarily a closed copy of
   $\omega_1$?
\end{prob}

We now prove Theorem \ref{surfaceomega_1compact} in a sequence of easy lemmas.
We refer to \cite{Kunen} for the definitions and basic properties of 
Aronszajn and Suslin trees.

\begin{lemma}\label{lemmaunpeuinutile}
  Let $T$ be an Aronszajn tree, and $A\subset T$ be unbounded. Then 
  there are incomparable $x,y\in T$ with $S\cap\{z\in T\,:\, z\ge x\}$ and
  $A\cap\{z\in T\,:\, z\ge y\}$ both unbounded.
\end{lemma}
\proof
  Suppose that there do not exist $x$ and $y$ as in the statement.
  Therefore, letting $E$ be the set of $y\in T$ such that
  there are uncountably many points of $S$ above it, $E$ intersects each level of 
  $T$ in exactly one point. But $E$ is then a chain of type $\omega_1$, a contradiction.
\endproof

If $X$ is a Type I manifold, 
recall from \cite{Nyikos:1984} that $\Upsilon(X)$ is the tree defined by the non-metrizable component 
boundaries of $X$, i.e.
$\sigma\in\Upsilon(X)$ if it is the boundary in $X$ of a non metrizable component of $X - X_\alpha$   (for some $\alpha$),
$\sigma\le\sigma'$ whenever $\sigma,\sigma'$ bound the components $S,S'$, respectively, and $S'\subset S$. 

\begin{lemma}\label{propAron}
  Let $M$ be a Type I manifold such that $\Upsilon(M)$ is an Aronszajn tree. Then, $M$ contains no direction.
\end{lemma}
\proof
  Let $D$ be club, we show that there are $D_1,D_2$ disjoint club subsets of $D$ and a map 
  $f:M\to\LL_{\ge 0}$ with $f\restrict{D_1}$ bounded and $f\restrict{D_2}$ unbounded. 
  Let 
  $\Upsilon(D)\subset\Upsilon(M)$
  be the subtree of boundaries $\partial C$ such that $C$ intersects $D$. Since $D$ is club, 
  $\Upsilon(D)$ is unbounded, so there must be two non metrizable components $C_1,C_2$ of $M$ 
  such that $D_1=D\cap C_1$ and $D_2=D\cap C_2$ are unbounded
  by Lemma \ref{lemmaunpeuinutile}. Let $s:M\to\LL_{\ge 0}$ be the slicer given by Lemma \ref{lemma0}. Let $U$ be a neighborhood of $\partial C_1$ in $C_1$.
  Define $g:M\to\LL_{\ge 0}$ such that 
  $g\restrict{M -  C_1}=s\restrict{M -  C_1}$, 
  $g\restrict{C_1 -  U}\equiv 0$, 
  using the fact that $\partial C_1,\partial C_2$ are disjoint and 
  included in the normal subspace $\wb{M_\alpha}$ for some $\alpha$. 
\endproof

\begin{rem}\label{rem:aron}
The previous lemma also holds under the weaker following assumptions: $X$ is regular and 
for each $\alpha\in\omega_1$, $X-\wb{X_\alpha}=\cup_{j\in J(\alpha)}U_{\alpha,j}$, where 
$\{U_{\alpha,j}\,:\,j\in J(\alpha)\}$ is a family of open unbounded subsets which are $2$-by-$2$ disjoints,
such that $\{U_{\alpha,j}\,:\,\alpha\in\omega_1,\,j\in J(\alpha)\}$ forms an Aronszajn tree for the reverse inclusion $\supseteq$.
\end{rem}

\begin{figure}[h]
  \begin{center}
  \epsfig{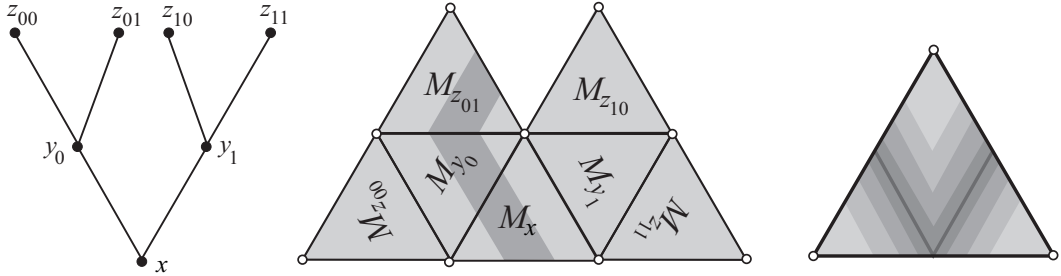}
  \caption{The surface $S_T$ and the wedges $W_{n,x}$.}\label{figureTsec3}
  \end{center}
\end{figure}

It is easy to build a surface $S_T$ having a prescribed Hausdorff $\omega_1$-tree $T$ as $\Upsilon$-tree. There are many constructions, let us 
quickly describe one. 
The construction is easier to describe when the tree is binary, we shall do this case since it is enough for our purposes.
First, build the road space
$R_{T}$ by inserting a line segment $(0,1)$ between each consecutive points in $T$. We then `fatten' $R_T$ to obtain $S_T$, 
which is the union of equilateral triangles (with deleted vertices) $M_x$ for $x\in T$ that are
glued together according to $T$, as in Figure \ref{figureTsec3} (left, middle), where $y_0,y_1$ are the immediate successors of $x\in T$, and similarly
for $z_{i0},z_{i1}$ and $y_i$ ($i=0,1$). Inside the triangles, the topology is the usual one, and neighborhoods of points
in the boundaries at successor heights are clear. A neighborhood of a point in the bottom boundary of $M_x$ for $x$ at a limit height
is a strip zigzagging the surface down to a smaller height, as the one partly depicted in darker grey in Figure \ref{figureTsec3} (middle).
Notice that each $M_x$ is a countable union of compact larger and larger `wedges' $W_{n,x}$ ($n\in\omega$) shown in Figure \ref{figureTsec3} (right). 
An uncountable discrete subset $E$ of $S_T$ has an uncountable subset that intersects $\cup_{x\in T}W_{n,x}$ for some fixed $n$. If there
is a sequence $x_k\in T$ converging to $x\in T$ such that $E\cap W_{n,x_k}\not=\varnothing$, then there is an accumulation
point in the bottom boundary of $E\cap W_{n,x}$. It follows that $S_T$ is $\omega_1$-compact iff $T$ is $\omega_1$-compact.
For more details, see \cite{meszigues+Nyikos} where a variant of this construction is given.

\proof[Proof of Theorem \ref{surfaceomega_1compact}]
  Start with a Hausdorff binary Suslin tree.
  A closed discrete subset of a Hausdorff tree is a countable union of antichains 
  (see Theorem 4.11 in \cite{Nyikos:trees}), so a Hausdorff Suslin tree is $\omega_1$-compact.
  Now $S_T$ described above is $\omega_1$-compact and
  satisfies $\Upsilon(S_T)=T$.
  We conclude via Lemma \ref{propAron}.
\endproof

Notice that Lemma \ref{propAron} shows that it is easy to obtain Type I surfaces without direction in {\bf ZFC}. 
Without additional axioms, the surface cannot be shown to be $\omega_1$-compact.
This construction cannot be generalized to obtain an $\omega$-bounded space, since $\omega$-bounded trees contain a copy of $\omega_1$.
On the other hand, Nyikos has announced the following theorem some years ago, which would give a consistent example
of an $\omega$-bounded Type I surface without direction,
but his construction is still unpublished.
Since some of our results rely on it, we quote it with a `$\ast$',
and stress that all our results depending on it do hold {\em only under the proviso that Nyikos construction is indeed correct}.

\begin{thm*}\label{AxN}
   (Nyikos)
   {\rm($\diamondsuit^+$)} There is an $\omega$-bounded surface $N$ containing no direction.
\end{thm*}

The next two sections show that such an example cannot be too simple.


\subsection{Closed preimages of $\omega_1$}\label{sec:closedpreimages}

A closed preimage of $\omega_1$ is a space $X$ such that there is
a closed map $\pi:X\to\omega_1$, with $\pi^{-1}(\{\alpha\})\not=\varnothing$
for all $\alpha$. Under some weak axioms,
it is not hard to find such spaces with very small fibers that do not contain a closed copy of $\omega_1$, 
for instance, it is consistent with {\bf MA$+\neg$CH} that there are 
closed $2$-to-$1$ preimages of $\omega_1$ (i.e. $|\pi^{-1}(\{\alpha\})|=2$ for all
$\alpha$) that do not contain a copy of $\omega_1$ (this result of Nyikos
is published in \cite[19.1]{FremlinPerfect}, see also \cite{Nyikoscoherent}). 
This cannot be shown in {\bf ZFC} by Theorem \ref{AxB}.
We shall show here that the fibers must be bigger for a space to contain no direction.
Let $\langle \mathfrak{C}_X,\subseteq\rangle$ denote the poset of club subsets of  
a Type I space $X$, ordered by the inclusion. 
The following lemma is well known, 
(see \cite{EisworthNyikos:2005}, Prop. 2.2, for instance).

\begin{lemma}\label{cpomegaclosed} 
   Let $X$ be a closed preimage of $\omega_1$. 
   Then the following holds.\\
   (a) 
   $\langle \mathfrak{C}_X,\subseteq\rangle$ is downward countably closed, i.e.
   if $\{C_i\}_{i\in\omega}$ is a sequence of clubs of $X$ with 
   $C_i\supseteq C_{i+1}$ for all $i$, then there is a club $C\subset X$
   with $C\subset C_i\,\forall i\in\omega$.\\
   (b)
   If $C\subset X$ is club, $\{\alpha\, :\, C\cap\pi^{-1}(\{\alpha\})\not=\varnothing\}$
   is club in $\omega_1$.
\end{lemma}

\begin{rem}
  Lemma \ref{cpomegaclosed} (a) also holds for countably compact Type I spaces.
\end{rem}

\begin{thm}\label{thmcpomega1} 
  Let $X$ be a closed preimage of $\omega_1$ and $S\subset\omega_1$ be stationary.\\
  (a) If $|\pi^{-1}(\{\alpha\})|<\omega$ 
  for $\alpha\in S$,
  then $X$ contains a direction.\\
  (b) Suppose that
  $2^\omega\ge\lambda >\omega_1$, with $\lambda$ regular, and that
  $|\pi^{-1}(\{\alpha\})| < \lambda$ for $\alpha\in S\subset\omega_1$. Then $X$ contains a direction.
\end{thm}

This theorem is an easy consequence of the following lemma. It is convenient to introduce some notation.
Recall that $2^{<\omega_1}$ is the set of sequences of $0$ and $1$ of length $<\omega_1$ ordered by the inclusion, which makes it a tree.
We say that a Type I space $X$ is {\em $2^{<\omega_1}$-spiked} if there are 
clubs $C_t\subset X$ for $t\in 2^{<\omega_1}$ such that
$C_t\cap C_s\not=\varnothing$ iff $s$ and $t$ are non-comparable for the inclusion,
and $C_t\subset C_s$ whenever $t\supset s$.

\begin{lemma}\label{lemmanodir}
  Suppose that $X$ is Type I with $\langle \mathfrak{C}_X,\subseteq\rangle$ countably closed.
  If $X$ does not contain a direction, 
  then any club in $X$ is $2^{<\omega_1}$-spiked.
\end{lemma}
\proof
  We proceed by induction. Let $C$ be a club in $X$ and set
  $C_\varnothing=C$.
  If $t\in 2^{<\omega_1}$ is at a limit level, set $C_t=\cap_{s\subset t}C_s$, which is club by assumption. 
  Given $C_t$ for $t\in 2^{<\omega}$, choose 
  $f_t:X\to\LL_{\ge 0}$ such that $f^{-1}(\{0\})\cap C_t$ is club
  and $f\restrict{C_t}$ unbounded. Set
  $C_{t\smallfrown 0}=f^{-1}(\{0\})\cap C_t$ and $C_{t\smallfrown 1}=f^{-1}([1,\omega_1))\cap C_t$.
\endproof

\proof[Proof of Theorem \ref{thmcpomega1}]
  Let $X$ be a closed preimage of $\omega_1$ containing no direction.
  By Lemmas \ref{cpomegaclosed} and \ref{lemmanodir} there is a collection $C_t$ of pairwise disjoint club subsets of cardinality the continuum.\\
  (a) Take a countably infinite subset $\{C_i\,:\,i\in\omega\}$ of the $C_t$, then there is some $\alpha\in S\cap\bigcap_{i\in\omega}\pi(C_i)$.
  Then $\pi^{-1}(\{\alpha\})$ contains at least $\omega$ points, a contradiction.\\
  (b)
  Let $K_t\subset\omega_1$ be the set of $\alpha$ such that 
  $C_t\cap\pi^{-1}(\{\alpha\})\not=\varnothing$. By Lemma
  \ref{cpomegaclosed} (b), $K_t$ is club, let thus $\alpha_t\in S\cap K_t$.
  Since $\lambda>\omega_1$ is regular, there is some $\alpha\in S$ such that
  $|\{t\in 2^\omega\,:\,\alpha_t=\alpha\}|=\lambda >\omega_1$.
  Thus, there are $\lambda$ disjoint sets intersecting 
  $\pi^{-1}(\{\alpha\})$, which contradicts $|\pi^{-1}(\{\alpha\})|<\lambda$.
\endproof

\begin{cor}
  If $\text{\rm cf}(2^\omega)>\omega_1$, any closed preimage of $\omega_1$ without direction
  must satisfy $|\pi^{-1}(\{\alpha\})|\ge 2^\omega$ for club many $\alpha$.
\end{cor}
\proof
  Otherwise $|\pi^{-1}(\{\alpha\})|< 2^\omega$ for stationary many $\alpha$.
  Taking $\kappa$ to be the supremum of $|\pi^{-1}(\{\alpha\})|$ for these $\alpha$,
  we have $\kappa < 2^\omega$, and the result follows.
\endproof


\subsection{Subsets of $2^{<\omega_1}$}\label{sec:subsets2omega_1}

The aim of this section is to show that there is no hope of finding an $\omega$-bounded Type I space without
direction by looking at subsets of the tree $2^{<\omega_1}$
endowed with the `accumulated product topology', i.e.
a neighborhood base for $t\in 2^\alpha$ is given by the sets
of the form
$\cup_{\beta<\gamma\le\alpha} B(\gamma ;\alpha_1,\dots,\alpha_n)$
where $\alpha_1 <\dots <\alpha_{n-1}<\alpha_n\le\beta<\alpha$ and
  $$
    B(\gamma ;\alpha_1,\dots,\alpha_n)=\{s\in 2^\gamma\,:\, s(\alpha_i)=t(\alpha_i)
    \text{ for }i=1,\dots,n\}.
  $$
Notice in passing that $2^\alpha$ is homeomorphic to the Cantor set in $2^{<\omega_1}$,
which is $\omega$-bounded, 
locally metrizable, $2^{<\omega_1}$-spiked (set 
$C_t=\{s\in X\,:\, s\ge t\}$), but contains a lot of directions.
We could hope that one of its unbounded
subsets has a closure that does not contain any direction. Not surprisingly this is not the case.

\begin{prop}\label{prop2<omega1} 
  Let $2^{<\omega_1}$ be endowed with the above topology. The following holds.\\
  (a) Let 
  $D$ be a direction in $2^{<\omega_1}$.
  Then $D$ is contained in a branch of $2^{<\omega_1}$
  outside of a compact subset.\\
  (b) Let 
  $A\subset 2^{<\omega_1}$ be unbounded. Then $\wb{A}$ contains a closed copy of
  $\omega_1$.
\end{prop}

\proof \ \\
  (a) 
  As in Lemma \ref{propAron}, if there are non-comparable $x,y\in 2^{<\omega_1}$ with
  $\{z\in D\,:\, z\ge x\}$ and $\{z\in D\,:\, z\ge x\}$ both unbounded then $D$ is not a direction.
  Therefore, for each $\alpha\in\omega_1$ there is a unique $x\in 2^\alpha\subset  2^{<\omega_1}$ with
  $\{z\in D\,:\, z\ge x\}$ unbounded.
  \\
  (b)
  Suppose that there is some unbounded $A\subset 2^{<\omega_1}$ with
  $\wb{A}$ containing
  no copy of $\omega_1$. 
  We can assume that $A\cap 2^\alpha =\{x_\alpha\}$ for each
  $\alpha\in\omega_1$. Indeed, 
  $\wb{A}$ must intersect $2^\alpha$ for a club $C=\{\alpha_\beta\,:\,\beta\in\omega_1\}$, 
  we take only one point in $A$ (or in $\wb{A}$ if there is none 
  available in $A$) in $2^{\alpha_\beta}$, and we can assume that $\alpha_0=0$.
  Any $\gamma\in\omega_1 -  C$ belongs to a unique interval
  $(\alpha_{\beta(\gamma)},\alpha_{\beta(\gamma)+1})$, set $x_{\gamma}= x_{\alpha_{\beta(\gamma)+1}\restrict{\gamma}}$.
  Then, setting $B=\{x_\beta\,:\,\beta\in\omega_1\}$, $y\in\wb{B}\cap 2^\alpha$ 
  if and only if
  $y=z\restrict{\alpha}$ for some $z\in\wb{\{x_\beta\,:\,\beta\in C\}}$.
  \\
  Let us thus $A$ be $\{x_\alpha\,:\,\alpha\in\omega_1\}$, with
  $x_\alpha\in 2^\alpha$ and $\wb{A}$ containing
  no copy of $\omega_1$.
  Then, $\wb{A}$ is bounded in every branch, i.e.:
  $\forall y\in 2^{\omega_1}$, $\exists \alpha(y)\in\omega_1$ such that
  $\forall \beta\ge\alpha(y)$, $y\restrict{\beta}\not\in\wb{A}$.
  This means that $\forall \beta\ge\alpha(y)$, we have:
  \begin{align*}
      \bullet & x_\beta\not= y\restrict{\beta}, \text{ if $\beta$ is successor,} \\
     \bullet & \exists \gamma(y,\beta)<\beta\text{ and a finite }F(y,\beta)\subset\gamma(y,\beta)\\
         \   &\text{such that } \forall\xi\text{ with } 
              \gamma(y,\beta)\le\xi<\beta,\,
         x_\xi\restrict{F(y,\beta)}\not= y\restrict{F(y,\beta)},\\
         \ & \text{ if $\beta$ is limit.}                                     
  \end{align*}
  Let us fix $y\in 2^{\omega_1}$. Since $\gamma(y,\beta)<\beta$ if $\beta$ is limit, 
  there is a stationary $S_0(y)\subset\omega_1$ such that
  $\gamma(y,\beta)$ is constant on $S_0(y)$, taking the value $\gamma(y)$.
  Since there are only countably many finite subsets of $\gamma(y)$, 
  there is a stationary $S(y)\subset S_0(y)$ such that $F(y,\beta)$
  takes the constant value $F(y)$ on $S(y)$. Thus, whenever $\beta\in S(y)$, given
  $\xi$ with $\gamma(y)\le\xi<\beta$,  we have that
  $x_\xi\restrict{F(y)}\not= y\restrict{F(y)}$.
  But since $S(y)$ is unbounded, this is true for all $\xi\ge\gamma(y)$.
  \\
  Now, endow $2^{\omega_1}$ with the product topology, and cover it with
  the open sets 
  $$
    U_y=\{z\in 2^{\omega_1}\,:\,z\restrict{F(y)}= y\restrict{F(y)}\}.
  $$
  By compactness there are $y_0,\dots,y_n$ such that $\{U_{y_i}\,:\,i=0,\dots,n\}$
  cover $2^{\omega_1}$.
  Letting $\gamma$ be the max of the $\gamma(y_i)$, this implies that
  for all $\xi\ge\gamma$, for all $i=1,\dots,n$,
  $x_\gamma\restrict{F(y_i)}\not=y_i\restrict{F(y_i)}$. But since for any $z$ in
  $2^{\omega_1}$ there is an $i$ with $y_i\restrict{F(y_i)}=z\restrict{F(y_i)}$,
  the $x_\xi$ for $\xi\ge\gamma$ take values different from any $z\in 2^{\omega_1}$
  on $\gamma$, which
  is a contradiction.
  \\
  Thus, there is some branch $y\in 2^{\omega_1}$ such that $\wb{A}$ intersects
  $y$ on a club set, so $\wb{A}$ contains a copy of $\omega_1$.
\endproof





\section{Constructions}\label{sec:constructions}

\subsection{Surfaces built with the first octant $\OO$ and their ufo-equivalent classes of directions}\label{secconstr}

Here we describe some constructions that yield different examples of $\langle{\mathfrak{D}}_X,\preceq\rangle$. 
Some already appeared in \cite{Nyikos:1984} and \cite[Appendix A]{mesziguessurf}.
All are surfaces (with or without boundary), most of them $\omega$-bounded.
Our goal is to show that there are `fairly complicated' examples even in this restricted case
and to lay the ground for the more elaborate constructions of Section \ref{sec:generalization}.

\begin{figure}[h]
  \begin{center}
  \epsfig{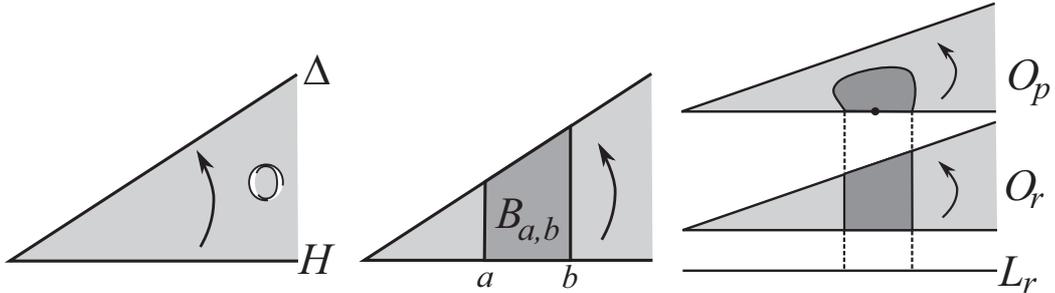}
  \caption{The octant (left), neighborhoods in Example \ref{ex1} (middle, right)}\label{figure1sec5}
  \end{center}
\end{figure}

We will use the first octant octant $\mathbb{O}$, which we recall is
$\{\langle x,y\rangle\in\LL_{\ge 0}^2\, : \,y\le x\}$, as a building brick.
Let $H=\{\langle x,0\rangle \,:\, x\in\LL_{\ge 0}\}\subset\mathbb{O}$ and $\Delta=\{\langle x,x\rangle \,:\, x\in\LL_{\ge 0}\}\subset\mathbb{O}$
be respectively the horizontal and the diagonal sides of $\mathbb{O}$. 
(Strictly speaking, $\OO$ has only one boundary component homeomorphic to the long line $\LL$, so
we call $H$ and $\Delta$ `sides' of $\OO$.)
We picture
$\OO$ as on Figure \ref{figure1sec5} (left), with an arrow going from $H$ to $\Delta$ when the orientation is relevant.

Recall Proposition \ref{prop1} (b): any horizontal in $\mathbb{O}$ is ufo-equivalent to $H$
and

\begin{equation}\label{O}
   \langle {\mathfrak{D}}_{\mathbb{O}},\preceq\rangle \simeq \langle\{H,\Delta\},    H\prec\Delta\rangle
   \simeq \langle\{0,1\},\le\rangle.
\end{equation}

For convenience, we compile in the next lemma some facts that we shall often use (sometimes implicitly). 
We omit the proof since we will prove more general statements later.

\begin{lemma}\label{lemmaclassical}
   Let $U\subset \OO$ be open. 
   If $U$ contains a club subset of some horizontal line in $\OO$, it 
   contains some strip $[\alpha,\omega_1)\times [a,b]$.
   If $U$ contains a club subset of $\Delta\subset\OO$, it contains 
   some $[\alpha,\omega_1)^2\cap\OO$.
   If $C\subset\OO$ is club, then it club-intersects $\Delta$ or some horizontal line.
\end{lemma}


\subsubsection{Obtaining the Cantor set}

\begin{ex}\label{ex1}
  An $\omega$-bounded surface $S$ with
  $\langle {\mathfrak{D}}_{S},\preceq\rangle=\langle K,\le\rangle$, the
  Cantor set embedded in the interval with the induced order. 
\end{ex}

In particular,
$\langle {\mathfrak{D}}_{S},\preceq\rangle$ has cardinality the continuum.
The idea is to start with $[0,1]\times\LL_{\ge 0}$ and to replace each $\{p\}\times\LL_{\ge 0}$ for 
$p$ dyadic with a copy of $\mathbb{O}$. Denote the dyadic numbers by $\mathcal{D}$. Let thus 
$$ 
 S=\left(\bigcup_{p\in[0,1]\cap\mathcal{D}}O_p\right)\cup
 \left(\bigcup_{q\in[0,1] - \mathcal{D}}L_q\right),
$$
where each $O_p$ is a copy of $\mathbb{O}$ and $L_q$ a copy of $\LL_{\ge 0}$.
The neighborhoods of a point in the interior of some $O_p$ are the usual neighborhoods.
If $E$ is a subset or a point of $\mathbb{O}$, we write $E_p$ for the corresponding subset (or point) of $O_p$,
and similarly for subsets of $\LL_{\ge 0}$. For $a,b\in\LL_{\ge 0}$, $a<b$, define
$B_{a,b}=\{\langle x,y\rangle \in\mathbb{O}\,:\,a<x<b\}$ (Figure \ref{figure1sec5}, middle).
If $z=\langle x,0\rangle_p$, a system of neighborhoods of $z$ is given 
by the $N_{p',a,b}$ with $p'<p$, $a<x<b$, where
$N_{p',a,b}$ is the union of the $(B_{a,b})_r\subset O_r$ and of $(a,b)_r\subset L_r$ for all $p'<r<p$
together with an open set in $O_p$ whose intersection with $H_p$ is the segment $(a,b)$,
see Figure \ref{figure1sec5} (right).
If $z=\langle x,x\rangle_p$, we define the neighboroods similarly, 
taking $p<r<p'$ and an open set intersecting $\Delta_p$ as $(B_{a,b})_p$ does.
Finally, if $x\in L_q$, take $q'<q<q''$, and define the neighborhoods in the same way again.
It is then not difficult to check that $S$ is an $\omega$-bounded surface.
Of course, each $L_q$ and each $H_p$ and $\Delta_p$ is a direction in $S$.
\begin{lemma}\label{lemmachiant0}
  A club subset of $S$ must club-intersect one of 
  the $L_q$'s or one of the $O_p$'s. 
  In the latter case it intersects either $\Delta_p$ or a horizontal line in $O_p$.
\end{lemma}
\proof
  Let $k:S\to[0,1]\times\LL_{\ge 0}$ be the map that collapses each $O_p$ by sending $\langle x,y\rangle$ to $x$.
  Since $k$ is perfect and $S$ is locally compact, $k$ is closed, so $k(C)$ is club in $[0,1]\times\LL_{\ge 0}$.
  The result follows from Lemma \ref{lemmaclassical}.
\endproof

It follows that a direction in $S$ must be ufo-equivalent
to $L_r$, $H_r$ or $\Delta_r$ for some $r$.
For $p\in[0,1]\cap\mathcal{D},$ define the map $f_p:S\to\LL_{\ge 0}$ by 
$f_p\restrict{L_r}=f_p\restrict{O_r}\equiv 0$ for $r<p$, 
$f_p(\langle x,y\rangle_p)=y$, $f_p(\langle x,y\rangle_{q})=x$ if $q>p$ with $q\in\mathcal{D}$,
and $f_p\restrict{L_q}=id$ if $q>p$ with $q\not\in\mathcal{D}$. 
Then $f_p$ is continuous and witnesses that
$$
  \left. 
  \begin{array}{c} \Delta_r \\ H_r \\ L_r \end{array}
  \right\}
  \quad \prec \quad H_p\quad\prec\quad\Delta_p\quad
  \left\{ 
  \begin{array}{c} \Delta_{r'}  \\ H_{r'} \\ L_{r'} \end{array}
  \right.
$$ 
for all $r<p<r'$.
It follows that $\langle {\mathfrak{D}}_{S},\preceq\rangle$
is indeed the Cantor Set.

\begin{ex}\label{ex2}
   An $\omega$-bounded surface $S$ with 
   $\langle {\mathfrak{D}}_{S},\preceq\rangle$ containing antichains of cardinality
   the continuum.
\end{ex}

As in Example \ref{ex1}, except
that one replaces the $\{p\}\times\LL_{\ge 0}$ for $p$ dyadic by a copy $L^2_p$ of $\LL_{\ge 0}^2$ instead of $\mathbb{O}$. 
The  
antichain is given by the diagonals of each inserted
$L^2_p$ and of the $\{q\}\times\LL_{\ge 0}$ for $q$ non-dyadic.
Note however that the diagonal in $L^2_p$ is $\succ$ the horizontal,
so some members of $\langle{\mathfrak{D}}_{S},\preceq\rangle$ are comparable.

One can also replace each $\{p\}\times\LL_{\ge 0}$ for $p$ dyadic by a copy of $S'$, $S'$ being made
of two copies of the surface of Example \ref{ex1} with their top boundaries identified. Then,
$\langle {\mathfrak{D}}_{S},\preceq\rangle$ contains copies of the Cantor set $K$ and antichains
of cardinality the continuum.

The next examples are slightly more complicated.


\subsubsection{Obtaining an antichain of cardinality the continuum (and additional properties)}\label{sec:exDominant}

\begin{ex}[P. Nyikos \cite{NyikosPC}]\label{exDominant}
   An $\omega$-bounded surface $R$ with 
   $\langle \mathfrak{D}_{R},\preceq\rangle$ an antichain of cardinality
   the continuum such that if $D,D'$ are ufo-equivalent directions in $R$, then neither $D-D'$ nor $D'-D$ contain a club.
\end{ex}
P. Nyikos used this example to prove Theorem \ref{thmnld} 
in Section \ref{sec:dominant} below. We give a variant of his construction, which
is a `stratisfied' version of Example \ref{ex2}.
We built $R=R^{\omega_1}$ by induction, defining $R^{\alpha}$ for each $\alpha\in\omega_1$. Each is an $\omega$-bounded surface with the
following properties:\\
(i) Whenever $\beta \le \alpha$, there is a closed map $k_{\beta,\alpha}:R^{\beta}\to R^{\alpha}$ that is bone-preserving. \\
(ii) Each $R^\alpha$ contains so-called {\em rays} (i.e. copies of $\LL_{\ge 0}$), some of which are added at stage $\alpha$,
and some rays will be affected during the process. 
Each ray crosses each bone of $S^\alpha$ in exactly one point,
but two rays may intersect on a bounded set.\\
(iii)
The unaffected rays  are dense in $R^{\alpha}$.
Moreover, for $\beta>\alpha+1$ 
the points that are crossed only by unaffected rays added at stage $\alpha$ are dense
in the $\beta$-th bone of $S^{\alpha+1}$,
and for limit $\alpha$ the same holds with unaffected rays added at stage $<\alpha$.

Start with $R^0=\LL_{\ge 0}\times\mathbb{S}^1$,
the rays in $R^0$ being $\LL_{\ge 0}\times\{x\}$ for each $x\in \mathbb{S}^1$, and none is affected.
At stage $\alpha+1$, 
take a 
countable set $\mathcal{D}$ of unaffected rays whose intersection 
with the $\alpha+1$-th bone of $S^\alpha$ is dense, and such that any point in this intersection is crossed only by 
unaffected rays that were inserted at the $\alpha$-th stage.
Above height $\alpha+1$, replace each 
ray in $\mathcal{D}$ by a copy of $[\alpha+1,\omega_1)^2\subset\LL^2$ similarly as in Example \ref{ex2}, as suggested in 
Figure \ref{figureAsec5} below.

\begin{figure}[h]
  \begin{center}
  \epsfig{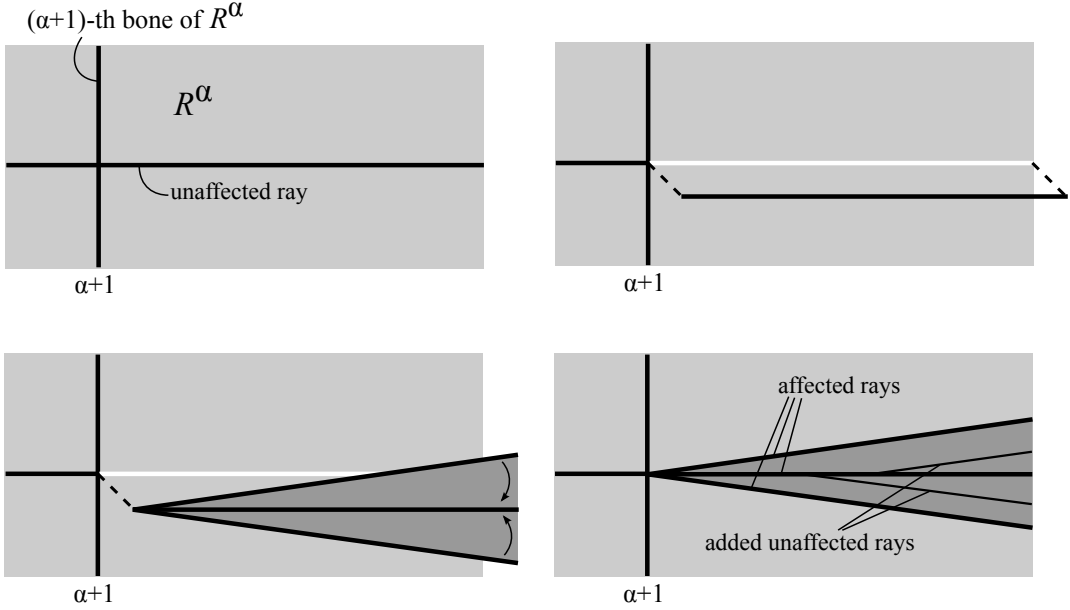}
  \caption{Affecting the rays.}\label{figureAsec5}
  \end{center}
\end{figure}

Then the rays in $\mathcal{D}$ become affected, and are thought to continue through the diagonal in the respective inserted pieces.
New rays appear: the horizontal and vertical lines in the inserted $[\alpha+1,\omega_1)^2$, which run along the affected ray before separating from it
(see Figure \ref{figureAsec5}).
Each is unaffected, except the ones corresponding to $[\alpha+1,\omega_1)\times\{\alpha+1\}$ and 
$\{\alpha+1\}\times[\alpha+1,\omega_1)$.
We use the generic notation $E({\alpha+1})$ for any copy of $[\alpha+1,\omega_1)^2$ inserted at this stage.
Then, the map $k_{\alpha+1,\alpha}:S^{\alpha+1}\to S^{\alpha}$ is defined by simply collapsing each inserted $[\alpha+1,\omega_1)^2$ to its diagonal,
we check as in Lemma \ref{lemmachiant0} that it is indeed closed, so (i) is verified, (ii)  and (iii) follow by construction.
\\
At a limit stage $\alpha$, take the inverse limit of the $R^\beta$ for $\beta<\alpha$, to which we add continuum many copies of $[\alpha,\omega_1)\subset\LL_{\ge 0}$
as follows. 
First, 
whenever $\beta\ge\alpha$,
$E(\alpha)^{\beta}$ denotes the subspace of $R^\beta$ obtained by the expansion of $E(\alpha)$ inside $R^\beta$ until stage $\beta\le\omega_1$.
Fix a
sequence $\alpha_n\nearrow\alpha$ ($n\in\omega$), and
take a nested sequence of $E(\alpha_n +1)$ (that is: $E(\alpha_n +k+1)^\alpha \subset E(\alpha_n +k)^\alpha$).
Then this sequence accumulate to a new copy of $[\alpha,\omega_1)$, as in Figure \ref{figureBsec5} below.

\begin{figure}[h]
  \begin{center}
  \epsfig{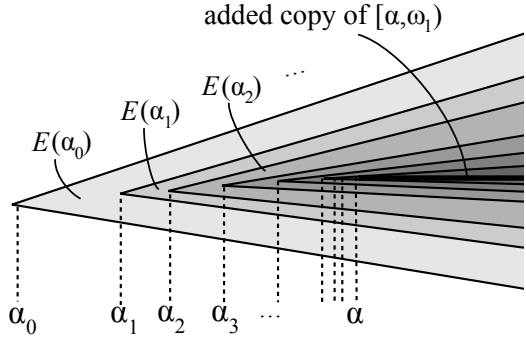}
  \caption{Adding copies of $[\alpha,\omega_1)$.}\label{figureBsec5}
  \end{center}
\end{figure}

This copy is attached to the rays below as indicated on the figure, and forms a new affected ray.
We use the generic notation $L(\alpha)$ for those rays inserted at limit stage $\alpha$.
The map $h_{\alpha,\beta}$ for $\beta\le\alpha$ is obtained by collapsing all pieces inserted at stages above $\alpha$,
and (ii)--(iii) hold by construction.
We let then $R=R^{\omega_1}$ be the inverse limit of the $R^\alpha$, and
define the map $k_{\omega_1,\alpha}:R\to R^\alpha$ accordingly. Then $k_{\omega_1,\alpha}$ is a closed map.
Notice that $R^1$ is essentially the surface given in Example \ref{ex2}.

\begin{lemma}
 \label{nld2}
  Any club subset of $R$ intersects some ray on a club set.
\end{lemma} 
\proof
  Let $C\subset R$ be club. 
  By construction, 
  for each successor $\alpha$, 
  the union of all $E(\alpha)^{\omega_1}$ is dense in $R$. Recall that there are countably many $E(\alpha)$ added at each stage.  
  If $C$ club-intersects some $L(\alpha)$, we are done.
  If there is a nested sequence $E(\alpha_n)$ ($n\in\omega$) such that $C$ club-intersects each $E(\alpha_n)^{\omega_1}$,
  then by construction $C$ club-intersects the $L(\alpha)$ to which these accumulate.
  If $C$ does not club-intersect any $E(\alpha)^{\omega_1}$ or any $L(\alpha)$, then $k_{\omega_1,0}(C)\subset R^0=\LL_{\ge 0}\times\mathbb{S}^1$ must
  club-intersect some ray which is unaffected, and thus $C$ club-intersects it as well.
  \\
  We may thus assume 
  that $C$ club-intersects some $E(\alpha)^{\omega_1}$ but that it does not club-intersect any 
  $E(\beta)^{\omega_1}\subset E(\alpha)^{\omega_1}$ for any $\beta>\alpha$, and does not club-intersect any $L(\beta)$ either.
  Recall that the diagonal in $E(\alpha)$ is affected at stage $\alpha$. 
  If $C$ does not club-intersect this diagonal, then for some $\gamma$,
  $k_{\omega_1,\alpha}(C)$ club-intersects $[\alpha, \gamma]\times [\alpha,\omega_1)\subset E(\alpha)$
  or its mirror image by the diagonal by Lemma \ref{lemmaclassical}.
  But since $k_{\omega_1,\alpha}(C)$ does not club-intersect some affected ray in $E(\alpha)$, 
  it must club-intersect some unaffected ray, and so must do $C$.
  This finishes the proof.
\endproof

\begin{lemma}
 \label{nld3}
  Given two different rays $C,D\subset R$, there is 
  a map $f_{C,D}:R\to\LL_{\ge 0}$ that is bounded on $C$ and unbounded on $D$.
\end{lemma} 

\proof
  Above some height $\alpha$, $C$ and $D$ will be disjoint.
  Take a ray that is `in between' them and that is affected at stage $\beta\ge\alpha$.
  In the corresponding $E(\beta)$, define $f_{C,D}$ as the projection on the second factor if $C$ is `to the right' of $D$ and
  on the first factor if it is `to the left'.
  Then $f$ can be easily extended to all of $R$.
\endproof

The claimed properties of $\langle \mathfrak{D}_{R},\preceq\rangle$ follow directly: a direction must club-intersect exactly one ray,
and two rays are $\preceq$-incomparable.
We notice in passing that under {\bf CH} a variant of this construction enables to affect all rays in the space.
Since any ray becomes a diagonal in some $E(\alpha)$ and a direction is essentially a ray, 
it then is not difficult to find a map $f:S\to\LL_{\ge 0}$ that is bounded
on each direction except a given one. We do not know whether we can build such a surface without {\bf CH}.


\subsubsection{Obtaining $\alpha<\omega_2$}\label{sec:Obtainingalpha}

\begin{figure}[h]
  \begin{center}
  \epsfig{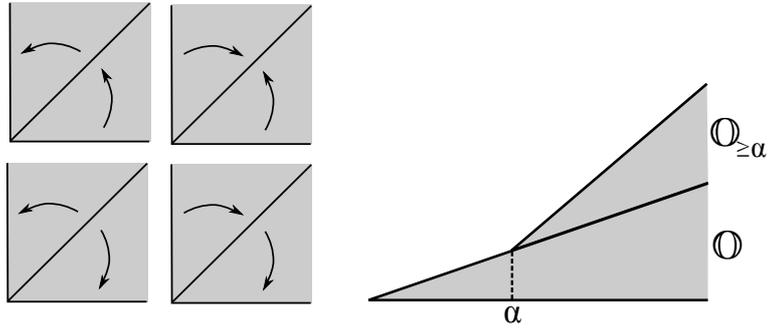}
  \caption{Gluings of the octant: $\uparrow\uparrow$, $\uparrow\downarrow$, $\downarrow\uparrow$, $\downarrow\downarrow$ (left), at height $\alpha$ (right).}
  \label{figure2sec5}
  \end{center}
\end{figure}

Given two copies of $\OO$, one can glue them together
along either
$H$ or $\Delta$, as described in Figure \ref{figure2sec5} (left). 
To each gluing corresponds a function $r:2\to\{\uparrow,\downarrow\}$ in the obvious way.
We can also glue a copy of 
$\OO_{\ge\alpha}=\{\langle x,y\rangle\in\LL_{\ge 0}^2\, :\, \alpha\le y\le x\}$ at height $\alpha$, as in Figure \ref{figure2sec5} (right).

\begin{figure}[h]
  \begin{center}
  \epsfig{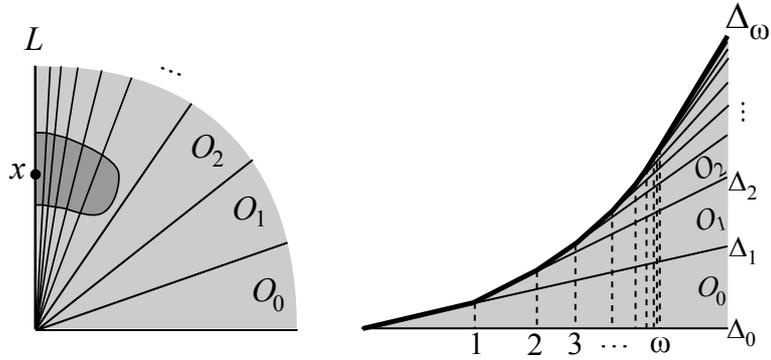}
  \caption{Piling up $\omega$ octants (left), $\Delta_\omega$ (right).}
  \label{figure3sec5}
  \end{center}
\end{figure}

If we pile up $\omega$ copies $O_n$ ($n\in\omega$) of $\OO$, it is then possible to add a copy $L$ of $\LL_{\ge 0}$ that serves as a `top side' to which the
copies accumulate. The neighborhoods of points in this added side are defined similarly as in Example \ref{ex1} (see Figure \ref{figure3sec5}, left).

We can also pile up $\omega_1$ copies of $\OO$ the following way.
We build the surface by induction. 
At each stage $\gamma<\omega_1$, 
the surface will have a free upper side homeomorphic to $\LL_{\ge 0}$, that we will call $\Delta_\gamma$.
During the constructions, we say that a point $x$ in some copy of $\OO$ already included is {\em at height $z\in\LL_{\ge 0}$} if 
$x=\langle z,y\rangle$ for some $y$, and similarly for points in copies of $\LL_{\ge 0}$.
Let $O_{0}=\LL_{\ge 0}$.
If $\gamma=\beta+1$, let $O_\gamma$ be a copy of
$\OO_{\ge \gamma}$ and glue $O_\gamma$ on the free side of $O_\beta$ as 
in Figure \ref{figure2sec5} (right), the orientation of the arrow
being given by $r(\gamma)$. 
The free side of the surface thus built is then the union of the part of the free side 
at stage $\beta$, up to the height $\gamma=\beta+1$, and of the side of $O_\gamma$ that was not glued.
At limit ordinals $\gamma<\omega_1$,
notice that there is a copy $I$ of $[0,\gamma]$ that remains free, but it does not continue above the height $\gamma$
since we added countably many octants `above it' (see Figure \ref{figure3sec5}).
Attach to $I$ a copy of $[\gamma,\omega_1)\subset\LL_{\ge 0}$ as suggested on Figure \ref{figure3sec5} (right) for $\gamma=\omega$, defining the neighborhoods  
as before. This defines $\Delta_\gamma$.
Proceeding up to $\omega_1$ yields an $\omega$-bounded surface.
Notice then that 
when we glue $O_{\beta+1}$ on $O_{\beta}$, the part of the free side of $O_\beta$ at height between $0$ and $\beta+1$ remains free, 
the union of all these free sides give a
copy of $\LL_{\ge 0}$ that we call $\Delta_{\omega_1}$.

\begin{figure}[h]
  \begin{center}
  \epsfig{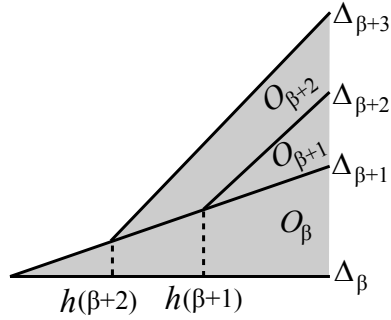}
  \caption{Gluing according to $h$.}
  \label{figure4sec5}
  \end{center}
\end{figure}

We can generalize a bit more and
pile up $\alpha$ copies of $\OO$ for any $\omega_1\le\alpha<\omega_2$ as follows.
First, fix a (non-continuous) $1$-to-$1$ map $h:\alpha\to\omega_1$, which gives the height at which the $\gamma$-th copy of $\OO$ is glued,
and again a map $s:\alpha\to \{\uparrow,\downarrow\}$, which gives the orientation of the gluing.
(The previous construction is the particular case $\alpha=\omega_1$ and $h$ the identity.)
That is, we glue $O_{\beta+1}$ on the free side of $O_\beta$ at height $h(\beta+1)$ (see Figure \ref{figure4sec5}), and define 
$\Delta_{\beta+1}$ similarly as before.
If $\gamma\le\alpha$ is limit, things are a bit more complicated, 
but since $h$ is $1$-to-$1$, at a countable height only countably many copies of the octant appear.
Let us look between height $\beta$ and $\beta+1$.
If for some $\xi<\gamma$, $h(\xi')>\beta$ for each $\xi\le\xi'<\gamma$, then 
the piece of free side at this height will remain untouched after stage $\xi$,
otherwise there will be no free side left at these heights.
(Notice that if $\text{\rm cf}(\gamma)=\omega_1$, we are always in the former case.)
If there is no free side left, we add a copy of $[\beta,\beta+1]$ to which the (countably many) octants below accumulate, as before. The union of 
the remaining parts of the free sides and the one added is then a copy of $\LL_{\ge 0}$ that we call again $\Delta_\gamma$.
A particular case is whenever $\gamma=\alpha$, which gives the upper side $\Delta_\alpha$, as in the case $\alpha=\omega_1$.
We denote by $S^{\alpha,h,r}$ the $\omega$-bounded surface (with boundary) thus obtained.
A canonical sequence $\cup_{\xi\in\omega_1}S^{\alpha,h,r}_\xi$
for $S^{\alpha,h,r}$ is given by open sets defined as the union of points at height $\le\xi$ for each $\xi<\omega_1$.

\begin{lemma}\label{lemmachiant}
   Let $\alpha,h$ and $r$ be as above, and $C\subset S^{\alpha,h,r}$ be club.
   Then either $C\cap\Delta_\gamma$ is club for some $\gamma\le\alpha$, or 
   $C\cap O_\gamma$ is club for some successor $\gamma\le\alpha$.
\end{lemma}
\proof
   For clarity, we omit the superscripts $h$ and $r$ when they are not absolutely necessary,
   and write $S^{\alpha}=\cup_{\xi\in\omega_1}S^{\alpha}_\xi$.
   If $\beta\le\alpha$, we denote by $S^{\beta}$ the subspace of $S^\alpha$ given by $S^{\beta,h',s'}$, where 
   $h',s'$ are the restrictions to $\beta$ of $h,s$.\\
   The proof is by induction on $\alpha$.
   If $\alpha=0$, the result is immediate.
   If $\alpha=\beta+1$ then either $C$ club-intersects 
   $O_{\beta+1}\subset S^{\alpha}$ or it club-intersects $S^{\beta}$.   
   If $\text{\rm cf}(\alpha)=\omega$, either $C$ club-intersects
   $\Delta_\alpha$, or it club-intersects some $S^{\beta}$ for a smaller $\beta$.
   In both cases we are done by induction.
   \\  
   Assume thus that $\text{\rm cf}(\alpha)=\omega_1$. Fix a closed increasing
   $\omega_1$-sequence $\alpha_\xi$ converging to $\alpha$.
   Since $h$ is $1$-to-$1$, there is a club subset of $\xi$
   such that $h(\alpha_{\xi'})\ge\xi$ whenever $\xi'\ge\xi$.   
   We can thus assume that this holds for each $\xi$, up to taking a subsequence.
   This implies that any member of $(S^\alpha - S^{\alpha_{\xi}})$ has height at least $\xi$.
   If $C\cap S^{\alpha_{\xi}}$ is club for some $\xi$ we are done by induction, so assume that 
   this intersection is bounded by $\mu(\xi)\in\omega_1$.
   Again, this implies that for a club subset $E\subset\omega_1$, $C\cap S^{\alpha_{\xi}} \subset S^\alpha_\xi$ when $\xi\in E$.
   \\
   Let $\xi(0)\in E$, we shall find a point at height at least $\xi(0)$ in $C\cap\Delta_\alpha$, this implies that this intersection is club.
   Given $\xi(n-1)$,
   Let $x_n$ be a point at height at least $\xi(n-1)+1$ in $C\cap(S^{\alpha} - S^{\alpha_{\xi(n-1)}})$.
   Let $\xi(n)$ be the smallest point in $E$ that is $\ge$ the height of $x$.
   By definition, the limit $x$ of a subsequence of the $x_n$ will be a member of $C\cap (S^{\alpha}_\xi - S^{\alpha_{\xi}})$, for 
   $\xi=\sup_n\xi(n)$. But by construction, it belongs to $\Delta_{\xi}$.
   Moreover, since $h$ does not go below $\mu$ after stage $\alpha_\xi$, $\Delta_{\xi}\cap S^\alpha_{\alpha_\xi}$ is included
   in $\Delta_\alpha$, so $x\in \Delta_\alpha$ as well.
\endproof

A corollary of Lemma \ref{lemmachiant} and Proposition \ref{prop1} (b) is:

\begin{cor}\label{corchiant}
   Any direction in
   $S^{\alpha,h,r}$ is ufo-equivalent to some $\Delta_\beta$, for $\beta\le\alpha$.
\end{cor}

\begin{rem}
We can weaken the assumption that $h:\alpha\to\omega_1$ is $1$-to-$1$ to 
$h$ being $\omega$-to-$1$ without affecting the construction. 
\end{rem}

Notice in passing:

\begin{lemma}\label{corchiant0}
   Let $\alpha<\omega_2$, $h:\alpha\to\omega_1$ be $1$-to-$1$ (non-continuous) and $r:\alpha\to\{\uparrow,\downarrow\}$.
   Then
   $S^{\alpha,h,r}$ is spanned by directions.
\end{lemma}
\proof
By Lemma \ref{lemmachiant},
a club in $S^{\alpha,h,r}$ either club-intersects some $\Delta_\gamma$, and thus contains a copy of $\omega_1$,
or it club-intersects the interior of some $O_\gamma$. In the latter case, it must club-intersect
some horizontal in $O_\gamma$ by Lemma \ref{lemmaclassical}, and thus also contains a copy of $\omega_1$.
\endproof

\begin{ex}\label{ex3} 
   For each $\alpha<\omega_2$, there is a
   surface $S$ with 
   $\langle {\mathfrak{D}}_{S},\preceq\rangle= \langle\alpha,\le\rangle$. 
   If $\alpha$ is successor, then $S$ can be made $\omega$-bounded.
\end{ex}

If $\alpha=\beta+1$,
take any (non-continuous) $1$-to-$1$ $h:\beta\to\omega_1$, and $s$ constant on $\uparrow$. 
Set $S=S^{\beta,h,s}$.
If $\alpha$ is limit, set $\alpha=\beta$, choose $h$ and $s$ as in the successor case, and set
$S=S^{\beta,h,s} - \Delta_\alpha$.
The claimed properties of $\langle {\mathfrak{D}}_{S},\preceq\rangle$ follow from Corollary \ref{corchiant}
and the next lemma.

\begin{lemma}\label{lemmachiant2}
    Let $\alpha<\omega_2$, $h:\alpha\to\omega_1$ be $1$-to-$1$ (non-continuous), and $r:\alpha\to \{\uparrow,\downarrow\}$ be constant on $\uparrow$.
    Then, $\Delta_{\gamma'}\prec\Delta_{\gamma}$ whenever $\gamma'<\gamma\le\alpha$.
\end{lemma}
\proof
    (As above, we omit the superscripts $h$ and $r$.)
    The fact that $\Delta_{\gamma'}\preceq\Delta_{\gamma}$ is proved by induction on $\gamma$.
    Proposition \ref{prop1} (b) gives the successor steps.
    For the limit steps, one first shows that if some open set $U$ contains $\Delta_{\gamma}$, then 
    there are $\eta<\gamma$ and $\xi\in\omega_1$ such that
    outside of $S^\alpha_\xi$, $U$ contains $O_{\gamma'}$ for each $\eta\le\gamma'\le\gamma$.
    Indeed, if it is not the case, the complement of $U$ is unbounded in 
    $S^{\gamma'}$ for each $\gamma'<\gamma$. 
    By Lemma \ref{lemmachiant} this implies that this complement intersects $\Delta_\gamma$, a contradiction.
    Thus, if $f:S^{\alpha}\to\LL_{\ge 0}$ is bounded by $\xi$ on $\Delta_\gamma$, there is some $\eta$ such that
    $f^{-1}([0,\xi+1))$ contains the terminal part of each $O_{\gamma'}$ for $\eta\le\gamma'\le\gamma$.
    This shows that $f$ is bounded on $\Delta_{\gamma'}$ for these $\gamma'$, and by induction $f$ is bounded on each 
    $\gamma'<\gamma$.\\
    To see that the inequality is strict, let $\gamma <\alpha$.
    Define $f$ to take the value $h(\gamma+1)$ on $S^{\gamma}$. By construction
    the octant $O_{\gamma+1}=\OO_{\ge h(\gamma+1)}$ is glued at height $h(\gamma+1)$. We may thus define $f$ to be the projection
    on the second coordinate in $O_{\gamma+1}$.
    Points in $O_{\eta}$ for $\eta>\gamma$ at height below $\gamma+1$ are sent by $f$ to $h(\gamma+1)$, and to their
    second coordinate if their height is above $\gamma+1$.
    Then, $f$ is continuous, bounded on $\Delta_{\gamma'}$ if $\gamma'\le\gamma$,
    and unbounded if $\gamma' >\gamma$.
\endproof

It is important to notice that if $\alpha\ge\omega_1$, we do {\em not} obtain
the reverse order $\alpha^*$ by simply reversing each $\uparrow$ to $\downarrow$, as 
the $\Delta_\beta$ for $\beta$ of cofinality $\omega_1$ will not be $\preceq$ those for smaller $\beta$,
see Theorem \ref{propufo}. 

\begin{prob}\label{pomega2}
  Is there a Type I (first countable) space (or manifold) $X$ with $\langle {\mathfrak{D}}_{X},\preceq\rangle\simeq \langle\omega_2,\le\rangle$~?
  One with $\langle {\mathfrak{D}}_{X},\preceq\rangle\simeq \langle\alpha^*,\le\rangle$ for some $\alpha\ge\omega_1$~?
\end{prob}


\subsubsection{Obtaining an Aronszajn or Suslin line}

The same type of ideas can be used to obtain more `exotic' ordered types for the ufo-equivalent classes of direction. 
Recall that $2^{<\omega_1}$ is the tree of $0,1$-sequences of length $<\omega_1$ ordered by the inclusion.
Let 
$T\subset 2^{<\omega_1}$ be a subtree (in this example we assume that
subtrees are downward closed in $2^{<\omega_1}$).
We let $\wt{T}\subset 2^{\le\omega_1}$ be $T$ union 
all its maximal branches, that is, given a maximal branch in $T$ we add to $T$ the union of its members (which defines a map $\alpha\to 2$ for 
some $\alpha\le\omega_1$).
Given $x,y\in\wt{T}$, we let $d(x,y)$ be the smallest $\alpha$ such that
$x(\alpha) \not= y(\alpha)$, 
and write $x \vartriangleleft y$ if $x(d(x,y)) < y(d(x,y))$. This defines a total order on $\wt{T}$.
 
\begin{ex}\label{ex4}
   Given an Aronszajn subtree $T$ of $2^{<\omega_1}$ 
   without maximal element,
   there is
   an $\omega$-bounded surface $A_T$ such that $\langle {\mathfrak{D}}_{A_T},\preceq\rangle$ 
   is isomorphic to the set of maximal branches of $T$ ordered as above.
\end{ex}
 
So, if $T$ is Suslin, ${\mathfrak{D}}_{A_T}$ is in fact a (compact) Suslin line.
We need $T$ to be binary 
(or at least that the number of immediate successors
of each point of $T$ is uniformly bounded by some $n\in\omega$)
for the surface to be $\omega$-bounded.
\\
The construction works as follows and is a kind of mix between Examples \ref{exDominant} and \ref{ex3}.  
Start with a copy $L_\varnothing$ of $\LL_{\ge 0}$.
To each $t\in T$ at successor level $\alpha$ will correspond a 
copy $O_t$ of $\OO_{\ge\alpha}$. 
Then $L_{t\smallfrown 0}$ is the horizontal $\{\langle x,\alpha\rangle\,:\, x\ge\alpha\}\subset O_t$ and
$L_{t\smallfrown 1}$ the diagonal $\Delta\subset O_t$. 
This defines $L_t$ for $t$ at successor levels once $O_t$ is defined.
We shall also define $L_x$ for $x\in \wt{T}$ at limit level $\alpha$. 
Each will be a copy of $[\alpha,\omega_1)\subset\LL_{\ge 0}$.
While the interior of each $O_t$ remains in the final construction as well as $L_x$ for $x\in\wt{T}-T$, $L_t$ for $t\in T$
will be affected as in Example \ref{exDominant}.
\\
Given $t\in\text{Lev}_{\alpha}(T)$, take out $L_t$ above height $\alpha$ and sew in 
$O_t$, such that the arrow $\uparrow$ points to the $L_x$  for $t\vartriangleleft x$ as in Figure \ref{figureXsec5} below.
This defines $L_{t\smallfrown 0}$ and $L_{t\smallfrown 1}$.
At a limit stage $\beta$, insert copies $L_x$ of $[\beta,\omega_1)\subset\LL_{\ge 0}$ for each $x\in \wt{T}$ at height $\beta$.
These copies are inserted according to $\vartriangleleft$: $L_x$ 
will be exactly between the $O_t$ and $O_s$ ($s,t$ at lower levels) such that
$t\vartriangleleft x \vartriangleleft s$, as in Example \ref{ex1}.

\begin{figure}[h]
  \begin{center}
  \epsfig{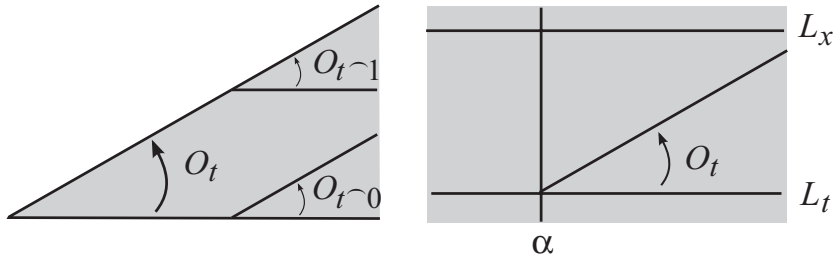}
  \caption{Inserting $O_t$.}\label{figureXsec5}
  \end{center}
\end{figure}

As in Example \ref{ex3}, we call $A_T^{\alpha}$ the space
obtained after completing the induction up to each $t$ at level $\alpha$, and $A_T=A_T^{\omega_1}$ is the inverse limit of these.
It should be clear that $A_T$ is indeed $\omega$-bounded.

\begin{lemma}
   Any club $C$ in $A_T$ club-intersects either the interior of some $O_t$ or some $L_x$.
\end{lemma}
\proof
   This is very similar to Lemma \ref{nld2} (but easier).
   First, for each $\alpha$ there is a closed map $k_{\omega_1,\alpha}:A_T^{\omega_1}\to A_T^\alpha$ that
   collapses the rays and octants introduced after stage $\alpha$.
   If $C$ does not club-intersect the interior of some $O_t$, then
   $k_{\omega_1,0}(C)$ club-intersects $L_0$ or $L_1$, say that the former holds.
   Thus, $C$ club-intersects the union of the $L_x$ for $x \supset 0$ in $\wt{T}$.
   By induction, define $t(\alpha)\in\text{\rm Lev}_\alpha (\wt{T})$ such that 
   $k_{\omega_1,\alpha}(C)$ club-intersects $L_{t(\alpha)}$. 
   At successor levels this is as for $L_0$, at limit levels we use the fact that
   a countable nested intersection of clubs is club in $A_T$ by $\omega$-boundedness.
   Since $T$ is Aronszajn, $t(\alpha)\in \wt{T} - T$ for some $\alpha$, and 
   thus $C$ club-intersects this $L_x$.
\endproof

The claim about $\langle {\mathfrak{D}}_{A_T},\preceq\rangle$ follows as in Example \ref{ex3}, so
we do not provide details.

\vskip .2cm
Notice that all the examples in Section \ref{secconstr} are surfaces with a top and a bottom side, which can be used in place of the octant 
to build $\omega$-bounded surfaces with more and more complicated ufo-equivalent classes of directions.


\subsection{Generalizations}\label{sec:generalization}

In this section we generalize some of the constructions done in Section \ref{secconstr}.


\subsubsection{$\Sigma(X,s)$}\label{sec:Sigma}

  Let $X$ be a Type I space with a slicer $s$.
  We define $\Sigma(X,s)$ as the space given by 
  $\bigl( X\times[0,1)\bigr) \sqcup L$, $L$ being a copy of $\LL_{\ge 0}$, where the points 
  in $X\times[0,1)$ have the usual topology,
  and neighborhoods of 
  $x\in L$ are given by 
  all $W(y,z,a)$ for $y<x<z\in L$, $0\le a < 1$, defined as the union of
  $(y,z)\subset L$ with $s^{-1}((y,z))\times (a,1)\subset X\times[0,1)$. Then
  $\Sigma(X,s)$ is a Type I space with canonical sequence
  $\Sigma(X,s)_\alpha=X_\alpha\times[0,1)\sqcup[0,\alpha)$, whose properties follow closely those of $X$, as the next few lemmas show.
  
\begin{rem}
  $\Sigma(X,s)$ may have fewer open
  sets than the mapping cylinder (as defined in classical homotopy theory) of the scale $s:X\to\LL_{\ge 0}$.
\end{rem}

For brevity, we shall often omit $s$ and write $\Sigma(X)$ in place of $\Sigma(X,s)$.

\begin{lemma} \label{lemmasigmaA}
  If $X$ is regular, so is $\Sigma(X)$.
\end{lemma}
\proof
  Let $x\in U\subset \Sigma(X)$, with $U$ open. If $x\not\in L$, by regularity of $X\times[0,1)$ there is 
  an open $V$ such that $x\in V\subset \wb{V}\subset U$.
  If $x\in L$, there are $y<x<z$ and $a<1$ such that $W(y,z,a)\subset U$.
  Take $y',z'$ with $y<y'<x<z'<z$ and $a'$ between $a$ and $1$, then the following inclusions hold:
  $$
   x\in W(y',z',a')\subset\wb{W(y',z',a')}\subset \left( s^{-1}([y',z'])\times[a',1)\right) \sqcup [y',z']\subset W(y,z,a)\subset U.
  $$
\endproof

\begin{lemma} \label{lemmasigmaB}
  If $X$ is locally second countable, then so is $\Sigma(X)_\alpha$ for each $\alpha\in\omega_1$.
\end{lemma}
\proof
  $\wb{X_\alpha}$ is Lindel\"of and locally second countable hence second countable, 
  thus so is $X_\alpha\times[0,1)$. Given a countable basis for $X_\alpha\times[0,1)$,
  one obtains a countable basis for $\Sigma(X)_\alpha$ by adding the sets containing the $W(x,y,a)$ for rational $x,y,a$.
\endproof

\begin{cor}\label{corsigmaC}
  If $X$ is regular and locally second countable, then $X_\alpha$ and $\Sigma(X)_\alpha$
  are metrizable for each $\alpha\in\omega_1$.
\end{cor}
\proof
  By Urysohn's theorem and the previous lemma.
\endproof

We also have:

\begin{lemma}\label{lemmasigmaD}
  $X$ is $\omega$-bounded ($\omega_1$-compact, respectively) iff $\Sigma(X)$ is $\omega$-bounded ($\omega_1$-compact, respectively).
\end{lemma}
\proof
  The reverse implications are immediate. 
  If $X$ is $\omega_1$-compact, then $\Sigma(X)$ is a countable union of $\omega_1$-compact closed subspaces, and is
  thus $\omega_1$-compact as well.
  Suppose that $X$ is $\omega$-bounded, and let $E\subset\Sigma(X)$ be countable, thus $E\subset\Sigma(X)_\alpha$ for some $\alpha$.
  Let a cover of $\wb{E}$ be given.
  Since $\wb{E}\cap L \subset [0,\alpha]$, it is compact, take a finite subcover of it. 
  Then for some $a<1$ the union of its members covers $\wb{E}\cap X\times(a,1)$.
  Since $X\times[0,a]$ is $\omega$-bounded, we take a finite subcover of $\wb{E}\cap (X\times[0,a])$ to conclude.
\endproof

Notice however that $\Sigma(X)$ might be non-locally compact even if $X$ is.  
If $X$ contains no direction, then $\Sigma(X)$ contains club subsets (copies of $X$) that do not contain a direction in itself.
Call {\em horizontals} the subsets of $\Sigma(X)$ of the form $X\times\{x\}$ or $L$, 
and call $L$ 
the {\em core} of $\Sigma(X)$.

\begin{lemma}\label{lemmadt3}
    Let $X$ be Type I such that uncountably many bones of $X$ are non-empty.
    If $U$ is open and contains a club subset $C$ of the core of
    $\Sigma(X)$, then $U\supset (X-X_\alpha)\times (a,1)$ for some $0\le a < 1$ and $\alpha\in\omega_1$. 
  \end{lemma}
  \proof
    Suppose that it is not the case.
    Since $C$ is club in $L$ it contains a club subset of $\omega_1$.
    Start with $\alpha_0=0$. Given $\alpha_n$, 
    choose $\alpha_{n+1}\in C$ such that $\alpha_{n+1}>\alpha_n$ and there is some $x_n\in \bigl(B_{\alpha_{n+1}}(X)\times (1-\frac{1}{n},1)\bigr) - U$.
    Then the $x_n$ converge to $\alpha=\sup_{n\in\omega}\alpha_nÊ\in C$, a contradiction.
  \endproof


\subsubsection{$\we{X,s}$, $\ve{X,s}$}\label{sec:UDDelta}

$\we{X,s}$ and $\ve{X,s}$ are generalized octants.

\begin{defi} 
   Let $X$ be a Type I space possessing a slicer $s:X\to \LL_{\ge 0}$.
   We set:
     \begin{align*}
          \we{X,s}&=\{\langle x,y\rangle\in X\times\LL_{\ge 0}\,:\, y\ge s(x)\}, \\
          \ve{X,s}&=\{\langle x,y\rangle\in X\times\LL_{\ge 0}\,:\, y\le s(x)\},\\
          \de{X,s}&=\{\langle x,y\rangle\in X\times\LL_{\ge 0}\,:\, y = s(x)\}.
     \end{align*}  
A slicer for any subspace of $X\times\LL_{\ge 0}$ is given
by $\wt{s}(\langle x,y\rangle)=\max\{s(x),y\}$.
\end{defi}

\begin{figure}[h]
  \begin{center}
  \epsfig{figure=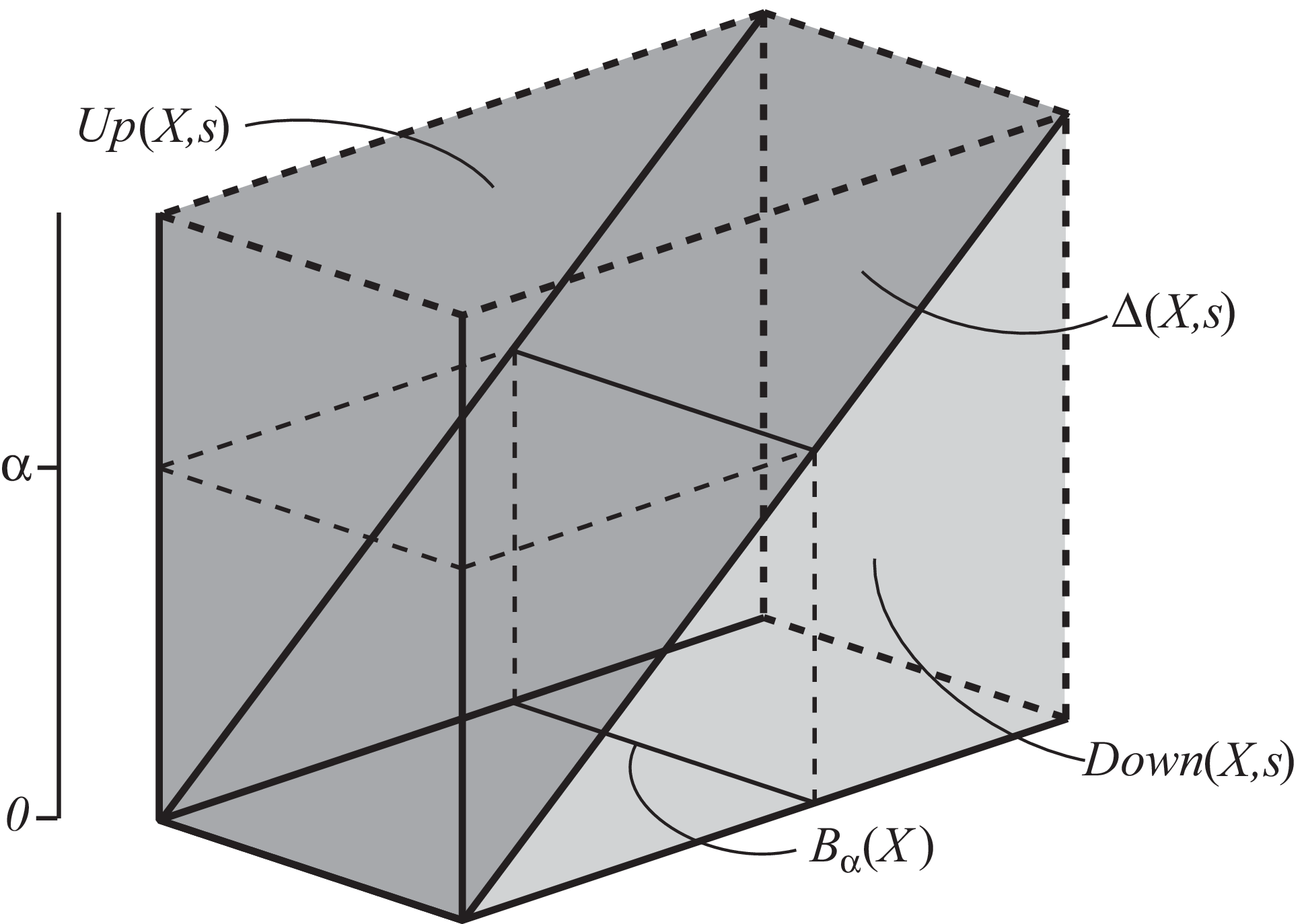, width=8cm}\,
  \epsfig{figure= 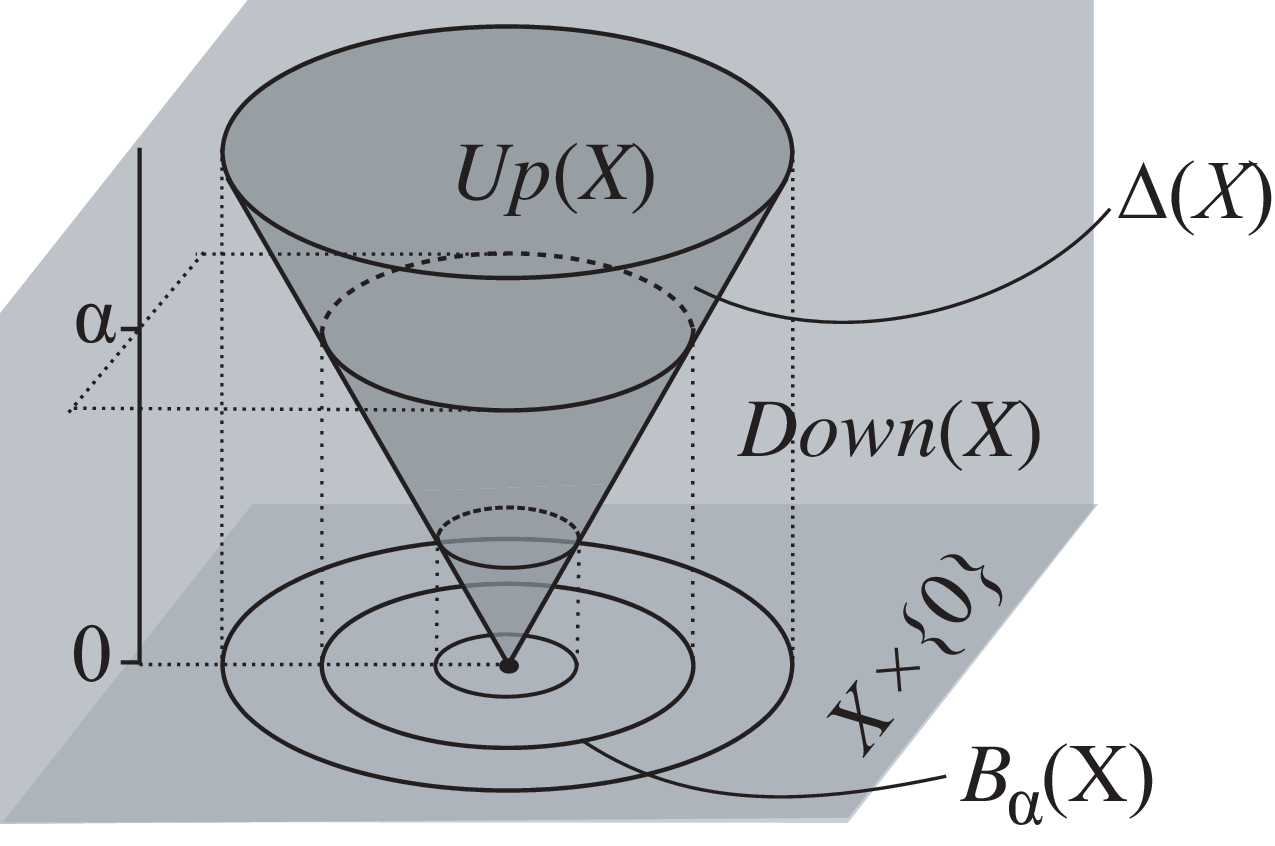, width=6cm}
  \caption{$\we{X}$, $\ve{X}$ and $\de{X}$.}\label{figureWsec6}
  \end{center}
 \end{figure}

They are depicted in Figure \ref{figureWsec6}.
The righthandside shows a somewhat accurate picture when $X$ is a surface, while in the lefthandside we 
depicted $X$ as a strip, which may be less `realistic' (whatever that means) but makes some
constructions easier to describe graphically.
As for $\Sigma(X,s)$, we shall write $\we{X},\ve{X},\de{X}$ for brevity except when the choice of the 
slicer is relevant.
Notice that $\de{X}$, being the graph of $s$, is homeomorphic to $X$ through the map 
$\langle x,y\rangle\mapsto x$.
If $X$ is a manifold (in particular, a surface), then $\we{X}$ is a manifold with boundary $\de{X}$, while
$\ve{X}$ is `almost' a manifold with two boundaries
$X\times\{0\}$ and $\de{X}$, both homeomorphic to $X$. 
The only point where the space fails to be a manifold
is $\langle x_0,0\rangle$, where $x_0$ is the (only) $s$-preimage of $0$. 
By either deleting it, or through a small local change, one obtains a manifold as well.

\begin{lemma}\label{projclosed}
  Let $X$ be Type I with a slicer $s$. 
  Then the following holds.\\
  (a) If $X$ is countably tight, then the projection on the first factor $\pi:X\times\LL_{\ge 0}\to X$ is a closed map.
      Moreover, $\pi\restrict{\ve{X}}$ and $\pi\restrict{\we{X}}$ are closed maps.\\
  (b) If $X$ is locally compact, then $\pi\restrict{\ve{X}}:\ve{X}\to X$ is a closed map.
\end{lemma}
\proof
  \
  \\
  (a)
  Let $C\subset\ve{X}$ be closed, $x\in X$ be in the closure of $\pi(C)$ and
  let $A=\{x_n\,:\,n\in\omega\}\subset\pi(C)$ be such that $x\in\wb{A}$.
  There is thus a sequence $y_n\in\LL_{\ge 0}$ such that $\langle x_n,y_n\rangle\in C$.
  Let $y$ be an accumulation point of the $y_n$. 
  Then $\langle x,y\rangle$ is an accumulation point of the $\langle x_n,y_n\rangle$, and thus in $C$, so $x\in \pi(C)$.
  The `moreover' part follows from the fact that both $\we{X}$ and $\ve{X}$ are closed in $X\times\LL_{\ge 0}$.
  \\
  (b)
  Notice that $\pi\restrict{\ve{X}}$ is proper, that is, the preimage of a compact subset is compact.
  Indeed, a compact subset $K\subset X$ is contained in some $X_\alpha$, so 
  $\pi^{-1}(K)$ is a closed subset of $K\times[0,\alpha]$, and is thus compact. 
  If $X$ is locally compact, so is $\ve{X}$.
  A proper map from a locally compact space to a Hausdorff space is closed.
  \endproof

\begin{lemma}\label{nodirinveX}
   If $X$ contains no direction and is either countably tight or locally compact, then $\ve{X}$ contains no direction.
\end{lemma}
\proof
   Let $C$ be club in $\ve{X}$. 
   Then $\pi(C)$ is club as well by Lemma \ref{projclosed}, and
   let $f:X\to\LL_{\ge 0}$ witness the fact that $\pi(C)$ is not a direction. 
   So $f^{-1}(\{0\})\cap \pi(C)$ is unbounded while $f$ is unbounded on $\pi(C)$.
   Then, $g:\ve{X}\to\LL_{\ge 0}$, $\langle x,y\rangle\mapsto f(x)$ shows that $C$ is not a direction.
\endproof

\begin{lemma}\label{closeddiscreteweX}
   Let $X$ be either countably tight or locally compact, and
   $C\subset\ve{X}$. Then the following holds.\\
   (a) If $\pi(C)$ contains a closed 
   discrete subset, then so does $C$.\\
   (b) If $C$ is closed discrete, then so is $\pi(C)$.
\end{lemma}
\proof
   \ \\
   (a) Let $D\subset \pi(C)$ be closed discrete.
   For each $x\in D$, take $y(x)\in \pi^{-1}(\{x\})\cap C$.
   Then $E=\{y(x)\,:\,x\in D\}$ is closed: take any $z\not\in E$, if $\pi(z)\not\in D$ then 
   $\pi^{-1}(X-D)\cap E=\varnothing$, and if $\pi(z)=\pi(x)$ for some $x\in D$,
   there is some open interval $(a,b)\subset\LL_{\ge 0}$ such that
   $\pi^{-1}\bigl(X-(D-\{x\})\bigr)\cup \bigl(\{\pi(x)\}\times (a,b)\bigr)$ contains $z$ and does not intersect $E$.\\
   (b) Similar to Lemma \ref{projclosed}, we skip the proof.
\endproof

\begin{lemma}\label{lemma:UDomega_1cpct}
   Let $X$ be a Type I space with a slicer $s$.
   Then the following holds.\\
   (a) If $X$ is $\omega$-bounded, then so are $X\times\LL_{\ge 0}$, $\we{X}$ and $\ve{X}$.\\
   (b) If $X$ is $\omega_1$-compact and either locally compact or countably tight, then $\ve{X}$ is $\omega_1$-compact.\\
   (c) If $X$ is first countable and $\omega_1$-compact, then $X\times\LL_{\ge 0}$ and $\we{X}$ are $\omega_1$-compact.
\end{lemma}
\proof
\ \\
(a) $\omega$-boundedness is productive, hence $X\times\LL_{\ge 0}$ is $\omega$-bounded.
    Since $\we{X}$ and $\ve{X}$ are closed subsets of $X\times\LL_{\ge 0}$ they inherit the property.
\\
(b) By Lemmas \ref{projclosed} and \ref{closeddiscreteweX}.
\\
(c) Let $D$ be a closed uncountable discrete subset of $X\times\LL_{\ge 0}$.
    By Lemma \ref{projclosed} (a)
    $\pi(D)$ is closed.
    If $\pi(D)$ is unbounded, choose by 
    induction $x_\alpha\in D$ and $\beta_\alpha\in\omega_1$ such that
    $\pi(x_\alpha)\in X_{\beta_{\alpha+1}}-\wb{X_{\beta_\alpha}}$.
    Then $\{x_\alpha\,:\,\alpha\in\omega_1\}$ is closed and its projection is a club discrete subset of $X$
    by Lemma \ref{closeddiscreteweX}, which is impossible. 
    Thus $D\subset X_\alpha\times\LL_{\ge 0}$. By the next lemma, 
    $D$ has a club intersection with some $\{x\}\times \LL_{\ge 0}$ and cannot be discrete.
    Thus $X\times\LL_{\ge 0}$ 
    and $\we{X}$ (as a closed subset) are $\omega_1$-compact.
\endproof

The missing lemma is very similar to the sublemma after Lemma 2.2 in \cite{BGGG}.
The proof goes in part as that of Lemma \ref{lemmadt3}. We state a more general version than
the one just needed since we will use it again later.

\begin{lemma}\label{lemma:lemma2.2BGGG}
  Let $L$ be $\omega$-bounded and Type I, and $Y$ be first countable and Lindel\"of. Let
  $D\subset L$ be club.
  Then any open set $U\subset Y\times L$ containing $\{y\}\times D$ contains $V\times Y$ for some open $V\subset Y$,
  and any club in $Y\times L$ club-intersects $\{y\}\times L$ for some $y\in Y$.
\end{lemma}
\proof
Since $D$ is also an $\omega$-bounded Type I space, we show the first assertion with $D=L$.
Let $U\supset \{y\}\times L$ and suppose it does not contain $V\times L$ for any open $V\subset Y$.
Take a countable base of neighborhoods $V_n$ whose intersection is $\{y\}$ and
choose a point $\langle y_n,x_n\rangle\in (V_n\times L) - U$ for each $n$.
Let $x$ be an accumulation point of the $x_n$, then $\langle y,x\rangle$ is an accumulation point of the sequence
$\langle y_n,x_n\rangle$ and thus does not belong to $U$, a contradiction.\\
Let now $C\subset Y\times L$ be closed,
if $C\cap \{y\}\times L$ is bounded for each $y\in Y$, then
$L-C\supset V(y) \times (L-L_{\alpha(y)})$ for some open $V(y)\ni y$.
Take a countable subcover $V(y_n)$, then $C\subset L_\alpha$ for $\alpha=\sup_n\alpha(y_n)$.
\endproof

\begin{lemma}\label{lemma3d1}
  Let $X$ be a first countable $\omega_1$-compact Type I space possessing a slicer $s$, and let $E\subset\we{X}$ be club.
  Then either there is an $\alpha$ such that $E\subset (X_\alpha\times\LL_{\ge 0})$, or 
  $E\cap \de{X}$ is club.
\end{lemma}  
\proof 
$\we{X}$ is $\omega_1$-compact by Lemma \ref{lemma:UDomega_1cpct} (c).
Let $f:\we{X}\to\LL_{\ge 0}$ be given by $f(x,y)=s(x)$. 
If $f$ is unbounded on $E$, the subset $C(f,E)\subset E$ given by Lemma \ref{lemmastar}
is club, and by definition if $\langle x,y\rangle\in C(f,E)$ then $y=s(x)$ so this point belongs to $\de{X}$.
If $f$ is bounded by $\alpha$ on $E$, then $E\subset X_\alpha\times\LL_{\ge 0}$.
\endproof

\begin{lemma}\label{lemma3d2}
  Let $X$ be a first countable, $\omega_1$-compact Type I space possessing a slicer $s$.
  If $U\subset \we{X}$ is a large open set, it contains a copy of $\omega_1$ located in $X_\alpha\times\LL_{\ge 0}$ for some $\alpha$.
\end{lemma}
\proof
  Let $C\subset U$ be club. If $C\cap (\wb{X_\alpha}\times\LL_{\ge 0})$ is unbounded, by Lemma \ref{lemma:lemma2.2BGGG} 
  $C$ club-intersects $\{x\}\times\LL_{\ge 0}$
  for some $x\in \wb{X_\alpha}$.
  If $C\cap\wb{X_\alpha}\times\LL_{\ge 0}$ is bounded for each $\alpha$, then
  $C\cap\de{X}$ is club by Lemma \ref{lemma3d1}.
  Now if $U\cap\wb{X_\alpha}\times\LL_{\ge 0}$ is bounded for each $\alpha$, 
  build a sequence $x_\gamma\not\in U$ as in Lemma \ref{lemmastar} to obtain a point in $C - U$, a contradiction.
\endproof

\begin{lemma}\label{lemma3d3}
  Let $X$ be a first countable $\omega_1$-compact Type I space with a slicer $s$ and possessing no direction.
  If $D\subset \we{X}$ is a direction, it is contained in $X_\alpha\times\LL_{\ge 0}$ for some $\alpha$.
\end{lemma}

\proof
  Otherwise, by Lemma \ref{lemma3d1} its intersection with $\de{X}$ is club. 
  Let $h:X\to\de{X}$ be the homeomorphism $h(x)=\langle x,s(x)\rangle$, let $g:X\to\LL_{\ge 0}$ witness the fact that
  $h^{-1}(D\cap\de{X})$ is not a direction. Then $f:\we{X}\to\LL_{\ge 0}$ given by $f(x,y)=g(x)$ 
  shows that $D$ is not a direction.
\endproof


\subsubsection{$\Omega(X,A,s)$}\label{sec:Omega}

\begin{figure}
  \begin{center}
  \epsfig{figure=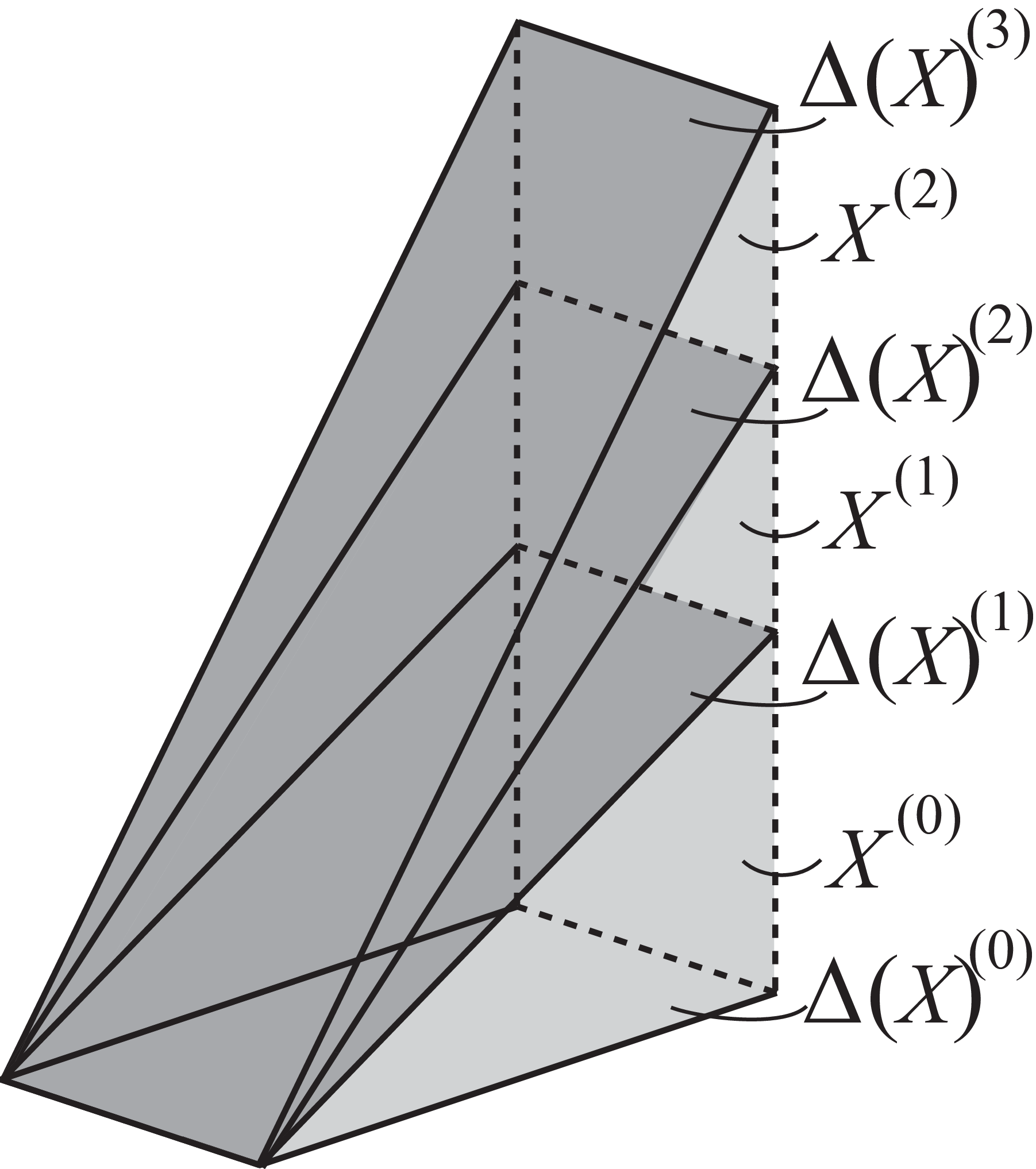, width=6cm}
  \caption{Piling up copies of $\ve{X}$.}\label{figureOmegaX}
  \end{center}
 \end{figure}

We may pile up copies of $\ve{X}$ like we did with the octant as on Figure \ref{figureOmegaX}.
(As in the previous section, we often write $\ve{X}$ instead of $\ve{X,s}$.)

We define $X^{(\omega)}$ as the union of countably many copies $ X^{(n)}$ ($n\in\omega$) of 
$\ve{X}$, where we identify $\langle x,0\rangle\in X^{(n+1)}$ with $\langle x,s(x)\rangle\in X^{(n)}$. By abuse of notation
we still write $ X^{(n)}$ for the images of these subspaces in $X^{(\omega)}$ (they are thus not disjoint).
The map $ s^{(\omega)}$ is defined as being equal to $\wt{s}$ on each $ X^{(n)}$, it is a slicer for $X^{(\omega)}$.
We denote by $\de{X}^{(n)}$ the copies of $X$ given by the images of $X\times\{0\}\subset  X^{(n)}$ 
(identified with $\de{X}\subset X^{(n-1)}$ if $n\ge 1$) in the resulting space.
A canonical sequence $(X^{(\omega)})_\alpha$ follows naturally given one for $X$.
If $X=\LL_{\ge 0}$ with the slicer given by the identity, this construction yields exactly
$(S^{\omega,h,r} - \Delta_\omega)$ with $h(n)=0$ and $r(n)=\uparrow$ for each $n$.

\begin{lemma}
   In $X^{(\omega)}$, $\de{X}^{(m)}\not\preceq\de{X}^{(n)}$ if $n < m$.
\end{lemma}
\proof
   Let $n\in\omega$.
   Define $g:X^{(\omega)}\to\LL_{\ge 0}$ as
   being identically $0$ on $X^{(k)}$ for $k<n$, 
   equal to $\langle x,y\rangle\mapsto y$ on $X^{(n)}$, and to $\wt{s}$ on $X^{(\ell)}$ for $\ell>n$.
   Then $g$ is continuous, unbounded on each $\de{X}^{(\ell)}$ for $\ell>n$ and identically $0$ on $\de{X}^{(k)}$ for $k\le n$.
\endproof

Now, set $\Omega(X)$ to be the disjoint union of $X^{(\omega)}$ and a copy $L$ of $\LL_{\ge 0}$. 
The topology on $X^{(\omega)}$ 
is the one described above. Given $x<y$ in $\LL_{\ge 0}$ and $n\in\omega$, set
$$
      V(x,y,n)= ( s^{(\omega)})^{-1}\bigl( (x,y)\bigr) - \cup_{k<n} X^{(k)}.
$$
The neighborhoods of $z\in L$ are then given by the union of 
$V(x,y,n)$ with $(x,y)\subset L$, for $x<z<y$ and 
$n\in\omega$. 
In words: points with the same slicer value $z$ in the $X^{(n)}$ accumulate to $z\in L$ when $n$ grows.
Again, a canonical sequence $\Omega(X)_\alpha$ follows automatically from one for $X$.
By definition the induced topology on $L$ is the usual one.
It is easy to check that Lemmas \ref{lemmasigmaA}--\ref{lemmasigmaD} also hold for $\Omega(X)$,
in fact $\Omega(X)_\alpha$ and $\Sigma(X)_\alpha$ are essentially homeomorphic (the only difference being
at the preimage of $0$ under the slicer).
Given a subset $A$ of $L$, the space $\Omega(X,A)\subset\Omega(X)$ is defined as the subspace $X^{(\omega)}\sqcup A$.
When the choice of the slicer $s$ has some importance, we write more accurately $\Omega(X,A,s)$ or $\Omega(X,s)$.

\begin{lemma}   \label{lemma:Omega1}
    Let $X$ by a Type I space possessing a slicer, such that $B_\alpha(X)\not=\varnothing$ for a stationary set of $\alpha$.
    Let $A\subset L$ be unbounded.
    Then $A$ is a direction in $\Omega(X,A)$.
\end{lemma}    
\proof
   Let $f:\Omega(X,A)\to\LL_{\ge 0}$, and set $F_0=f^{-1}([0,1])$, $F_1=f^{-1}([2,\omega_1))$.
   For $n\in\omega$ and $x\in\LL_{\ge 0}$, set $N(x,n)=\bigl( ( s^{(\omega)})^{-1}(\{x\}) - \cup_{k<n}X^{(k)}\bigr)$.
   Notice that if $\alpha\in\omega_1$, $N(\alpha,n)=B_\alpha(X)\times( X^{(\omega)} - \cup_{k<n}X^{(k)})$.
   If $F_0$ and $F_1$ are both unbounded in $A$, then for some $n$ the sets 
   $C_i=\{x\in\LL_{\ge 0}\,:\,N(x,n)\subset F_i\}$ ($i=0,1$) are both unbounded.
   Notice that if $\ell> n$ then $N(x,n)\supset \wt{s}^{-1}(\{x\})$ in $X^{(\ell)}$.
   If $x_m\in\LL_{\ge 0}$ is a sequence converging to $\alpha\in\omega_1$ and
   $N(x_m,\ell)\subset F_i$, then $N(\alpha,\ell)\subset F_i$ 
   and when $\ell>n$ $F_i\cap X^{(\ell)}$ contains $B_\alpha\times[0,\alpha]$.
   Thus, the sets of $\alpha$ such that either $B_\alpha(X)=\varnothing$ or $\alpha\in C_i\subset L$ 
   is club for $i=0,1$, so their intersection must contain some 
   $\alpha$ for which $\varnothing\not=N(\alpha,n)\subset C_i$ for  
   both $i$, which is impossible.
\endproof 

We will show that a direction in $\Omega(X,A)$ is forced to be inside $A$
in good cases. Before, a small technical lemma.

\begin{lemma}\label{lemma:Omega1bis}
  Let $f:\ve{X}\to\LL_{\ge 0}$ satisfy $f(x)\le \wt{s}(x)$ for all $x\in \ve{X}$ and $n\in\omega$ be given.
  Then there is $g:\Omega(X,A)\to\LL_{\ge 0}$ such that $g\restrict{X^{(n)}}=f$ and $g\restrict{A}=id$.
\end{lemma}
\proof
  Set $g\restrict{A}=id$ and
  define $g\restrict{X^{(k)}}$ as follows. 
  Let $\langle x,y\rangle^{(k)}$ denote the point $\langle x,y\rangle\in X^{(k)}$.
  If $k>n+1$, set $g(\langle x,y\rangle^{(k)})=\wt{s}(\langle x,y\rangle)$.
  If $k<n-1$, set $g(\langle x,y\rangle^{(k)})=0$. Then,
  set $g(\langle x,y\rangle^{(n+1)})=\max\{f(\langle x,s(x)\rangle),y\}$,
  $g(\langle x,y\rangle^{(n)})=f(\langle x,y\rangle)$ and finally 
  $g(\langle x,y\rangle^{(n-1)})=\min\{f(\langle x,0\rangle),y\}$,
  $g$ is continuous and satisfies the required properties.
\endproof

\begin{lemma}   \label{lemma:Omega2}
    Let $X$ be a
    type I space which possesses a slicer and contains no direction.
    If $X$ is locally compact or countably tight, then a direction in $\Omega(X,A)$ 
    is contained in $A$ outside of some bounded set.
\end{lemma}  
\proof
  Let $C$ be club in $\Omega(X,A)$. Suppose that $C\cap X^{(\omega)}$ is unbounded.
  Thus, for some $n$, $C\cap X^{(n)}$ is unbounded as well. But by Lemma \ref{nodirinveX},
  there is some $h: X^{(n)}\to\LL_{\ge 0}$ witnessing that $C\cap X^{(n)}$ is not a direction.
  Notice that $f(x)=\min\{h(x),\wt{s}(x)\}$ is unbounded on $C\cap X^{(n)}$ as well
  and is bounded on some unbounded subset of it.
  Then $g:\Omega(X,A)\to\LL_{\ge 0}$ given by the previous lemma
  shows that $C$ is not a direction.
\endproof


\subsubsection{$\Xi(X,A,s)$}\label{sec:Xi}

This is similar to the construction in Section \ref{sec:exDominant} but done in $\we{X,s}$ instead of $\LL_{\ge 0}\times[0,1]$.

\begin{figure}[h]
  \begin{center}
  \epsfig{figure=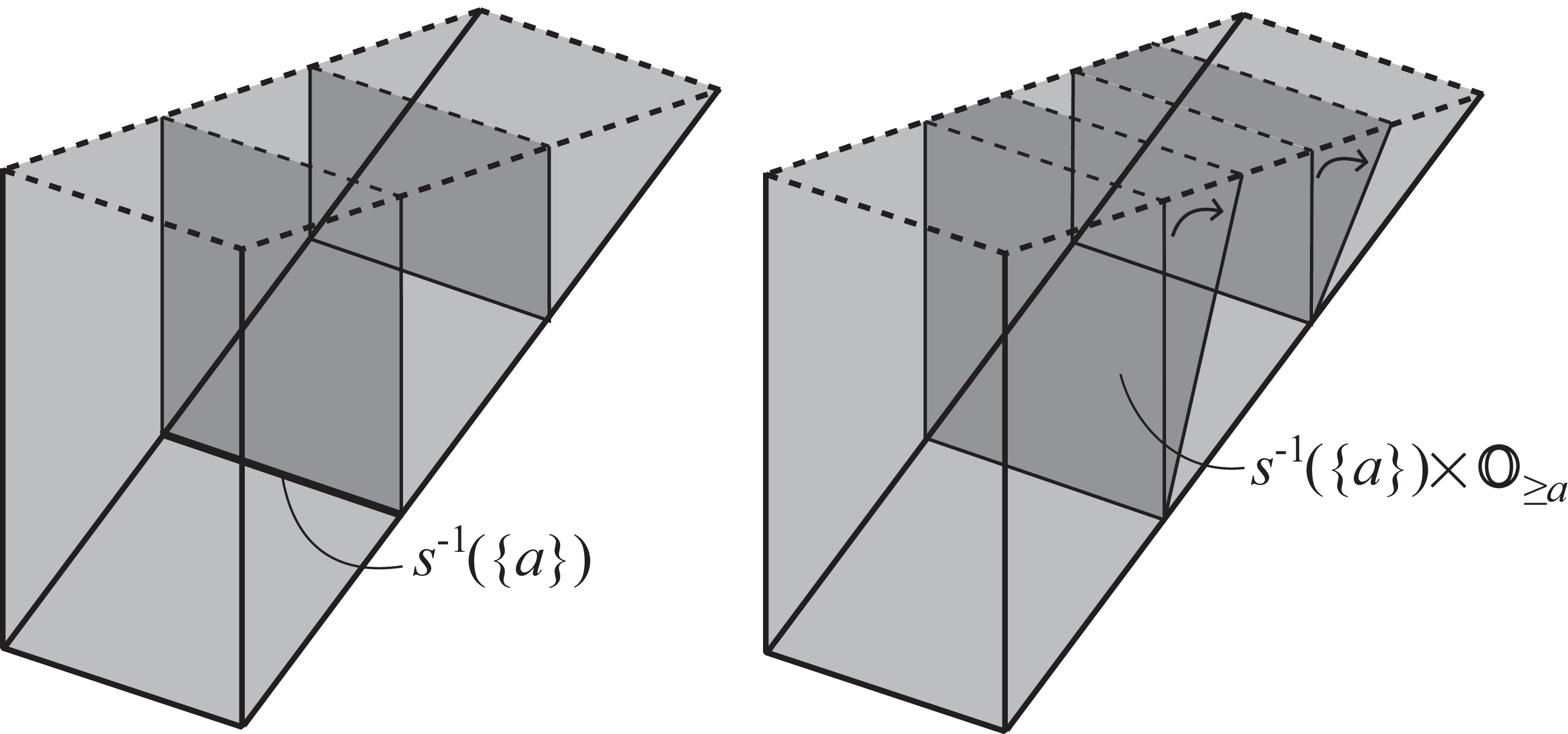, width=12.5cm}
  \caption{Inserting $s^{-1}(\{a\})\times\OO_{\ge a}$ in $\we{S}$.}\label{figureZsec7}
  \end{center}
 \end{figure}

Let $X$ be Type I with a slicer $s$. Let $A\subset\LL_{\ge 0}$ be discrete.
$\Xi(X,A,s)$ is the space obtained with the following procedure. Start with $\we{X}$.
For each $a\in A$, take out $s^{-1}(\{a\})\times[a,\omega_1)\subset\we{X}$, 
and sew in the scar a copy $E_{a}$ of $s^{-1}(\{a\})\times\OO_{\ge a}$ exactly as in the construction of Example \ref{exDominant}.
The orientation of the added octants (times $s^{-1}(\{a\})$) is chosen such that the arrow points `in the direction of the diagonal $\de{X}$',
as in Figure \ref{figureZsec7}. 

This can be made precise 
by specifying explicitly the neighborhoods of the points in the horizontals and diagonals of the octants sewed in, but
this is quite tedious to describe in detail. 
As in Example \ref{exDominant}, if $a_n$ is an increasing sequence in $A$ converging to some $a\in\LL_{\ge 0}$,
thus $s^{-1}(\{a_n\})\times\LL_{\ge a_n}$ accumulates to $s^{-1}(\{a\})\times\LL_{\ge a}$,
then the inserted $E_{a_n}$ also accumulate to it.
The $\Xi(X,A,s)_\alpha$ can be defined again as in Example \ref{exDominant}
to give $\Xi(X,A,s)$ a structure of Type I space. Notice that the `diagonal' $\de{X}$ has not been altered by the
construction, and we write $\de{X}\subset \Xi(X,A,s)$ by a slight abuse.
Denote by $p:\Xi(X,A,s)\to\we{X,s}$ the map which collapses the inserted $E_a=s^{-1}(\{a\})\times\OO_{\ge a}$ by projection on the horizontal
in $\OO_{\ge a}$. Then define the slicer $s_\Xi:\Xi(X,A,s)\to\LL_{\ge 0}$
by $s_\Xi=\wt{s}\circ p$.
If $E\subset X$, we can build $\Xi(E,A,s\restrict{E})$ and may consider it as a subspace of $\Xi(X,A)$.
This holds also if $E$ is bounded, as the construction uses only the continuity of the slicer $s$
and no other property.

\begin{lemma}\label{lemma:Xi0}
  The following holds.\\
  (a) If $X$ is $\omega$-bounded, then so is $\Xi(X,A,s)$.\\
  (b) If $X$ is $\omega_1$-compact and first countable, then then so is $\Xi(X,A,s)$.
\end{lemma}
\proof 
  This is very similar to Lemma \ref{lemma:UDomega_1cpct}, so we do not provide much detail.
  \\ 
  (a) Let $\{x_n\,:\,n\in\omega\}\subset\Xi(X,A,s)$ be given.
      If infinitely many $x_n$ belong to some $E_a$, there is an accumulation point since $E_a$ is $\omega$-bounded.
      Otherwise we may assume that there is at most one $x_n$ in each $E_a$.
      As in Lemma \ref{projclosed} (b) $p$ is proper and thus closed whenever $X$ is $\omega$-bounded,
      so $\{p(x_n)\,:\,n\in\omega\}$ has an accumulation point $y\in\we{X,s}$.
      By construction the $x_n$ must accumulate to some preimage of $y$ under $p$.
        \\
  (b)  
      By Lemma \ref{lemma:UDomega_1cpct} (c), $s^{-1}(\{a\})\times\LL_{\ge a}^2$ is $\omega_1$-compact and so 
      is $E_a$ (as a closed subset).
      Given $D\subset\Xi(X,A,s)$ closed discrete, then 
      either $D\cap E_a$ is uncountable for some $a$ and we are done,
      or there is an uncountable $D'\subset D$ such that each $|D\cap E_a|\le 1$ for each $a\in A$.
      In the latter case $p(D')$ is closed discrete, applying Lemma \ref{lemma:UDomega_1cpct} (c) again
      we obtain an accumulation point $y$ and as in (a) some preimage of $y$ must be an
      accumulation point of $D'$.
\endproof

The next two lemmas are the $\Xi$-versions of Lemmas \ref{lemma3d1} and \ref{lemma3d2}.
Their proofs are similar and we skip them.

\begin{lemma}\label{lemma:Xi1}
   Let $X$ be a Type I space that is either $\omega$-bounded or $\omega_1$-compact and first countable. 
   Let $A\subset\LL_{\ge 0}$ be discrete.
   If $D\subset \Xi(X,A,s)$ is club, then either it is contained in $\Xi(X_\alpha,A,s\restrict{X_\alpha})$
   for some $\alpha$, 
   or its intersection with $\de{X}$ is club. 
\end{lemma}

\begin{lemma}\label{lemma:Xi4}
   Let $X$ be a Type I $\omega_1$-compact and first countable space. Let $A\subset\LL_{\ge 0}$ be discrete.
   If $U\subset \Xi(X,A,s)$ is a large open set, it contains a copy of $\omega_1$
   located in $\Xi(X_\alpha,A,s\restrict{X_\alpha})$ for some $\alpha$.
\end{lemma}

Finally, we have the following version of Lemma \ref{lemma3d3}.

\begin{lemma}\label{lemma:Xi2}
   Let $X$ be a Type I space that is either $\omega$-bounded or $\omega_1$-compact and first countable.
   Let $A\subset\LL_{\ge 0}$ be unbounded and discrete.
   If $X$ does not possess a direction, then 
   there is no maximal direction in $\Xi(X,A,s)$.
\end{lemma}
\proof
   Let $D\subset\Xi(X,A,s)$ be a direction.
   If $D\cap\de{X,s}$ is club, as in Lemma \ref{lemma3d3} we may find a map witnessing that $D$ is not a direction,
   a contradiction.
   Thus by Lemma \ref{lemma:Xi1}, $D\subset\Xi(X_\alpha,A,s\restrict{X_\alpha})$ for some $\alpha$.
   Take $a\in A$ bigger than $\alpha$ and  denote $s^{-1}([0,a))$ by $X_a$.
   Define $f:\Xi(X,A,s)\to\LL_{\ge 0}$ as follows.
   On $\Xi(X_a,A,s\restrict{X_a})$ it takes the value $a$.
   In $E_{a}=s^{-1}(\{a\})\times\OO_{\ge a}$, the point
   $\langle x, y_1,y_2\rangle$ is sent to 
   $y_1$. Finally, a point $x$ elsewhere is sent to $s_\Xi(x)$.
   This defines a continuous map which is bounded on $D$ but unbounded
   on any `vertical' $\{y\}\times\LL_{\ge s(y)}$ whenever $s(y)>a$.
   Thus $D$ is not ufo-maximal.
\endproof




\section{Bounded and unbounded maps}\label{secgenprop3}

Let $X$ be Type I and $f:X\to\LL_{\ge 0}$ be unbounded. Is there a direction on which $f$ is unbounded~?
Any space without direction provides a trivial counterexample, 
$X=\OO-\Delta$, $f$ being the projection on the second coordinate, is another.
The question is more interesting for spaces spanned with directions.

\begin{thm}[{\bf PFA}]\label{thm3.0}
  Let $X$ 
  be a Type I countably tight locally compact space which is spanned by directions, and 
  let $f:X\to\LL_{\ge 0}$ be an unbounded map. Then 
  there is a direction in $X$
  on which $f$ is unbounded. 
\end{thm}

Notice that $\omega_1$-compactness is not included in the assumptions.

\begin{thm}\label{3d} 
    The following holds.
    \\
    (a) Under $\neg${\bf SH}, there is an $\omega_1$-compact $3$-manifold $M$ spanned by directions and an unbounded function
       $f:M\to\LL_{\ge 0}$ bounded on each direction in $M$.\\
    $\ast$(b) {\rm ($\diamondsuit^+$)} There is an $\omega$-bounded $3$-manifold as in (a).
\end{thm}

\proof[Proof of Theorem \ref{thm3.0}]
Let $f:X\to\LL_{\ge 0}$ be unbounded.
Let us suppose that there is some $\omega_1$-compact $Y\subset X$, closed in $X$, on which $f$ is unbounded. 
Then $Y$ is unbounded,
therefore the set $C(f)=C(f,X)$ introduced in Lemma \ref{lemmastar} is club.
As before, by Theorem \ref{thmomega_1compact}
there is a (closed) copy of $\omega_1$ in $C(f)$ on which $f$ is thus unbounded, and we are done.
\\
Assume now that any club of $X$ on which $f$ is unbounded possesses a closed uncountable discrete subset. 
Up to changing the canonical sequence, we can assume that $f(X_\alpha)\subset[0,\alpha)$ $\forall\alpha\in\omega_1$.
Set $\alpha_0$ to be $0$, choose by induction $x_\beta\in X$ and $\alpha_{\beta+1}\in\omega_1$ with
$f(x_\beta)\in (\alpha_\beta,\alpha_{\beta+1})$, and set $\alpha_\beta=\sup_{\gamma<\beta}\alpha_\gamma$ when
$\beta$ is limit. 
Let $Y$ be the closure in $X$ of $\{x_\beta\,:\,\beta\in\omega_1\}$. Since $f$ is unbounded on $Y$, 
there is a closed uncountable discrete subset $Z\subset Y$ which 
is thus unbounded
(by Lindel\"ofness of each $\wb{X_\alpha}$). Thus, $Z$ is club
and $f$ is unbounded on each club subset of $Z$.
Let $\gamma_0=0$. Given $\beta$, choose $y_{\beta+1}\in Z$ and $\gamma_\beta\in\omega_1$ such that
$\gamma_\beta < f(y_{\beta+1}) < \gamma_{\beta+1}$ and 
$y_\beta \in X_{\gamma_{\beta+1}} - \wb{X_{\gamma_\beta}}$.
Set then $V=\cup_{\beta\in\omega_1}V_\beta$, where
$$
  V_\beta =
  (X_{\gamma_{\beta+1}} - \wb{X_{\gamma_\beta}})\cap f^{-1}((\gamma_\beta, \gamma_{\beta+1})).
$$
Then $V$ is a large open set since it contains each $y_\beta$, and
since $X$ is spanned by directions, one must be included in $V$. By definition of the $V_\beta$, $f$ must be unbounded on it.
\endproof

For Theorem \ref{3d}, we will use the spaces $\we{X}$ defined in Section \ref{sec:UDDelta}.

\proof[Proof of Theorem \ref{3d}]
\
\\
(a) Let $S_T$ be a fattening of a Suslin tree $T$ as defined in Section \ref{sec:omega_1cpct}.
  By Lemma \ref{lemma3d2}, for any large open set $U$ of $\we{S_T}$ 
  there is some $y\in S_T$ such that $U\cap\{y\}\times\LL_{\ge 0}$ contains a club, and thus a copy of $\omega_1$. 
  Hence $\we{S_T}$ is spanned by directions. But the unbounded map $f:\we{S_T}\to\LL_{\ge 0}$ given by $f(x,y)=s(x)$ is bounded on any
  direction by Lemma \ref{lemma3d3}. If one prefers a boundaryless example, take the double of $\we{S_T}$ (i.e. take two copies
  of $M$ and identify pointwise
  the boundaries). 
\\
(b) The same construction than in (a) starting with the $\omega$-bounded surface given by Theorem \ref{AxN}.
\endproof

\begin{rem}
  We can do the same construction starting with $Q$ given by Theorem \ref{thm:betaomega}.
  However, $\we{Q,s}$ might not be spanned by directions because Lemmas \ref{lemma3d2}--\ref{lemma3d3} 
  may not hold. See Proposition \ref{prop:Ppoint} below and recall that $\beta\omega$ has P-points (consistently).
  $\we{Q,s}$ is then an $\omega$-bounded space which possesses directions but with an unbounded map 
  which is bounded on each direction. However, the disjoint union of $\LL_{\ge 0}$ and $Q$
  is a much simpler example.
\end{rem}



\section{Directions in themselves}\label{secdirthem}

\subsection{Relations with monoclinic and monolithic spaces}\label{sec:monolithic}

Following Nyikos \cite{Nyikos:1992, NyikosPC}, we say that a Type I space $X$ is {\em monoclinic} if any unbounded $0$-set 
(i.e. inverse image of $0$ for some real
valued function) contains the bones $B_\alpha(X)$ (which are assumed to be non-empty) for a club set of $\alpha$. 
Notice that the inverse image of a closed subset of $\R$ is a $0$-set. 
$X$ is {\em monolithic} if any club subset of $X$ contains the bones (non-empty again) on a club set.
(Caution: the word `monolithic' is sometimes used to denote a completely different property.)
A real valued map on a Type I space $X$ is said to be eventually constant iff it is constant
outside some $X_\alpha$.
The following was partially stated in \cite{NyikosPC}:

\begin{thm}\label{thm:monoclinetc}
  Let $X$ be a Type I space. 
  The items below are related as follows. 
  $$
  \begin{array}{cc}
  (c) \longleftrightarrow (d) \longleftrightarrow (e) & \text{ and }\quad
  \begin{CD} 
     (a)  @>>> (b) @>>> (c) & & \\ 
     & & @VVV @VVV  \\
     & & (f) @>>> (g) @>>>  (h)
  \end{CD}
  \end{array}
 $$
  Moreover, if $X$ possesses a slicer then (f) $\longleftrightarrow$  (b) and (g) $\longleftrightarrow$  (c) hold, and if $X$ is normal then 
  items (b) to (g) are equivalent.
  \\
  (a) $X$ is monolithic.\\
  (b) If $A,B\subset X$ are club, then $A\cap B$ is club.\\
  (c) $X$ is monoclinic.\\
  (d) If $A,B\subset X$ are unbounded $0$-sets, then $A\cap B$ is club.\\
  (e) Every real valued function on $X$ is eventually constant. \\
  (f) Any club set of $X$ is a direction in itself.\\
  (g) Any unbounded $0$-set of $X$ is a direction in itself.\\
  (h) $X$ is a direction in itself. 
\end{thm}
\proof
  Let $X$ be Type I. Let us first prove (c) $\leftrightarrow$ (d) $\leftrightarrow$ (e).\\
  (c) $\rightarrow$ (d) Immediate.\\
  (d) $\rightarrow$ (e) 
      (We phrased this proof so as to emphasize that we use only
      that $\R$ is regular, Lindel\"of and has countable pseudocharacter.)
      Let $f:X\to\R$. For each 
      $\alpha$ pick $x_\alpha\in X-X_\alpha$ and let 
       $z_\alpha=f(x_\alpha)$, $Z=\{z_\alpha\,:\,\alpha\in\omega_1\}$.     
      By Lindel\"ofness, some $z\in\R$ is a point of complete accumulation of $z$ (i.e. any neighborhood
      of $z$ has an intersection with $Z$ of cardinality $\aleph_1$).
      If $f$ is not eventually constant taking value $z$, choose $y_\alpha\in X-X_\alpha$
      such that $f(y_\alpha)\not= z$.
      Take a sequence $\{U_n\,:\,n\in\omega\}$ of neighborhoods of $z$ whose intersection is $\{z\}$.
      Then for some $n$, $\R -U_n$ contains uncountably many $f(y_\alpha)$.   
      By regularity there is some open $U\ni z$ such that $\wb{U}\subset U_n$.
      Then $f^{-1}(\wb{U})$ and $f^{-1}(X-U_n)$ are disjoint unbounded $0$-sets, a contradiction.\\
  (e) $\rightarrow$ (c) If $X$ is not monoclinic then for some real valued $f$, $f^{-1}(\{0\})$ is unbounded but does not contain
  all the bones, and thus $f$ is not eventually constant.\\
  We now prove the other implications.\\
  (a) $\rightarrow$ (b) $\rightarrow$ (d) and (f) $\rightarrow$ (g) $\rightarrow$ (h) are immediate. Let us
  first prove (e) $\rightarrow$ (h), which is not in the diagram but is useful for other implications.\\
  (e) $\rightarrow$ (h) Given $f:X\to\LL_{\ge 0}$ unbounded, set 
  $g(x)=\min\{f(x),1\}$, then
  for some $\beta$, $g$ is constant on $X-X_{\beta}$. Since $f$ is unbounded, $g$ must be constant on $1$.
  Thus $f^{-1}(\{0\})\subset X_\beta$ is bounded.\\
  (b) $\rightarrow$ (f)
  Since (b) is hereditary for club subsets, a club in $X$ satisfies (e). We conclude with (e) $\rightarrow$ (h).\\
  (c) $\rightarrow$ (g) We show (e) $\rightarrow$ (g). 
  Let $C=f^{-1}(\{0\})$ be unbounded for some $f:X\to\R$. By (e), there is $\alpha$ such that
  $C\supset X-X_\alpha$. But $X-X_\alpha$ is a direction in itself by (e) $\rightarrow$ (h).\\
  (f) $\rightarrow$ (b) if $X$ has a slicer.
  If $A,B\subset X$ are disjoint  club sets, the map sending $A$ to $0$ and which is equal to the restriction 
  of a slicer in $B$ is continuous on $A\cup B$, which shows that $A\cup B$ is not a direction in itself.\\
  (g) $\rightarrow$ (d) if $X$ has a slicer. As (f) $\rightarrow$ (b).\\
  Equivalence of items (b) to (g) if $X$ is normal.
  Since $X$ is regular, (g) $\rightarrow$ (c) and it is enough to show (c) $\rightarrow$ (b).
  Assume that (c) holds. Let $A,B$ be disjoint club sets of $X$.
  By normality there is a continuous $f:X\to\R$ such that $f^{-1}(\{0\})\supset A, f^{-1}(\{1\})\supset B$,
  and thus $f^{-1}(\{0\}), f^{-1}(\{1\})$ are disjoint unbounded $0$-sets, contradicting (c).
\endproof

Notice that we also have:
\begin{lemma}
  Let $X$ be a Type I space. Then (e) in Theorem \ref{thm:monoclinetc} is equivalent to:
  \\
  (e') Any map $f:X\to Y$ where $Y$ is regular Lindel\"of and has countable pseudocharacter is eventually constant.
\end{lemma}
\proof
  (e') $\rightarrow$ (e) is immediate. 
  Let $f:X\to Y$ be given and assume that (e) and thus (d) above holds.
  The proof of (d) $\rightarrow$ (e) shows that there are disjoint closed $E,F\subset Y$ 
  with unbounded preimages in $X$. Since $Y$ is regular and Lindel\"of it is normal, so there
  is a map $g:Y\to\R$ which takes the value $0$ on $E$ and $1$ on $F$.
  Then $g\circ f:X\to\R$ is a non-eventually constant map, a contradiction.
\endproof

The implication charts for regular and normal spaces are thus reduced to:
$$
  \begin{array}{ccc}
  (a) \rightarrow (b) \rightarrow (g) \rightarrow (h) & \quad & (a) \rightarrow (g) \rightarrow (h)\\
  \text{regular spaces}& &\text{normal spaces}
  \end{array}
$$

None of the other implications holds in the relevant classes.
Neither (a) $\rightarrow$ (b) nor (g) $\rightarrow$ (h) does reverse for normal spaces.
For (h) $\not\rightarrow$ (g), notice that the normal space $\LL_{\ge 0}\times[0,1]$ is a direction in itself
while
$\LL_{\ge 0}\times\{0,1\}$ is not.

\begin{ex}
  A $2$-to-$1$ closed preimage of $\omega_1$ which is not monolithic but satisfies (b) of Theorem \ref{thm:monoclinetc}.
\end{ex}
Such a space is thus locally countable, locally metrizable, locally compact, $\omega$-bounded and normal 
(even collectionwise) 
by Lemma \ref{lemma:monolithic} below.
This is a simple modification of a 1979 construction by Nyikos \cite[19.1]{FremlinPerfect}.  
Let $X=\omega_1\times\{0,1\}$ and $\pi$ be the projection on the first factor.
Denote by $\Lambda$ the set of limit ordinals in $\omega_1$, and by $\Lambda_1$ the set of limit of limits.
Points $\langle\alpha,i\rangle$ are isolated whenever $\alpha\not\in\Lambda_1$.
  The neighborhoods of points at levels in $\Lambda_1$ are defined as:
  $$ 
     \begin{array}{ccl}
       \langle\alpha,0\rangle & : & \{\langle\alpha,0\rangle\}\cup\{ \langle\gamma,i\rangle\,:\,\beta<\gamma<\alpha,\,\gamma\in\Lambda,\,i=0,1\},\\
       \langle\alpha,1\rangle & : & \{\langle\alpha,1\rangle\}\cup\{ \langle\gamma,i\rangle\,:\,\beta<\gamma<\alpha,\,\gamma\not\in\Lambda,\,i=0,1\}     
     \end{array}
  $$
  for $\beta<\alpha$. 
  It is immediate that $X$ is Hausdorff, $\omega$-bounded and locally countable and that $\pi$ is closed.
  Since $\pi^{-1}([0,\alpha])$ is compact, by Urysohn's theorem it is metrizable and $X$ is locally compact.
  \\
  Now, given two club subsets $A,B\subset X$, by Lemma \ref{convenientlemma} (b) both have non-empty intersection with the bones
  of $X$ on a club set $C\subset\omega_1$.
  But then both contain $\langle\alpha,0\rangle$ for any $\alpha\in\Lambda_1\cap C$, so $X$ satisfies (b).
  The club subset $\Lambda_1\times\{0\}$ does
  not contain any bone (it is actually a copy of $\omega_1$), so $X$ is not monolithic.

We now show that (g) $\not\rightarrow$ (b) 
for regular spaces.
\begin{ex}\label{ex:NyikosClass7}(Nyikos \cite{Nyikos:1992})
  There is a regular Type I space satisfying (g) but not (b) of Theorem \ref{thm:monoclinetc}.
\end{ex}
We take $X$ to be the positive part of a tangent bundle of $\LL_+$ of class 7 as in \cite[Section 6]{Nyikos:1992}
where
it is shown that $X$ satisfies (e) and thus (g) but is not normal.
The next lemma 
(whose proof is extracted from Lemma 9.5 and Theorem 4.22 in the same paper)
shows that $X$ does not satisfy (b).
Recall in passing
that a metrizable space is collectionwise normal.

\begin{lemma}\label{lemma:monolithic} (Nyikos \cite{Nyikos:1992}, in effect.)
   Let $X$ be a regular Type I space satisfying (b) of Theorem \ref{thm:monoclinetc}. 
   Then $X$ is $\omega_1$-compact and normal. If moreover $\wb{X_\alpha}$ is collectionwise
   normal for each $\alpha$, then so is $X$.
\end{lemma}
We refer to \cite{Nyikos:1992} for the definition of collectionwise normality.
\proof
   Given disjoint closed sets $A,B$, one of them must be bounded, say $A\subset X_\alpha$.
   By Lemma \ref{lemma0}, $\wb{X_{\alpha+1}}$ is normal, so we can separate $A$ and $B$ inside it by $U,V$.
   Then $U,V\cup (X-\wb{X_\alpha})$ separate $A$ and $B$. This shows that $X$ is normal.
   If $C$ is uncountable closed discrete, then for each $\alpha$, $\wb{X_\alpha}\not\supset C$ by Lindel\"ofness.
   Thus $C$ is club. But by discreteness we can write $C$ as the union of two disjoint clubs, a contradiction. It follows that
   $X$ is $\omega_1$-compact.
   \\
   Suppose now that $\{A_\beta\,:\,\beta\in\kappa\}$ is a discrete collection of closed subsets of $X$.
   If $\kappa>\omega_1$, for some $\alpha$ uncountably many $A_\beta$ have a non-empty intersection with $\wb{X_\alpha}$ 
   contradicting its Lindel\"ofness, so $\kappa\le\omega_1$.
   If some $A_\beta$ is unbounded then the others are bounded and as above we may separate all of them using the collectionwise normality of $\wb{X_\alpha}$
   as in the proof of normality. The same holds if $A=\cup_{\beta\in\kappa}A_\beta$ is bounded.
   If $A$ is unbounded and each $A_\beta$ is bounded, 
   then $\kappa=\omega_1$ and we may divide it in two disjoint
   sets such that the union of the $A_\beta$ for $\beta$ in each set is unbounded. 
   But then $A$ is the union of two disjoint club subsets, a contradiction.
\endproof


\subsection{In products}\label{sec:inproducts}

If $X$ is a direction in itself, then it is easy to see that its product with $\R$ and $[0,1]$ is a direction in 
itself as well. Of course, its product with a non-connected space cannot be a direction in itself.
Also, the product of two Type I spaces one of which possessing a slicer is never a direction in itself.
Note also the following.

\begin{ex}\label{ex:compactnotdir}
  A compact connected space $X$ such that $\LL_{\ge 0}\times X$ is not a direction in itself.
\end{ex}
\proof
  Set $X$ to be $\LL_{\ge 0}\sqcup\{\omega_1\}$ with the obvious topology (i.e. neighborhoods of $\omega_1$ are 
  terminal parts of $\LL_{\ge 0}$). Then
  use the map $f(x,y)=\min\{x,y\}$. 
\endproof

Recall that a subset $A\subset X$ is a {\em P-set} 
whenever the intersection of any countable family of neighborhoods of $A$ is a neighborhood of $A$.
An easy generalization of Example \ref{ex:compactnotdir} is the following,
whose proof is very similar to that of \cite[Theorem 3.1]{JuhaszVanMillWeiss:2013}. 

\begin{prop}\label{prop:Ppoint}
  Let $Y$ be a normal space containing a closed P-set. Then
  $\LL_{\ge 0}\times Y$ is not a  direction in itself.
\end{prop}
\proof
  Let $A$ be a closed P-set.
  By normality of $Y$, one chooses a sequence $\{U_\alpha\,:\,\alpha\in\omega_1\}$ of open neighborhoods of $A$
  such that $U_\alpha\supset\wb{U_\beta}$ whenever $\alpha<\beta$, and
  as in Lemma \ref{lemma0} we define $r:Y\to\LL_{\ge 0}\cup\{\omega_1\}$
  such that $r(\wb{U_\alpha}-U_{\alpha+1})\subset[\alpha,\alpha+1]$ and $r(\cap_{\alpha\in\omega_1}U_\alpha)=\{\omega_1\}$.
  Then, set $f:\LL_{\ge 0}\times Y\to\LL_{\ge 0}$ as $\langle x,y\rangle\mapsto\min\{x,r(y)\}$.
  This provides a function which is unbounded on $\LL_{\ge 0}\times Y$ but bounded on $\LL_{\ge 0}\times\{y_0\}$ for each 
  $y_0\not\in\cap_{\alpha\in\omega_1}U_\alpha$.
\endproof

On the positive side, we have the following proposition.

\begin{prop}\label{lemmadoublestar}
   Let $X=\cup_{\alpha\in\omega_1}X_\alpha$ be a Type I space which is a direction in itself, 
   and $Y$ connected and such that $\wb{X_\alpha}\times Y$ is Lindel\"of for each $\alpha\in\omega_1$.
   Then $X\times Y$ with canonical sequence $\langle X_\alpha\times Y\,:\alpha\in\omega_1\rangle$ is a direction in itself
   if either\\
   (a) $Y$ is separable, or\\
   (b) $X$ is $\omega$-bounded and $Y$ is first countable.
\end{prop}

The assumption on $\wb{X_\alpha}\times Y$ implies that $Y$ is Lindel\"of and is satisfied if either $Y$ or each $\wb{X_\alpha}$ is
productively Lindel\"of (i.e. its product with any Lindel\"of space is Lindel\"of),
for instance $\sigma$-compact spaces. 

\proof
  Given $f:X\times Y\to\LL_{\ge 0}$, set
  $$ 
  \begin{array}{l}
      A_0=\{y\in Y\,:\, f\text{ is bounded on }X\times\{y\}\},\\
      A_1=\{y\in Y\,:\, f\text{ is unbounded on }X\times\{y\}\}.
  \end{array}
  $$
  (a)
  Let $D\subset Y$ be countable and dense
  and set $\gamma=\sup_{y\in A_0\cap D}\sup f\restrict{X\times\{y\}}<\omega_1$.
  By continuity, $f$ is bounded by $\gamma$ on $X\times\wb{D\cap A_0}$. 
  By 
  Lemma \ref{pulpita} (e), since $X$ is a direction
  in itself there is a club
  $C\subset\omega_1$ such that for all $y\in A_1\cap D$ 
  and all $\alpha\in C$, $f((X-X_\alpha)\times\{y\})\subset[\alpha,\omega_1)$. 
  If $y\in\wb{D\cap A_0}\cap\wb{D\cap A_1}$, then for $\alpha \in C$ we have
  by continuity that $f(x,y)\in [\alpha,\omega_1)$ whenever $x\in X-X_\alpha$, a contradiction when $\alpha >\gamma$.
  Thus $\wb{D\cap A_0}\cap\wb{D\cap A_1}=\varnothing$, since $D$ is dense and $Y$ is connected
  it follows that one of $\wb{D\cap A_0}$, $\wb{D\cap A_1}$ is empty. If the latter is empty,
  then $f$ is bounded by $\gamma$ on $X\times Y$. If $f$ is unbounded,
  $\wb{D\cap A_1} = Y$, and by continuity $f((X-X_\alpha)\times Y)\subset [\alpha,\omega_1)$
  whenever $\alpha\in C$, which shows that $X\times Y$ is a direction in itself.\\
  (b) As in (a), one sees that $\wb{A_0}=A_0$ and $\wb{A_1}=A_1$, using the
  fact that if $y\in\wb{A_0}$ there is a countable $B\subset A_0$ with $y\in\wb{B}$.
  Thus either $A_0=\varnothing$ or $A_1=\varnothing$.
  In the latter case,
  for each $y\in A_0=Y$ there is $\gamma(y)$ such that
  $f^{-1}([0,\gamma(y)))\supset X\times\{y\}$. By lemma \ref{lemma:lemma2.2BGGG}, there is an open
  $V(y)\ni y$ with $f^{-1}([0,\gamma(y))\supset X\times V_y$. A countable family of $V_y$ covers $Y$,
  hence the supremum of the corresponding $\gamma(y)$ is a bound for $f$ on $X\times Y$.
  If $f$ is unbounded then $A_1=Y$. Since $X$ is a direction in itself, given $\gamma\in\omega_1$
  for each $y\in Y$ there is $\gamma(y)$ such that
  $f^{-1}((\gamma,\omega_1))\supset (X-X_\gamma)\times\{y\}$, and we conclude similarly.
\endproof


\subsection{No direction in itself}\label{sec:nodirinitself}

Does a direction in some space always contain a direction in itself~? 
For $\omega_1$-compact countably tight locally compact spaces, the answer is yes under {\bf PFA}
by Theorem \ref{thmomega_1compact}, but consistently `no' for a similar class of spaces.

\begin{thm}\label{thm:dirinitself}
  The following holds.\\
  (a) There is an $\omega$-bounded Type I space which is a direction in itself and which contains
      a club subset containing no direction in itself.  \\
  (b) {\rm ($\neg${\bf SH})} There is a Type I locally metrizable locally separable $\omega_1$-compact
  space which is a direction in itself and which contains
      a club subset containing no direction in itself.  \\
  $\ast$(c) {\rm ($\diamondsuit^+$)} There is an $\omega$-bounded space as in (b).
\end{thm}

  In (b) our example is unfortunately not locally compact.

\proof[Proof of Theorem \ref{thm:dirinitself}]
  We will use $\Sigma(X)$ for spaces $X$ without direction in the relevant three classes.
  Since $X\times\{0\}\subset\Sigma(X)$ contains no direction in itself, we
  just need to show that $\Sigma(X)$ is a direction in itself.
  \\
  (a) Let $Q$ be given by Theorem \ref{thm:betaomega}. 
      Suppose that $f:\Sigma(Q)\to\LL_{\ge 0}$ is given with $f^{-1}(\{0\})$ unbounded.
      If $f^{-1}(\{0\})\cap Q\times[0,1)$ is bounded then by density $f$ is bounded on all of $\Sigma(Q)$.
      We may thus assume that $f^{-1}(\{0\})\cap Q\times[0,1)$ is unbounded. 
      By Lemma \ref{lemma:lemma2.2BGGG} $f^{-1}(\{0\})\cap Q\times\{a\}$ is club for some $a\in[0,1)$.
      Let $Z=\{z\in Q\,:\,\langle z,a\rangle\in f^{-1}(\{0\})\cap Q\times\{a\}\}$.
      By Lemma \ref{lemma:lemma2.2BGGG} again $f^{-1}([0,1))\cap Q\times[0,1)$
      contains $Z\times (a_1,b_1)$ for some $a_1<a<b_1$.
      By continuity $f^{-1}([0,1])\supset Z\times[a_1,b_1]$.
      Proceeding by induction, we find $a_\alpha,b_\alpha$ indexed by ordinals
      with $a_\alpha <a_\beta <a < b_\beta <b_\alpha$ whenever $\beta<\alpha$
      such that $f^{-1}([0,\alpha))\supset Z\times(a_\alpha,b_\alpha)$.
      Then $b_\alpha=1$ for some $\alpha<\omega_1$, and thus $f^{-1}([0,\alpha+1))$ contains a club subset
      of the core of $\Sigma(Q)$ (namely, the set of $\alpha$ such that $B_\alpha(Q)\cap Z\not=\varnothing$).
      By Lemma \ref{lemmadt3}, $f^{-1}([0,\alpha+1))\supset (Q-Q_\gamma)\times(c_0,1)\sqcup L$
      for some $c_0$.
      By induction again, we find a sequence $c_\beta$ that `goes down', and 
      for some $\beta,\gamma\in\omega_1$, $f^{-1}([0,\beta+1))\supset \Sigma(Q)-\Sigma(Q)_\gamma$.
      Since $f$ is bounded on $\Sigma(Q)_\gamma$, $f$ is bounded on all of $\Sigma(Q)$ which is thus a direction in itself.
  
  (b) Let $T$ be a Suslin tree with slicer $s:T\to\omega_1\subset\LL_{\ge 0}$ given by the height in $T$. 
  $T$ is regular, locally second countable and $\omega_1$-compact, thus $\Sigma(T)$ is 
  locally metrizable and $\omega_1$-compact by Corollary \ref{corsigmaC} and Lemma \ref{lemmasigmaD}.
  It is moreover locally separable since $T\times[0,1)$ is such.
  We now show with a sequence of lemmas that       
  $\Sigma(T)$ is a direction in itself. 
  The idea is similar to that in (a): go `up and down' until you fill up
  the entire space, but since Lemma \ref{lemma:lemma2.2BGGG} may fail
  (at least the part about large open sets), we must work harder and use some properties of Suslin trees.
  
  \begin{lemma}\label{lemmadtSuslin}
    Let $T$ be a Suslin tree,
    and $f:T\times\R\to\LL_{\ge 0}$. Let $a,b\in\R$, then $f$ is bounded on $T\times\{a\}$ iff $f$ is bounded on $T\times\{b\}$.
  \end{lemma}
  This lemma does not hold if $T$ is an $\R$-special tree, see \cite{meszigues+Nyikos}.
  \proof
    This will follow from the fact that the road space of a Suslin tree is not contractible.
    We can assume that $a=0$, $b=1$, and that $f$ is constant on $0$ on $T\times\{0\}$.
    Suppose that $f$ is unbounded on $T\times\{1\}$.
    By Lemma \ref{lemmastar}, the subset 
    $T'=\{u\in T\,:\,u\in\text{Lev}_\alpha(T)\text{ and }f(u,1)=\alpha\}$ 
    is club in $T$, thus $T'$ with the induced order is a Suslin tree.
    We now look at the function only in $T'\times[0,1]$.
    Consider the road space
    $R_{T'}$ of $T'$.
    There is an obvious way to define a slicer $s:R_{T'}\to\LL_{\ge 0}$ that corresponds to the height in ${T'}\subset R_{T'}$, and
    to extend the domain of $f$ to $R_{T'}\times[0,1]$.   
    For $\langle u,t\rangle\in R_{T'}\times[0,1]$ set $h(u,t)\in R_{T'}$ to be the element $v\le u$ such that $s(v)=f(u,t)$,
    then $h:R_{T'}\times[0,1]\to R_{T'}$ is continuous, $h(u,1)$ is the identity and $h(u,0)$ is constant on $0$.
    But such a map does not exist because the road space of a Suslin tree is not contractible,
    by Corollary 6.3 in \cite{meszigues+Nyikos}.
  \endproof

  \begin{lemma}\label{lemmadt2}
    Let $T$ be a Suslin tree,
    and $f:\Sigma(T)\to\LL_{\ge 0}$ be such that for some $\beta\in\omega_1$, 
    $f^{-1}([0,\beta])$ intersected with some horizontal is unbounded.
    Then $f$ is bounded on $\Sigma(T)$.
  \end{lemma}
  \proof
    Suppose that
    for some $a\in[0,1)$ and an unbounded $T'\subset T$, $T'\times\{a\}=f^{-1}([0,\beta])\cap (T\times\{a\})$.
    But $T'$ is thus a Suslin subtree of $T$, and by the previous lemma, $f$ is bounded on $T'\times\{y\}$ for each $y\in\R$.
    Take a countable dense subset of $[0,1)$, then the supremum of the bounds for 
    the corresponding horizontals gives a bound for $f$ on $T'\times [0,1)$, and 
    $f$ is thus bounded on an unbounded subset of the core $L$.
    Since $L$ is a direction in itself, it follows that $f$ is bounded on it, say by $\alpha$.
    By Lemma \ref{lemmadt3} this implies that $f^{-1}([0,\alpha+1))$ contains the terminal part of 
    some horizontal, and thus $f$ is bounded on this horizontal. The previous lemma enables to conclude,
    using again a countable dense subset of $[0,1)$, that $f$ is bounded on $T\times[0,1)$. 
  \endproof

  We now finish the proof of (b).
  Let $f:\Sigma(T)\to\LL_{\ge 0}$ be such that $f^{-1}(\{0\})$ is unbounded. 
  For uncountably many $\alpha$ there are $a_\alpha,b_\alpha\in\R$ and $x_\alpha$ at height $\alpha$ in $T$
  such that $\{x_\alpha\}\times(a_\alpha,b_\alpha) \subset f^{-1}([0,1))$.
  Thus for some $a,b\in\Q$ there are uncountably many $\alpha$ with $\{x_\alpha\}\times(a,b)\subset f^{-1}([0,1))$.
  This implies that $f^{-1}([0,1])$ intersects the horizontals at height between $a$ and $b$ unboundedly.
  By Lemma \ref{lemmadt2} this implies that $f$ is bounded on $\Sigma(T)$.
  
  (c) 
  Let $N$ 
  be given by Theorem \ref{AxN}.
  As in (b), $\Sigma(N)$ is locally metrizable locally separable and $\omega$-bounded by
  Corollary \ref{corsigmaC} and Lemma \ref{lemmasigmaD}.
  The proof that it is a direction in itself is verbatim as in (a).
\endproof

Let us make some remarks about $\Sigma(X)$.
If one starts with $S_T$ (defined in Section \ref{sec:omega_1cpct}) in place of $T$, then $S_T\times[0,1)$ is a $3$-manifold, but $\Sigma(S_T)$ is not.
It is possible to obtain a manifold if one starts 
with an $\omega$-bounded surface $S$ that possesses a slicer with good properties. 

\begin{defi}
   Let $S$ be a Type I surface. 
   Then $s:S\to\LL_{\ge 0}$ is a gentle slicer if $s$ is a slicer and
   for each $\alpha\in\omega_1$ there is a continuous $\theta:S_\alpha\to\mathbb{S}^1$
   such that $x \mapsto \langle \theta(x), s(x) \rangle$ is a homeomorphism
   from $S-s^{-1}(\{0\})$ to $\mathbb{S}^1\times(0,\alpha)$.
\end{defi}

\begin{figure}[h]
  \begin{center}
  \epsfig{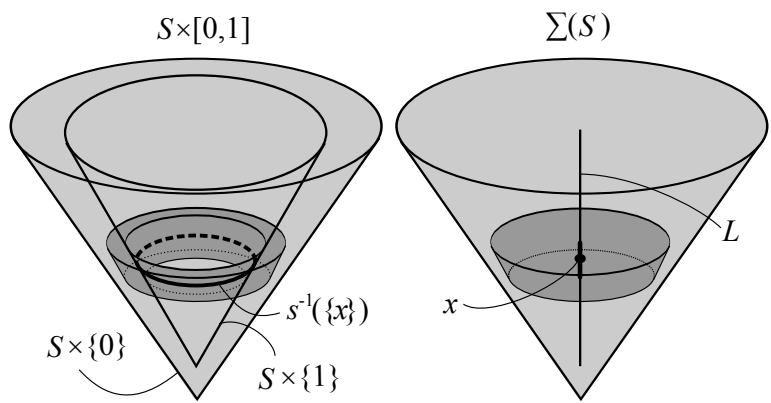}
  \caption{Neighborhoods in $\Sigma(S)$ for a $S$ a surface with a gentle slicer.}\label{figureNsec3}
  \end{center}
\end{figure}

This ensures that the points in the added copy of $\LL_{\ge 0}$ 
do possess neighborhoods homeomorphic to the $3$-ball, as suggested on Figure \ref{figureNsec3},
which implies the following corollary:

\begin{cor*} {\rm($\diamondsuit^+$)} If the surface given by Theorem \ref{AxN} possesses a gentle slicer, there is an $\omega$-bounded
             $3$-manifold which is a direction in itself and which possesses a club
             subset that is not a direction in itself.
\end{cor*}

\section{Discrete directions}\label{sec:directiondiscrete}

If a club subset of a Type I space is discrete, then it contains a club subset disjoint from some skeleton of $X$, as easily seen.
So, can a direction be disjoint from the skeleton (and thus contain a discrete direction)~? 
In general the answer is positive but not for
$\omega_1$-compact spaces: 

\begin{lemma}\label{lemma:firstsec7}
  Let $X$ be $\omega_1$-compact and Type I, and $D\subset X$ be club. Then,
  $D\cap Sk(X)$ is a direction in $X$ iff $D$ is a direction in $X$.
\end{lemma}

\proof 
  Let $D$ be a direction. By Lemma \ref{convenientlemma}
  $D\cap Sk(X)$ is club and thus a direction. Conversely,
  let $f:X\to\LL_{\ge 0}$ be continuous with $f\restrict{D}$ unbounded.
  If $f\restrict{D\cap Sk(X)}$ is bounded by $\alpha$, say, then $f^{-1}([\alpha+1,\omega_1))\cap D$ is a club
  that misses all the bones, which is impossible since $X$ is $\omega_1$-compact, so 
  $f\restrict{D\cap Sk(X)}$ is unbounded. So $f^{-1}(\{0\})\cap D\cap Sk(X)$ is bounded, say by $\beta$.
  Thus, 
  $f^{-1}(\{0\})\cap D$ must be bounded since it misses all the bones above $\beta+1$.
\endproof

If the space is not $\omega_1$-compact, both implications may fail. 
The space given by two copies of $\omega_1$ with respectively the discrete topology and the order topology
is not a direction in itself but its skeleton is.
The other implication may fail badly, as shown on the simple Example \ref{ex:spanneddiscrete} below.

While directions in a space can be discrete, and the space itself be highly non-$\omega_1$-compact (see Theorem \ref{thmdirectiondiscrete}),
the presence of a direction forces some kind of `indiscreteness', at least for regular spaces.
Recall that 
a collection of subsets of a space $X$ is said {\em discrete} iff any point of $X$ has
a neighborhood intersecting at most one element of the collection. This implies that given any subcollection,
the union of the closures of its members is closed.
If the sets are singletons, the union is thus a closed discrete set.

\begin{lemma}\label{lemmaTychonov2}
   Let $D$ be a direction in a regular Type I space $X$. Let $U$ be a neighborhood of $D$, then the subset of
   $\alpha$ such that
   $\wb{U}\cap B_\alpha(X)\not=\varnothing$ is stationary. 
\end{lemma}
\proof 
  If there is a neighborhood $U$ of $D$ whose closure misses a club set of bones, up to
  changing the canonical sequence we can assume that it misses all the bones.
  We can also assume that $U_\alpha\cap D$ is non-empty for each $\alpha$, where
  $U_\alpha = U\cap(X_{\alpha+1}-\wb{X_\alpha})$. Notice that the collection of the $U_\alpha$ is discrete.
  Choose
  $x_\alpha\in D\cap U_\alpha$. 
  Partition $\omega_1$ into two disjoint uncountable subsets $E_0,E_1$, then
  $D_i=\{x_\alpha\,:\,\alpha\in E_i\}\subset D$ are club for $i=0,1$.
  By Lemma \ref{lemma0} $X$ is Tychonoff, thus
  for $\alpha\in E_1$ there is $f_\alpha:\wb{U_\alpha}\to\LL_{\ge 0}$ that takes the value $0$ on $\wb{U_\alpha}-U_\alpha$ and $\alpha$ on $x_\alpha$.
  Define then $f:X\to[0,\alpha]\subset\LL_{\ge 0}$ as taking the value $0$ 
  outside of $\cup_{\alpha\in E_1}\wb{U_\alpha}$, and being equal to $f_\alpha$ on $\wb{U_\alpha}$.
  Then $f$ is continuous and unbounded on $D$ but $f^{-1}(\{0\})\supset D_0$ is unbounded as well, so $C$ is not a direction.
\endproof

A space is {\em (strongly) collectionwise Hausdorff} iff
any closed discrete set $D$ can be expanded to a disjoint (resp. discrete) collection of open sets.
That is, for each $x\in D$ there is an open $U_x\ni x$ such that the collection of the $U_x$ is disjoint (resp. discrete).
We abbreviate collectionwise Hausdorff by cwH. It is well known that a normal cwH space is strongly cwH.
The following is very similar to the previous lemma:

\begin{lemma}\label{lemmaTychonovI}
   Let $X$ be Type I, regular and strongly cwH.
   Then a direction in $X$ cannot be discrete.
\end{lemma}
\proof
   Let $C\subset X$ be closed discrete (and thus unbounded). 
   We may list $C=\{x_\alpha\,:\,\alpha\in\omega_1\}$.
   Expand $C$ to a collection of discrete open sets $U_\alpha\ni x_\alpha$.
   Then proceed as in Lemma \ref{lemmaTychonov2} to show that $C$ is not a direction.
\endproof

Example \ref{ex:spanneddiscrete} shows that `strongly' cannot be omitted in the statement. 
The first part of the following lemma is
an easy corollary.

\begin{lemma}\label{lemma:cwH}
  The following holds.\\
  (a)
  Let $X$ be Type I, regular, strongly cwH and spanned by directions. Then $X$ is $\omega_1$-compact.\\
  (b) Let $X$ be Type I, regular, $\omega_1$-compact and such that each $\wb{X_\alpha}$ is cwH.
       Then $X$ is strongly cwH.
\end{lemma}
\proof
  \ \\
  (a)
  If $X$ is not $\omega_1$-compact, there is a club discrete subset which expands to a discrete collection of open sets.
  The union of these open sets must contain a direction and thus a discrete one, which is impossible by the previous lemma.\\
  (b) A closed discrete subset must be countable and thus contained in some $\wb{X_\alpha}$.
       Since $\wb{X_{\alpha+1}}$ is normal, it is strongly cwH and we may expand the discrete set in it.
\endproof

\begin{cor}{\rm ({\bf PFA})}
  Let $X$ be Type I, regular, strongly cwH, locally compact, countably tight and spanned by directions. 
  Then any club subset of $X$ contains a copy of $\omega_1$.
\end{cor}
\proof
 By Theorem \ref{thmomega_1compact} and the previous lemma.
\endproof

\begin{ex}\label{ex:spanneddiscrete}
  A locally metrizable locally separable cwH Type I space which is a direction in itself and is spanned by discrete directions.
\end{ex}
\proof
Take $X=(\LL_{\ge 0}\times\R) - (\omega_1\times\Q)$.
The proof that $X$ is a direction in itself is similar to that of Theorem \ref{thm:dirinitself} (a).
Let $f:X\to\LL_{\ge 0}$ be given such that $f^{-1}(\{0\})$ is club.
For $r\in\R$ set $E_r=f^{-1}(\{0\})\cap(\LL_{\ge 0}\times\{r\})$.
If $f^{-1}(\{0\})\cap\LL_{\ge 0}\times(\R-\Q)$ is unbounded, by Lemma \ref{lemma:lemma2.2BGGG}
$E_r$ is club for some $r\in\R-\Q$, else 
$E_r$ is club for some $r\in\Q$. In both cases there are
$p_1<r<q_1$ such that $X\cap\{x\}\times(p_1,q_1)\subset f^{-1}([0,1))$ for an unbounded set of $x$.
Thus $f^{-1}([0,1])\supset[\beta_1,\omega_1)\times [p_1,q_1]\cap X$ for some $\beta_1\in\omega_1$.
Proceeding by induction as in Theorem \ref{thm:dirinitself} (a), we find $\alpha\in\omega_1$
such that $f^{-1}([0,\alpha])\supset[\beta_\alpha,\omega_1)\times\R\cap X$
and $f$ is thus bounded on $X$.
The proof that $X$ is spanned with discrete directions is very similar and we skip it.
To see that it is cwH, notice that if $D\subset X$ is closed discrete 
then by Lemma \ref{lemma:lemma2.2BGGG} again its intersection with $\LL_{\ge 0}\times(\R -\Q)$ must be bounded, say by $\alpha$.
A metrizable space is cwH, so inside $X_\alpha$ we may separate the points of $D$.
The remaining points of $D$ lie in $\cup_{\beta\in\omega_1}(\beta,\beta+1)\times\R$, and we may separate $D$ in 
each piece.
\endproof

It is possible to go beyond this example and find spaces (in various classes) such that any direction contains a discrete one.

\begin{thm}\label{thmdirectiondiscrete}
   The following holds: \\
   (a) There is a Type I locally metrizable locally separable space that possesses a discrete direction, and 
      such that any club in $X$ contains a club discrete subset. \\    
   (b) Under $\diamondsuit^*$, there is a cwH space as in (a).\\
   (c) {\rm ($\neg${\bf SH})} There is a Type I locally metrizable locally separable
       space possessing a discrete direction and a dense $\omega_1$-compact subspace, 
      and
      such that any direction in the space contains a club discrete subset. \\
   $\ast$(d) {\rm ($\diamondsuit^+$)} There is a Type I space as in (c), such that the
      dense $\omega_1$-compact subspace is a manifold which is a countable union of $\omega$-bounded submanifolds. 
      If the surface given by Theorem
      \ref{AxN} possesses a gentle slicer, the entire space is actually a $3$-manifold.\\
\end{thm}

For (a) and (b) a subset of $\Sigma(X)$ for well chosen $X$ would suffice, but
something else
is needed for (c) and (d) since any $\Sigma(X)$-variant would yield a direction in itself, hence any club
in it is a direction, which would contradict the existence of a club $\omega_1$-compact subset of the space
by Lemma \ref{lemma:firstsec7}.
We shall use $\Omega(X,A)$, which
gives also the result for (a) and (b).

\proof[Proof of Theorem \ref{thmdirectiondiscrete}]
   \ \\
   (a) Let $T$ be an $\R$-special Aronzsajn tree with a slicer given by the height in $T$, and $A$ the subset of successor ordinals in $L$. 
   Any club subset of $T$ contains a club discrete subset (see for instance \cite[Theorem 2]{DevlinNoteBaumgartner}).
   Lemmas \ref{lemma:Omega1}--\ref{lemma:Omega2}
   show that $A$ is a direction in $\Omega(T,A)$, 
   we just need to show that any club subset in $C\subset \Omega(T,A)$ contains a closed uncountable discrete subset.
   If $C\cap A$ is unbounded, we are done.
   If not, then $C\cap T^{(n)}$ is club for some $n$,
   and we apply Lemma \ref{closeddiscreteweX} (recall that $T^{(n)}$ is a copy of $\ve{T}$).
   \\
   (b) This is a consequence of (a) and a result of Devlin and Shelah \cite{DevlinShelah} who used $\diamondsuit^*$
   to build an $\R$-special cwH tree $T$. 
   We check that this implies that $\Omega(T,A)$ is also cwH, so let $C\subset \Omega(T,A)$ be closed
   discrete. For each $n$, by Lemma \ref{closeddiscreteweX} the projection $\pi(C\cap T^{(n)})$
   on the first factor of $T^{(n)}$
   is discrete. Since $T$ is cwH, $\pi(C\cap T^{(n)})$ avoids a club set $E_n\subset\omega_1$ of
   bones. This follows easily by the pressing down lemma and holds in any Type I cwH space $Z$ such that $Z_\alpha$
   is second countable (or just ccc) for each $\alpha$. Then $C$ avoids the bones in $E=\cap_{n\in\omega}E_n$, up to changing
   the canonical sequence assume $E$ to be $\omega_1$.
   Since $\Omega(T,A)_{\alpha+1}-\wb{\Omega(T,A)_{\alpha}}$ is metrizable, it is cwH and we may
   separate $C$ in each piece, yielding the result.
   \\
   (c) Start with $T$ a Suslin tree and $A$ as in (a). 
      Then, by Lemma \ref{lemma:UDomega_1cpct}, $T^{(\omega)}$ is $\omega_1$-compact, 
      and the result follows.\\
   (d) Take $N$ given Theorem \ref{AxN} and $A=\LL_{\ge 0}-\omega_1$ and proceed as in (c).
   If $N$ possesses a gentle slicer, then $\Omega(S,A)$
   is a $3$-manifold (see the discussion after the proof of Theorem
   \ref{thm:dirinitself}).
\endproof

Some remarks.
If {\bf MA + $\neg$CH} holds, then any Aronszajn tree is special by a result of Baumgartner, and
this easily implies that $T$ (and thus $\Omega(T,A)$) is not cwH.

\begin{prob} 
   Does Point (b) of Theorem \ref{thmdirectiondiscrete} hold in {\bf ZFC}?
\end{prob}

Notice that none of our examples is spanned with directions.

\begin{prob}\label{prob:directiondiscrete2}
   Is there a (first countable) Type I space spanned by directions, such that any direction in it contains a discrete club subset?
\end{prob} 

We believe that a construction mixing $\Omega(X,A)$ and Example \ref{exDominant} might work for an appropriate $X$, but
it remains to be checked.
We notice that Example \ref{ex:NyikosClass7} gives a related (though much weaker) result in {\bf ZFC}:
a Type I surface which is a direction in itself, such that no unbounded subset of $X$ is 
countably compact. Indeed, it is quasi-developable, and hence any countably compact subset is compact and thus
bounded
  (see Lemma 6.8 and the discussion p. 162 in \cite{Nyikos:1992}).





\section{More on the ufo $\langle{\mathfrak{D}}_X,\preceq\rangle$}\label{sec:ufo}

\subsection{Ufo-maximal and ufo-minimal directions}\label{secufomaxmin}

This section treats the question:
Given a Type I space $X$ spanned by directions, are there ufo-maximal or minimal directions in 
$\langle{\mathfrak{D}}_X,\preceq\rangle$~? 
Notice that Theorem \ref{thmomega_1compact} implies trivially (under {\bf PFA}) that an $\omega_1$-compact countably tight locally compact space
is spanned by directions, and even better: the direction is located inside the club.

We are going to use the $S^{\alpha,h,r}$ defined in Section \ref{sec:Obtainingalpha} in various examples.
Our first result is the following. 

\begin{thm}\label{thmufonominnomax}
  There is an $\omega_1$-compact surface $S$ spanned with directions 
  such that $\langle\mathfrak{D}_S,\preceq\rangle$ is isomorphic to $\langle\Z,\leq\rangle$.
  In particular, there is neither an ufo-min nor an ufo-max direction in $S$.
\end{thm}

\proof
  Take $T=S^{\omega,h,r}-\Delta_\omega$, with $h:\omega\to\omega_1$ the inclusion and
  $r(n)=\uparrow$ for all $n\in\omega$. Then $T$ is $\omega_1$-compact, because any uncountable discrete set
  has to be unbounded, and any club $C\subset T$ must intersect one of the (closed) $O_n$ in a club way. 
  By Lemma \ref{corchiant0}, $T$ is spanned by directions.
  Since it is a particular case of Example \ref{ex3}, 
  $\langle{\mathfrak{D}}_T,\preceq\rangle\simeq\langle\omega,\le\rangle$.
  Taking $r(n)=\downarrow$ for all $n$ yields $T'$ with $\langle{\mathfrak{D}}_{T'},\preceq\rangle\simeq\langle\omega^*,\le\rangle$. 
  Finally, glue
  $T$ and $T'$ along their remaining free side to obtain the desired $S$.
\endproof

In the $\omega$-bounded framework, all countable chains do possess minimal and maximal elements for sequential spaces
(in other words, $\langle\mathfrak{D}_M,\preceq\rangle$ is up- and downwards $\omega_1$-closed):

\begin{lemma}
  \label{ufo} 
  Let $X=\cup_{\alpha\in\omega_1}X_\alpha$ be Type I, $\omega$-bounded and sequential, $C_i\subset X$ $i\in\mathbb{Z}$ be directions with
  $C_i\preceq C_{i+1}$. Then there are directions $D,D'\subset X$ such that $D\preceq C_i\preceq D'$ for all 
  $i\in\mathbb{Z}$ 
\end{lemma}

\proof
If $\mathcal{A}=\{D_j\,:\,j\in \omega\}$ is a countable set of directions, we let $D_{\mathcal{A}}$ be the set of
accumulation points of sequences $x_{j_n}$, $j_n$ a strictly increasing sequence in $\omega$, 
with $x_{j_n}\in D_{j_n}$ for each $n\in\omega$. 
$D_{\mathcal{A}}$ is closed since $X$ is sequential, and
if $X$ is $\omega$-bounded, it is unbounded.
Suppose that $D_{\mathcal{A}}$ is not a direction, there is thus a function $f:X\to\LL_{\ge 0}$ unbounded on $D_{\mathcal{A}}$ with
$f^{-1}(\{0\})\cap D_{\mathcal{A}}$ unbounded. 
For each $\alpha$ there are therefore $k(\alpha)$ and $\ell(\alpha)$ in $\omega$, 
and
$x_\alpha\in D_{k(\alpha)} -  X_\alpha$, $y_\alpha\in D_{\ell(\alpha)} -  X_\alpha$ 
such that $f(x_\alpha)\in[0,1)$ and 
$f(y_\alpha)\in[\alpha,\omega_1)$. Let $k,\ell$ 
be such that $k(\alpha)=k$ for uncountably many $\alpha$, and $\ell(\beta)=\ell$ for uncountably many $\beta$.
Then, $f$ is unbounded on $D_k$ and bounded on an unbounded subset of $D_\ell$, and thus, since $D_\ell$ is a direction, $f$ is 
bounded on $D_\ell$. Moreover, $k$ and $\ell$ can be arbitrarily large in $\omega$ since 
$D_{\mathcal{A}}=D_{\mathcal{A'}}$ with $\mathcal{A'}=\{D_j\,:m\le j<\omega\}$. Thus, if $D_j\preceq D_{j+1}$ for
all $j$ or if $D_j\succeq D_{j+1}$ for all $j$, there will be $k_1,k_2\in\omega$, with $D_{k_1}\preceq D_{k_2}$
and $f$ unbounded on $D_{k_1}$ and bounded on $D_{k_2}$, a contradiction. 
Taking $\mathcal{A}=\{C_j\,:\,j\le 0\}$, $\mathcal{A'}=\{C_j\,:\,j\ge 0\}$, 
$D=D_{\mathcal{A}},D'=D_{\mathcal{A'}}$ are thus directions. 
\\
We now show that $D\preceq C_i$ for all $i\in\mathbb{Z}$. Let $f:X\to\LL_{\ge 0}$ be bounded on $C_i$, then it is bounded on
all $C_j$ with $j\le i$. Let $\alpha$ be a bound for $f$ on $\cup_{j\le i}C_i$, then by continuity $\alpha$ is a bound
for $f$ on $D$, proving that $D\preceq C_i$. If now $f$ is unbounded on some $C_i$, it is unbounded on all $C_j$ with $j\ge i$.
Given $\alpha$, let $x_j\in C_j$ ($j\ge i$) be such that $f(x_j)\in[\alpha,\omega_1)$. Taking a convergent subsequence
we obtain a $x\in D'$ with $f(x)\in[\alpha,\omega_1)$, proving that $f$ is unbounded on $D'$, hence that 
$C_i\preceq D'$.
\endproof

It happens that it is not difficult to find an $\omega$-bounded manifold without minimal direction:

\begin{thm}\label{propufo} 
  There is an $\omega$-bounded surface $S$ spanned by directions for which  
  $\langle\mathfrak{D}_S,\preceq\rangle$ does not possess a minimal direction.
\end{thm}

\begin{figure}[h]
  \begin{center}
  \epsfig{figure=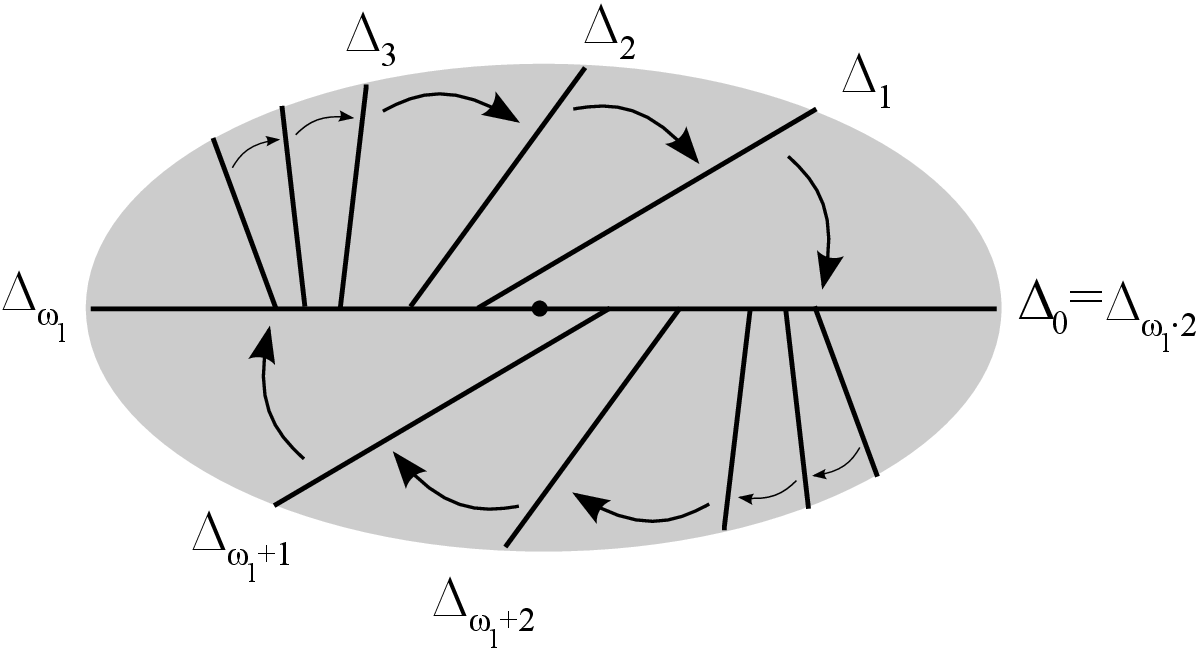, width=8cm}
  \caption{Theorem \ref{propufo}.}\label{figure5sec5}
  \end{center}
 \end{figure}

\proof[Proof of Theorem \ref{propufo}]
 Consider $S^{\omega_1\cdot 2,h,r}$, with $r(\alpha)=\downarrow$
 for all $\alpha\in\omega_1\cdot 2$, and $h(\beta)=h(\omega_1 + \beta)=\beta$ for each $\beta<\omega_1$. 
 Define $S$ by identifying $\Delta_0$ and $\Delta_{\omega_1\cdot 2}$
 (see Figure \ref{figure5sec5}).
 
 Then, if $\alpha<\beta<\omega_1$,
 $\Delta_{\alpha}\succ\Delta_\beta$ and  
 $\Delta_{\omega_1+\alpha}\succ\Delta_{\omega_1+\beta}$.
 But 
 $\Delta_{\omega_1}\not\preceq\Delta_{\alpha}$ for any $\alpha<\omega_1$, and neither
 $\Delta_0=\Delta_{\omega_1\cdot 2}$
 is $\preceq\Delta_{\omega_1+\alpha}$.
 Indeed, for $\Delta_{\omega_1}$ it is easy to find a continuous $f:S\to\LL_{\ge 0}$ such that $f\restrict{\Delta_{\omega_1}}$ is the identity,
 while $f\restrict{\Delta_{\alpha}}$ is constant on $\alpha$ after height $\alpha$ for $\alpha<\omega_1$ and on $0$ for $\omega_1<\alpha\le\omega_1\cdot 2$.
 Since this would be the only possible 
 lower bound, the chain $\langle\Delta_\alpha\,:\,\alpha<\omega_1\rangle$ does not possess such an lower bound. 
 We may find a similar map for $\Delta_{\omega_1\cdot 2}=\Delta_0$.
 Moreover any direction in $S^{\omega_1,h,r}$ 
 is ufo-equivalent to some $\Delta_\alpha$ by Corollary \ref{corchiant}, and
 therefore $\langle\mathfrak{D}_S,\preceq\rangle$ is isomorphic to two copies of $\omega_1^*$ such that the members lying in different copies
 are not comparable.
\endproof

The case of the maximum is trickier.
The naive idea of taking $S^{\omega_1,h,r}$ with $r(\alpha)=\uparrow$ for all $\alpha$ does not work:
$\Delta_{\omega_1}$ is ufo-maximal (and $\succeq\Delta_\alpha$ for each $\alpha$). 
We suspect that {\bf PFA} implies that an $\omega$-bounded manifold (or Type I locally compact first countable space) always contains an ufo-maximal
direction, but we have no evidence for this claim. 
(Notice that in that case, 
a positive result for $\omega$-bounded spaces {\em do not} imply the same result
for $\omega_1$-compact spaces.)
On the other hand, we have the following:

\begin{thm} 
  \label{thm:nomax}
  The following holds.\\
  (a) There is an $\omega$-bounded space $X$ containing directions without ufo-maximal direction.\\
  (b) {\rm($\neg${\bf SH})} There is an $\omega_1$-compact $3$-manifold $M$ spanned with directions without ufo-maximal direction,
      actually $\langle\mathfrak{D}_M,\preceq\rangle$ is a Suslin tree.\\
  $\ast$(c){\rm($\diamondsuit^+$)}
  There is an $\omega$-bounded $3$-manifold as in (b), with $\langle\mathfrak{D}_M,\preceq\rangle\simeq\langle\omega_1,\le\rangle$.
\end{thm}

\begin{figure}[h]
  \begin{center}
  \epsfig{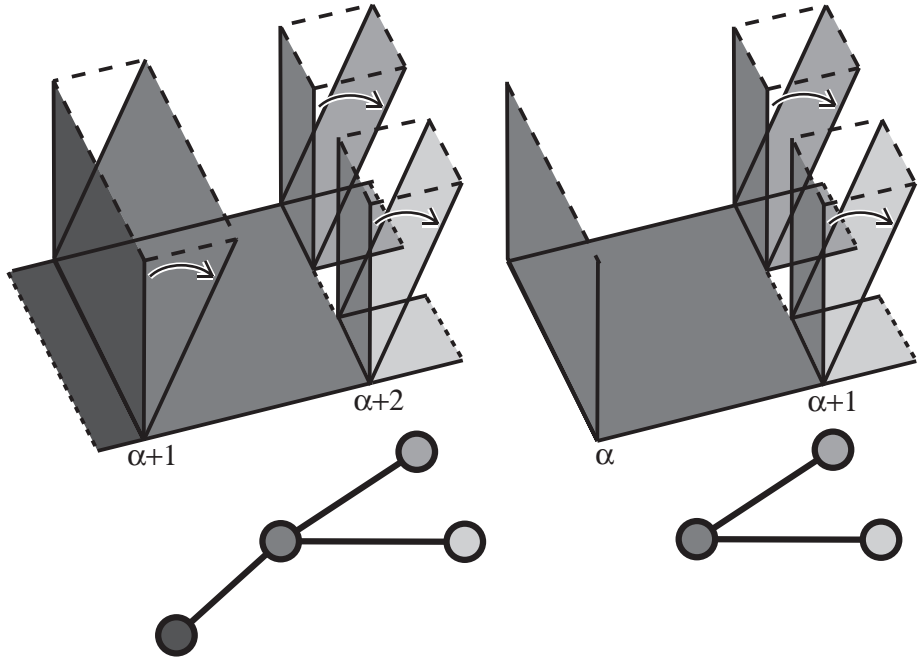}
  \caption{Reconstructing $T$ in $\langle\mathfrak{D}_{\Xi(S_T,A,s)},\preceq\rangle$.}\label{figure:XiSuslin}
  \end{center}
 \end{figure}

\proof
The proof uses the spaces $\Xi(X,A,s)$ of Section \ref{sec:Xi}. In each case $A$ will be the subset of successor ordinals in $\LL_{\ge 0}$.
For
(a), consider the space $\Xi(Q,A,s)$ with $Q$ given by Theorem \ref{thm:betaomega} 
and apply Lemma \ref{lemma:Xi2}.
The proof for (b) and (c) is very similar: we consider $\Xi(S_T,A,s)$ and $\Xi(N,A,s)$
where $T$ is a Suslin tree and $N$ is given by Theorem \ref{AxN}, with a discrete unbounded $A$.
In both cases the space is spanned by directions by Lemma \ref{lemma:Xi4}, and
Lemma \ref{lemma:Xi2} shows that there is no ufo-max direction.
With a bit more work, we will show that the ufo-equivalent classes of directions is 
precisely the tree $T$ in the first case and $\langle\omega_1,\le\rangle$ in the latter case.\\
Let $M=\Xi(S_T,A,s)$.
Call {\em verticals rays} the subsets of $\Xi(S_T,A,s)$ which are either $\{y\}\times\LL_{\ge s(y)}$ whenever $s(y)\not\in A$,
or an horizontal or diagonal ray in some inserted $\OO_{\ge \alpha+1}$ (which end up being `vertical' because
of the first factor $s^{-1}(\{\alpha+1\})$).
Notice that any direction must have a club intersection with a vertical ray: its
intersection with $\Delta(S_T)$ cannot be club, then by Lemma \ref{lemma:Xi1} it is contained in
$\Xi((S_T)_\alpha,A,s\restrict{(S_T)_\alpha})$.
If its intersection with one of the inserted $E_a$ is club, we are done, otherwise
we apply Lemma \ref{lemma:lemma2.2BGGG} to conclude.
\\
We will show later that we can choose $s$ such that
each component $\sigma$ of $s^{-1}(\{\alpha+1\})$ is homeomorphic to $\R$,
hence $\sigma\times\LL_{\ge \alpha+1}$ is a direction in itself.
Figure \ref{figure:XiSuslin} (left) below 
explains what happens `just above' this component. The only vertical parts depicted
are the vertical rays in the inserted pieces. All the rays which lie in zones of the
same shade of grey are ufo-equivalent, and $\preceq$ can be read by following the
lines from left to right on the 
small piece of tree below.
In particular, the three central vertical sides are ufo-equivalent (and so are the vertical rays between them).
The righthandside picture shows the limit level $\alpha$,
since no $E_\alpha$ is inserted we depicted the vertical rays at height $\alpha$ and
the two immediate successors (where there are inserted pieces).
We have thus constructed `jumps' between the directions in a perfect emulation of the tree $T$,
and thus the ufo-equivalent classes of direction are precisely $T$.
The case of $\Xi(N,A,s)$ is similar but easier: we show below that $s^{-1}(\{\alpha+1\})$ has only one component (homeomorphic to the circle),
so the inserted `jumps' between directions do form a chain of type $\omega_1$.

We now
show that $\Xi(S_T,A,s)$ and $\Xi(N,A,s)$ are indeed $3$-manifolds when $s$ is well chosen.
To this end we will for once
emulate the theory of foliations, and we thus start with some terminology (which is
not standard and valid just for this section).

\begin{defi}\label{def:A-foliation}
  Let $s:S\to\LL_{\ge 0}$ be a slicer on a a Type I surface $S$ and $A\subset\LL_{\ge 0}$.
  Then $s$ is an $A$-foliation of $S$ if for all $x\in s^{-1}(A)$, there are
  $a,b\in\LL_{\ge 0}$ with $a< s(x) < b$, an open subset $U\ni x$ and a map $\varphi: s^{-1}((a,b))\cap U\to \R$
  such that 
  $$\langle s,\varphi\rangle : s^{-1}((a,b))\cap U\to(a,b)\times\R$$ 
  is an homeomorphism.
  If $A=\LL_{\ge 0}$ we say that $s$ is a foliation of $S$.
\end{defi}

If $s$ is an $A$-foliation and $a\in A$,
then $s^{-1}(\{a\})$ is an at most countable disjoint union of $1$-dimensional metrizable submanifolds (called {\em leaves}) which are thus copies of
$\mathbb{S}^1$ or $\R$. It is convenient to have tubular neighborhoods (i.e. the space `between two leaves') that behave
well around such components. 

\begin{defi}
  An $A$-foliation on a Type I surface given by a slicer $s:S\to\LL_{\ge 0}$ is tame
  if  for each $a\in A$, $s^{-1}(\{a\})$ is a discrete collection of components, and for each component 
  $\sigma\subset s^{-1}(\{a\})$, $U$ can be chosen such that
  $U\supset\sigma$ in Definition \ref{def:A-foliation}.
\end{defi}

\begin{lemma}
  The following holds.\\
  (a) If $T$ is a binary Aronszajn tree, then there is a tame foliation on $S_T$.\\
  (b) If $S$ is an $\omega$-bounded Type I surface, there is a tame $A$-foliation on $S$ for $A$ the subset of
      successors ordinals.
\end{lemma}
\proof
\
\\
(a) Figure \ref{figure:S_Tfoliation} below shows what the foliation looks like in each piece of $S_T$
    and should convince the reader that the foliation is indeed tame.
\\
(b) Nyikos's bagpipe theorem \cite[Thm 5.14]{Nyikos:1984} shows that there is some $\alpha$ such
    that $S-S_\alpha$ is a finite union of longpipes, that is, Type I surfaces $P$
    such that for each $\beta$, $P_{\beta+1}$ is homeomorphic to the cylinder $\mathbb{S}^1\times(0,1)$
    and the boundary of $P_{\gamma+1}$ in it is homeomorphic to the circle whenever $\gamma<\beta$.
    We may thus take a small strip around these circles at successor height, define the scale
    accordingly in those strips, and then extend to the rest of the space as in Lemma \ref{lemma0}.
\endproof

Notice that we can assume that $N$ of Theorem \ref{AxN} is a longpipe: take one of the longpipes given by the 
bagpipe theorem and sew in a disk to close it. This gives an $\omega$-bounded surface without direction.

\begin{figure}[h]
  \begin{center}
  \epsfig{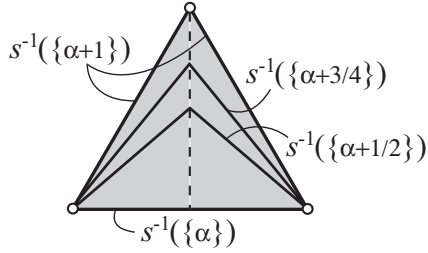}
  \caption{Foliating $S_T$.}\label{figure:S_Tfoliation}
  \end{center}
 \end{figure}

\begin{lemma}
  If $s$ is a tame $A$-foliation of a surface $S$ and $A$ is discrete, then 
  $\Xi(X,A,s)$ is a $3$-manifold (with boundary).
\end{lemma}
\proof
  Recall that $\we{S,s}$ is a $3$-manifold.
  Let $x\in\Xi(X,A,s)$. 
  We shall find a neighborhood of $x$ homeomorphic to $\R^3$ (or to $\R^2\times\R_{\ge 0}$).
  If $x$ lies in the interior of some inserted $E_a$, we are done.
  If $x$ is on the boundary of some inserted $E_a$, by discreteness of $A$
  we may easily describe such a neighborhood.
  Suppose that $x$ is not in some $E_a$.
  Consider $p(x)\in\we{S,s}$, then $p(x)\in\we{S,s}_\alpha$ for some $\alpha$.
  Notice that since $\we{S,s}$ is the subspace of $S\times\LL_{\ge 0}$ which is above
  the graph of $s$, there is an homeomorphism 
  from $\we{S,s}_\alpha$ to $S_\alpha\times[0,1)$ that leaves the first coordinate fixed,
  so we can assume that $p(x)\in S_\alpha\times[0,1)$.
  We choose some neighborhood $V$ of its projection on $S_\alpha$
  which is homeomorphic to $\R^2$ and 
  denote by $\phi:V\to\R^2$ an homeomorphism.
  Let us see how $V\times[0,1)$ changes when we insert the $E_a$ (for $a<\alpha$).
  In $S\times[0,1)$
  these insertions  boil down to taking out $s^{-1}(\{a\})\times[0,1)$
  and sewing in $s^{-1}(\{a\})\times P$, where $P$ is $\{\langle x,y\rangle\in[0,1)^2\,:\,y\le x\}$.
  Let $\mathcal{L}$ be the set of leaves in $\displaystyle\cup_{a\in A\cap[0,\alpha)}s^{-1}(\{a\})$.
  Notice that each component in $s^{-1}(\{a\})$ is closed in $S$ since the $A$-foliation is tame.
  Thus, by discreteness of $A$,
  for each leaf $\sigma\in\mathcal{L}$ 
  we may choose a neighborhood $U_\sigma$ 
  with an homeomorphism 
  which sends $\sigma$ to a vertical line. We can ensure that the $U_\sigma$ are pairwise disjoint.
  It is then easy to change the homeomorphism $\phi\times id:V\times[0,1)\to\R^2\times[0,1)$
  inside each $(V\cap U_\sigma)\times[0,1)$ to `absorb' 
  the inserted $E_a$, as shown in Figure \ref{figure:ultima} below. 
  The northwest picture shows a leaf $\sigma$ times $[0,1]$ together with
  a tubular neighborhood while the southwest picture shows the effect of $\phi\times id$.
  The northeast and southeast pictures show the situation after the insertion of $E_a$.
\endproof

This finishes the proof of Theorem \ref{thm:nomax}.
\endproof

 \begin{figure}[h]
  \begin{center}
  \epsfig{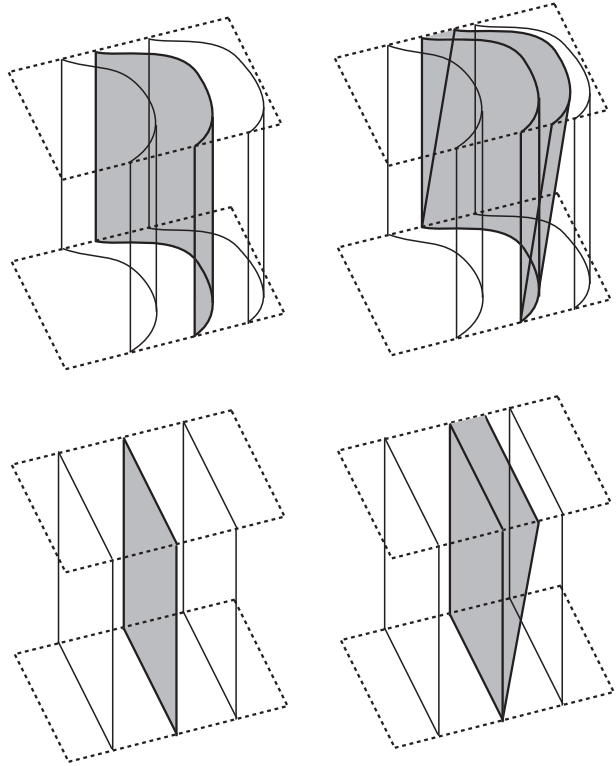}
  \caption{Looking at the insertions through the foliation.}\label{figure:ultima}
  \end{center}
 \end{figure}


\subsection{A partial order on maps $X\to\LL_{\ge 0}$ and its relations with the ufo}\label{sec:partialordermaps}

\begin{defi}
  Let $X$ be Type I and $f,g:X\to\LL_{\ge 0}$
  be both unbounded. We say that $f\trianglelefteq g$ iff
  for all clubs $C\in X$, `$f\restrict{C}$ is bounded'
  implies `$g\restrict{C}$ is bounded'. If $f\trianglelefteq g\trianglelefteq f$, we say that
  $f$ and $g$ are club-equivalent.
  We denote the strict order by $\vartriangleleft$ and by $\langle\mathfrak{F}_X,\trianglelefteq\rangle$ the poset of club equivalent classes.
\end{defi}
(As for $\preceq$, $\trianglelefteq$ is actually a preorder.)
Notice that $\trianglelefteq$ goes the other way than $\preceq$.
The idea is that a function which is unbounded only on directions that are `big' for $\preceq$ must be big for $\trianglelefteq$.
This is made more precise in Proposition \ref{propmapmax} below.

\begin{exs}\label{exs:F}
   \ \\
   $\bullet$ 
      $f:\LL_{\ge 0}^2\to\LL_{\ge 0}$, $(x,y)\mapsto\min\{x,y\}$ is $\trianglelefteq$-maximal
      in $\mathfrak{F}_{\LL_{\ge 0}^2}$.
   \\
   $\bullet$
      A slicer $X\to\LL_{\ge 0}$ is always $\trianglelefteq$-minimal.
   \\ 
   $\bullet$
      Let
      $X=S^{\omega,h,r}$, with $r(n)=\uparrow$ for all $n\in\omega$
      and $h:\omega\to\omega_1$ be the inclusion. Then $\mathfrak{F}_X$ does not
      possess a $\trianglelefteq$-maximal map.
      Indeed, an unbounded map $f$ must be unbounded on $\Delta_\omega$,
      but then for some $n$ it will be unbounded on all $\Delta_k$ for $k\ge n$ and bounded on the others.
      One then easily finds a map $g$ that has the same properties but for $n+1$ so $f\vartriangleleft g$.
   \\
   $\bullet$
      $X=S^{\omega_1,h,r}$, with $r(\alpha)=\uparrow$ for all
      $\alpha\in\omega_1$
      and $h:\omega_1\to\omega_1$ the identity. Then $\mathfrak{F}_{X}$ has a 
      $\trianglelefteq$-maximal 
      map that is the identity on $\Delta_{\omega_1}$ and bounded on all
      $\Delta_\alpha$ with $\alpha<\omega_1$ (recall Theorem \ref{propufo}).
\end{exs}

The connections between $\preceq$ for clubs subset and $\trianglelefteq$ for maps may be quite
evasive in general, but are
relatively explicit in the $\omega_1$-compact framework, and we shall restrict our attention to this case.
The first step is to find some kind of `minimal' club set on which a function is unbounded. 
This is provided by the $C(f,D)$ of Lemma \ref{lemmastar}. If $X$ is the entire space,
we denote $C(f,X)$ simply by $C(f)$. An immediate consequence of the definition of $C(f)$ is:

\begin{lemma}\label{encoreunlemme}
  If $X$ is a Type I $\omega_1$-compact 
  space, 
  $f:X\to\LL_{\ge 0}$ continuous, $D\subset X$ club, then $f\restrict{D}$ is unbounded iff $D\cap C(f)$ is club.
\end{lemma}

Whether $f\trianglelefteq g$ can be read directly in $C(f)$ and $C(g)$:

\begin{lemma} \label{lemmamaps2} 
    Let $X$ be an $\omega_1$-compact Type I space and $f,g:X\to\LL_{\ge 0}$ be
    continuous. Then the following holds.\\
    (a) $f\trianglelefteq g$ iff there is no club $C$ in 
    $C(g) -  C(f)$,\\
    (b) $f\vartriangleleft g$ iff $f\trianglelefteq g$ and 
    $C(f) -  C(g)$ contains a club.   
\end{lemma}

\proof
   \ \\
   (a) By Lemma \ref{encoreunlemme}, 
   $f\trianglelefteq g$ iff for all clubs $C\in X$, $C\cap C(g)$ is club
   implies that $C\cap C(f)$ is club, which is equivalent to saying that
   there is no club in $C(g) -  C(f)$.
   Indeed, if there is a club $C\subset (C(g) -  C(f))$, then $C\cap C(g)$ is
   club while $C\cap C(f)=\varnothing$. Conversely, suppose that there is no club in
   $C(g) -  C(f)$ and $C\cap C(g)$ is club. Then if $C\cap C(f)$ is bounded,
   say included in $X_\alpha$, by considering $C'=C -  X_\alpha$
   we have $C'\cap C(f)=\varnothing$, but then 
   $C(g) -  C(f)\supset (C'\cap C(g))$. Therefore 
   $C\cap C(f)$ must be unbounded.
   \\
   (b) $f\vartriangleleft g$ iff $f\trianglelefteq g$ and
   there is a club $C$ with $f\restrict{C}$ unbounded and $g\restrict{C}$ bounded,
   thus $C\cap (C(f) - C(g))$ is club and we conclude through
   Lemma \ref{encoreunlemme}.
\endproof

\begin{lemma}\label{lemmamin} 
   Let $X$ be an $\omega_1$-compact Type I space and $f,g:X\to\LL_{\ge 0}$ be continuous.
   Then $h:X\to\LL_{\ge 0}$ defined by $h(x)=\min\{f(x),g(x)\}$ (written $h=\min(f,g)$) is continuous and
   $C(h)= C(g)\cap C(f)$.
\end{lemma}

\proof
  Continuity is immediate, and $C(h)$ is the union of the points in the bones
  $B_\alpha(X)$ that are sent in $[\alpha,\omega_1)$ by $h$, so the second claim
  follows.
\endproof

\begin{prop}\label{propmapmax}
  Let $X$ be an $\omega_1$-compact Type I space, and $f:X\to\LL_{\ge 0}$ be continuous. Then the following holds.\\
  (a) If $C(f)$ is a direction, then $f$ is $\trianglelefteq$-maximal, and 
  $C(f)$ is a $\preceq$-maximal direction. \\
  (b) If $f$ is $\trianglelefteq$-maximal, then $C(f)$ is a $\preceq$-maximal direction in $X$.
\end{prop}

\proof 
  \ \\
  (a) 
  Suppose that there is some $g$ with $f\vartriangleleft g$. 
  Then $C(f) -  C(g)\supset C$ for some club $C$ and
  $g\restrict{C}$ is bounded. But since $C(g) -  C(f)$ contains no club, $C(g)\cap C(f)$ is club, so 
  $g\restrict{C(f)}$ is unbounded, which contradicts the fact that $C(f)$ is a   
  direction.
  Suppose now that there is some direction $D$ with $C(f)\prec D$. Let
  $g:X\to\LL_{\ge 0}$ satisfy $g\restrict{C(f)}$ bounded and 
  $g\restrict{D}$ unbounded, hence 
  $C(g)\cap D$ is club. 
  Up to taking $C(g)\cap D$, we can assume that $D\subset C(g)$.
  Since $C(g)\cap C(f)$ is bounded by Lemma \ref{encoreunlemme}, so is $D\cap C(f)$, which implies
  $f\restrict{D}$ bounded. But this contradicts $C(f)\preceq D$,
  since $f$ is unbounded on $C(f)$.
  \\
  (b) 
  We show that $C(f)$ is a direction, its maximality then follows from
  (a).
  Suppose it is not the case, there is therefore $g:X\to\LL_{\ge 0}$ with
  $g^{-1}(\{0\})\cap C(f)$ and $C(g)\cap C(f)$ both club in $X$.
  Set $h=\min(f,g)$, then $C(h) -  C(f)$ contains no club by Lemma   \ref{lemmamin},
  and $h^{-1}(\{0\})\supset g^{-1}(\{0\})$, therefore 
  $h^{-1}(\{0\})\cap C(f)$ is club in $X$, so 
  $C(f) -  C(h)$ contains a club. 
  By Lemma \ref{lemmamaps2} this means that $f\vartriangleleft h$, contradicting
  the maximality of $f$.
\endproof

Of course, there are spaces that do contain directions (and even maximal ones)
that do not contain maximal maps, as $S^{\omega,h,r}$ in Example \ref{exs:F}.
More generally, if there is a sequence of maps
$f_i$ ($i\in\omega$), with $f_j\vartriangleleft f_i$ when $j<i$, such that there is no
$g$ with $f_i\vartriangleleft g$ for all $i$, while there is no $\trianglelefteq$-maximal map,
it can be thought that $\cap_{i\in\omega}C(f_i)$ is a direction (if non-empty). But
there are easy examples where this is not true: take
$S^{\omega+1,h,r}$, where $h$ is the inclusion, $r(n)=\uparrow$
for all $n\in\omega$ and $r(\omega)=\downarrow$.
Take $f_i$ ($i\in\omega$) such that $f_i$ is unbounded on $\Delta_j$ for
all $j\ge i$ (and on $\Delta_{\omega+1}$),
and bounded on the others. It is easy to see that there is no
$g$ bigger than each $f_i$, but that $\cap_{i\in\omega}C(f_i)$ intersects
both $\Delta_\omega$ and $\Delta_{\omega+1}$ in a club way, and is therefore not
a direction.
The following lemma shows that we can obtain a direction if we impose more 
restrictions on the sequence $f_i$.

\begin{lemma}
  Let $X$ be an $\omega$-bounded Type I space. Let $f_i:X\to\LL_{\ge 0}$ ($i\in\omega$)
  be a sequence with $f_i\preceq f_j$ when $i\le j$.
  Then the following holds.
  \\
  (a) If $\cap_{i\in\omega}C(f_i)$ is a direction, then there is no $g:X\to\LL_{\ge 0}$
  with $f_i\vartriangleleft g$ for all $i\in\omega$, and $\cap_{i\in\omega}C(f_i)$ 
  is $\preceq$-maximal. \\
  (b) Suppose that
  there is no $\trianglelefteq$-increasing sequence $g_i$ ($i\in\omega$)
  satisfying: $f_i \vartriangleleft g_i$ for all $i$, and for $j\ge i$ 
  $g_i\not\trianglelefteq f_j$.
  Then $\cap_{i\in\omega}C(f_i)$ is a $\preceq$-maximal direction. 
\end{lemma}

It is not very clear however whether there are examples where condition (b) holds.

\proof 
  If the sequence of $f_i$ stagnates after some $i$ then Proposition \ref{propmapmax}
  applies. Assume that $f_i\vartriangleleft f_j$ when $i<j$. By Lemma \ref{cpomegaclosed}, $\cap_{i\in\omega}C(f_i)$ is club.
  \\
  (a) The proof is the same than that of Proposition \ref{propmapmax} once 
  one notes that if $g\restrict{\cap_{i\in\omega}C(f_i)}$ is bounded, then 
  $g\restrict{C(f_i)}$ must be bounded for all
  big enough $i$.\\
  (b)
  If 
  $C=\cap_{i\in\omega}C(f_i)$ is not a direction, there is a $g$ unbounded on
  $C$ with
  $g^{-1}(\{0\})\cap C$. Define $g_i=\min(g,f_i)$ and derive a contradiction
  as in Proposition \ref{propmapmax}. 
\endproof

It is easy to see that a map that is club-equivalent to a slicer is indeed a slicer for another
canonical sequence.
Moreover, a slicer has the property that $f^{-1}(\{x\})$ is bounded for each $x\in\LL_{\ge 0}$, but 
in general this property does not define slicers:
Take
$X$ to be the disjoint union of $\omega_1$ with the usual topology and $\omega_1$ with the discrete topology.
Then the map $f$ given by the identity on the first copy and by some $1$-to-$1$ map $\omega_1\to[0,1]$ on the second
copy has the alleged property but is not club-equivalent to a slicer.
We shall see below that both notions are equivalent if the space is $\omega$-bounded.

\begin{defi}
  Let $X$ be Type I, and $f:X\to\LL_{\ge 0}$. 
  We set
   \begin{align*}
     \mathsf{Bd}(f)&=\{x\in\LL_{\ge 0}\, :\, f^{-1}(\{x\})\text{ is bounded}\}\\
     \mathsf{Ubd}(f)&=\{x\in\LL_{\ge 0}\, :\, f^{-1}(\{x\})\text{ is unbounded}\}.
   \end{align*}
\end{defi}

First we note the following:

\begin{lemma}\label{lemmaBdUbd}
   If $X$ is Type I and $\omega$-bounded, and $f:X\to\LL_{\ge 0}$,
   then $\mathsf{Ubd}(f)$ is closed in $\LL_{\ge 0}$.
\end{lemma}

\proof
  Since $\LL_{\ge 0}$ is first countable we need only to check that $\mathsf{Ubd}(f)$ is sequentially closed.
  Fix an arbitrary $\alpha\in\omega_1$, and let
  $\{y_i\,:\,i\in\omega\}$ be a sequence in $\mathsf{Ubd}(f)\subset \LL_{\ge 0}$
  converging to $y$. Let $x_i\in f^{-1}(\{y_i\}) -  X_\alpha$. 
  Then for any accumulation point $x$ of 
  the $x_i$, $f(x)=y$, and $x\in X -  X_{\alpha}$.
\endproof

\begin{lemma}
   Let $X$ be Type I, and $f:X\to\LL_{\ge 0}$. Then,
   (a) and (b) below are equivalent and imply (c). If $X$ is $\omega$-bounded, then (c) implies (a) and (b) as well.\\
   (a) $f$ is club-equivalent
   to a slicer, \\
   (b) There is a club $C\subset\omega_1$ such that
   $f$ is a slicer for the canonical sequence $\langle X_\alpha:\alpha\in C\rangle$,\\
   (c) $\mathsf{Ubd}(f)=\varnothing$.
\end{lemma}

\proof
   Note that a map is club equivalent to a slicer if and only if it is unbounded on
   each club of $X$.
   \\
   (a) $\rightarrow$ (b).
   If $f$ is club-equivalent to a slicer, then
   $f^{-1}([0,\alpha])$ is bounded for all $\alpha\in\omega_1$.
   By a leapfrog argument this implies that the set of $\alpha$ such that
   $f^{-1}([0,\alpha))\subset[0,\alpha)$ is club. Since the set of $\alpha$ 
   with $f(\wb{X_\alpha})\subset[0,\alpha]$ is club too, 
   letting $C$ be their intersection yields that $f$ is a slicer for the canonical sequence
   given by $C$.
   \\
   (b) $\rightarrow$ (a). Immediate.\\
   (b) $\rightarrow$ (c).
   If $C=f^{-1}(\{x\})$ is unbounded, then $f$ cannot be a slicer for any canonical sequence, otherwise $f(C)$ would be unbounded.
   \\
   (c) $\rightarrow$ (a) if $X$ is $\omega$-bounded.
   Suppose that there is some club $C\subset X$ with $f(C)\subset [0,\alpha]$ for some
   $\alpha$. We define inductively closed intervals 
   $I_i\subset I_{i-1}\subset\cdots\subset I_0= [0,\alpha]$ by bisecting 
   $I_i$ and choosing $I_{i+1}$ as one of the two pieces such that 
   $f^{-1}(I_i)\cap C$ is unbounded. We may choose $I_i$ such that
   $\cap_{i\in\omega}I_i=\{x\}$, and since $X$ is countably compact,
   $\langle\mathfrak{C}_X,\subseteq\rangle$ is countably closed (cf Lemma \ref{cpomegaclosed}), so 
   $\cap_{i\in\omega}(f^{-1}(I_i)\cap C)=f^{-1}(\{x\})\cap C$ is club,
   which shows that $x\in\mathsf{Ubd}(X)\not=\varnothing$, a contradiction.
\endproof




\section{Dominant subsets}\label{sec:dominant}

There is a natural strengthening of the notion of direction which is the following.

\begin{defi} 
Let $X$ be Type I.
A club $D\subset X$ is {\em locally dominant} if there is an open set $U\supset D$ such that for
all $f:X\to\LL_{\ge 0}$, all club $C\subset U$ and all club $D'\subset D$, $f\restrict{D'}$ bounded implies $f\restrict{C}$ bounded. 
If $U=X$, we say that $D$ is globally dominant.
\end{defi}

The requirement about the club subsets $D'\subset D$ is here to prevent $X$ itself to 
be always globally dominant. A locally dominant subset is trivially a direction (take $C=D$), but not
conversely. A locally dominant subset of a Type I manifold is not always ufo-maximal:
the horizontals are locally dominant in $\OO$ (take an open horizontal strip). 
The reader may check that $\Delta_{\omega_1}$ and $\Delta_0$ are both ufo-maximal in the
surface
$S^{\omega_1\cdot 2,h,r}$ described in Theorem \ref{propufo}
but are not locally dominant (see Theorem \ref{strongprob3false} below for a proof in a similar case),
so local dominance and ufo-maximality, though similar, are logically unrelated.

\begin{defi}
A Type I space is spanned by dominant subsets if any large open subset contains a locally dominant subset.
\end{defi}

We note that the surface of Example \ref{exDominant} gives us the following.

\begin{thm}[P. Nyikos \cite{NyikosPC}]\label{thmnld}
  There is an $\omega$-bounded surface spanned by directions which does not contain any locally dominant subset.
\end{thm}
 
\proof
  By Lemma \ref{nld2} any large open set of $R$ contains the terminal part of a ray, and the rays
  are not locally dominant by construction.
\endproof

Notice as well that given $\alpha<\omega_2$ and $h,r$ as in Lemma \ref{corchiant0}, the surfaces
$S^{\alpha,h,r}$ are spanned by dominant subsets.
Indeed, given a club $C$ inside an open subset $U$, either $C$ contains a locally dominant subset or
$C$ club-intersects some $\Delta_\alpha$ for a limit $\alpha$. 
In the latter case, 
$U$ contains the terminal part of $\Delta_\beta$ for some (smaller) successor $\beta$ which is locally dominant
(see the proof of Lemma \ref{lemmachiant2}).
Thus, the surfaces given by Theorems \ref{thmufonominnomax} and \ref{propufo} are in fact spanned by dominant subsets.
The counterpart of Theorem \ref{3d} for dominant subsets is easy to prove in {\bf ZFC}:

\begin{thm}\label{strongprob3false}
   There is an $\omega$-bounded surface $S$ spanned by dominant subsets such that there is an unbounded 
   $f:S\to\LL_{\ge 0}$ that is bounded on each locally dominant subsets of $S$.
\end{thm}

\proof[Proof of Theorem \ref{strongprob3false}]
  Let $S$ be $S^{\omega_1,h,r}$ with $h:\omega_1\to\omega_1$ 
  the identity
  and
  $r:\omega_1\to\{\uparrow,\downarrow\}$ being such that for uncountably many $\alpha$,
  $r(\alpha)=\uparrow$ and $r(\alpha+1)=\downarrow$.
  We show that no club subset $C$  
  of $\Delta_{\omega_1}$ can be locally dominant.
  Given an open $U$ containing $C$, there is a compact $K$ and a $\beta<\omega_1$ such that  
  $U$ contains $\Delta_\alpha - K$ whenever
  $\beta<\alpha\le\omega_1$.
  If $r(\alpha)=\uparrow$ and $r(\alpha+1)=\downarrow$ there is
  a map $f:S^{\omega_1,h,r}\to\LL_{\ge 0}$ which is unbounded on $\Delta_\gamma$ iff $\gamma={\alpha+1}$, 
  and thus $C$ is not
  locally dominant.
  It is now easy to find a map $f:S^{\omega_1,h,r}\to\LL_{\ge 0}$ such that $f\restrict{\Delta_\alpha}$ is 
  bounded by $\alpha$ for each $\alpha<\omega_1$, and the identity on $\Delta_{\omega_1}$,
  so $f$ is bounded on any locally dominant subset of $S$. 
\endproof

Notice in passing that $\Delta_{\omega_1}$, being not comparable with any other $\Delta_\alpha$, is ufo-maximal.



\section{The role of $\LL_{\ge 0}$ and a Hausdorff peculiarity}\label{sec:Hausdorff}

In this section we investigate the question of the `universality' of $\LL_{\ge 0}$ in the definition of directions.
To be precise, let us define:

\begin{defi}
  Let $X,Y$ be Type I spaces, and let $D\subset X$ be club. Then $D$ is a $Y$-direction in $X$ iff
  given $f:X\to Y$ with $f$ unbounded on $D$, 
  then $f$ is unbounded any club subset of $D$ as well.
  If there exists an $f:X\to Y$ which is unbounded on $D$, the direction is said to be true. 
\end{defi}

Directions will thus be called $\LL_{\ge 0}$-directions in this section.
Almost everything done in Section \ref{sec:gendir} holds for $Y$-directions as well.
In particular, the equivalences in Lemma \ref{pulpita} hold, except clause (c) which may fail
if there is no map $f:Y\to Y$ with $f(Y_\alpha)=\{0\}$ and $f\restrict{Y-Y_{\alpha+1}}=id$.

\begin{lemma} \label{lemma:abcd}
  Let $X,Y$ be Type I spaces, then
  the following holds.\\
  (a) $\LL_{\ge 0}$ and $\omega_1$ are $Y$-directions in themselves. \\
  (b) If $X$ is regular and $D$ is an $\LL_{\ge 0}$-direction in $X$, then $D$ is a true $\LL_{\ge 0}$-direction in $X$.\\
  (c) If there is a slicer on $Y$, then any $\LL_{\ge 0}$-direction in $X$ is a $Y$-direction in $X$.\\
  (d) If there is an unbounded map $g:\LL_{\ge 0}\to Y$, then a $Y$-direction in $X$ is an $\LL_{\ge 0}$-direction in $X$. 
\end{lemma}
\proof
  \ \\
  (a) Given $f:\LL_{\ge 0}\to Y$, if $f^{-1}(Y_\alpha)$ is unbounded then $f^{-1}(\wb{Y_{\alpha}})$ is club.
      It follows that $f^{-1}(X-Y_{\alpha+1})$ cannot be club, hence $f$ is bounded. \\
  (b) By Lemma \ref{lemma0} there is a slicer $s:X\to \LL_{\ge 0}$ which is thus unbounded on $D$. \\
  (c) Given $f:X\to Y$ unbounded on $D$ and $\alpha\in\omega_1$, then $s\circ f(D-X_{\gamma(\alpha)})\subset [\alpha,\omega_1)$ for some $\gamma(\alpha)$,
      and thus $f(D-X_{\gamma(\alpha)})\subset Y-Y_\alpha$.\\
  (d) Let $f:X\to \LL_{\ge 0}$ be unbounded on $D$. If $g\circ f$ is bounded on $D$,
      then $g$ sends an unbounded subset of $\LL_{\ge 0}$ to a bounded subset of $Y$.
      By (a) $g$ must be bounded, a contradiction. Hence $g\circ f$ is unbounded.
      Take a club $C$ such that $g\circ f(X-X_\alpha)\subset Y-Y_\alpha$ and $g([0,\alpha))\subset Y_\alpha$ for each $\alpha\in C$.
      Then $f(X-X_\alpha)\subset [\alpha,\omega_1)$.
\endproof

Clause (c) shows that $\LL_{\ge 0}$ is a kind universal target for checking the `directionity' in regular spaces.
Though, the $Y$-direction might not be true, even if it is a true $\LL_{\ge 0}$-direction.
For instance if $Y$ is
the positive part of a tangent bundle of $\LL_+$ as in Example \ref{ex:NyikosClass7},
then no map $\LL_{\ge 0}\to Y$ is unbounded. By
quotienting by a suitable $\Z$-action
on the fibers, as explained in \cite{Nyikos:1992} and \cite[Section 6]{mesziguessurf},
this example can be made $\omega$-bounded.
If $Y$ is not (path) connected, then a $Y$-direction may well not be an $\LL_{\ge 0}$-direction,
for instance $\omega_1\times\LL_{\ge 0}$ is an $\omega_1$-direction in itself.

The aim of this section is to prove the following theorems, which show that the conclusions of (b), (c) above may fail in 
the class of Hausdorff connected spaces, while (d) may fail as well for normal connected spaces.

\begin{thm}\label{thm:truedir}
  There are Type I Hausdorff non-regular spaces $R,X$, with $R$ connected, such that the following holds:
  \\
  (a) $R$ is an $\LL_{\ge 0}$-direction in itself which is not true:
      any function into $\LL_{\ge 0}$ is constant.\\
  (b) $X$ contains two clubs $D_1,D_2$ such that $D_1$ is a true $\LL_{\ge 0}$-direction but not an $R$-direction,
      and $D_2$ is a true $R$-direction but not an $\LL_{\ge 0}$-direction.
\end{thm}

\begin{thm}\label{thm:notLdir}
  There are $\omega$-bounded Type I spaces $X,T$ which are normal, 
  connected and locally metrizable, and a club $D\subset X$ which is a $T$-direction but not 
  an $\LL_{\ge 0}$-direction in $X$. 
\end{thm}

$T$ is not path connected, which is somehow a cheat.
The idea is to find a space $Y$ such that a map $\LL_{\ge 0}\to Y$ cannot 
`move much', and to apply the following technical lemma.

\begin{lemma}\label{lemma:noZdir}
  Let $Z,Y$ be Type I. Suppose that for all $f:Z^2\to Y$ and $q\in Z^2$, 
  if $f(q)\in Y-Y_\alpha$, then $f(Z^2)\subset Y-Y_\alpha$.
  Let $D\subset Sk(Y)$ be a $Y$-direction in $Y$.
  Then $E=\cup_{\alpha\in\omega_1}(D\cap B_\alpha(Y))\times\wb{Z_\alpha}^2$ is a true $Y$-direction
  in $Y\times Z^2$, and not a $Z$-direction.
\end{lemma}
\proof 
  Write $X=Y\times Z^2$ and
  denote by $\pi_i$ the projection on the $i$-th factor. Then
  $\pi_3:Y\times Z^2\to Z$ shows that $E$ is not a $Z$-direction in $X$. 
  Let $f:X\to Y$ be unbounded on $E$. 
  If $f(E\cap Y\times\{q\})$ is bounded for each $q\in Z^2$, then $f(E)$ is bounded.
  (Otherwise for all $\alpha$ there are $y_\alpha,q_\alpha$ with
  $f(y_\alpha,q_\alpha)\in E\cap (Y-Y_\alpha)$, then
  $f(\{y_\alpha\}\times Z^2)\subset E\cap(Y-Y_\alpha)$ and $f(E\cap Y\times\{q\})$ is unbounded for each $q\in Z^2$.)
  Hence there is some $q$ for which $f(E\cap Y\times\{q\})$ is unbounded.
  Since $D$ is a $Y$-direction in $Y$, there is a club $C\subset\omega_1$ such that
  $f((D-Y_\alpha)\times\{q\})\subset Y-Y_\alpha$ for $\alpha\in C$. 
  Thus, $f((D-Y_\alpha)\times Z^2)\subset Y-Y_\alpha$.
  Since $E-X_\alpha = E\cap (D-Y_\alpha)\times Z^2$, 
  $E$ is a $Y$-direction, which is true since the projection $\pi_1:X\to Y$ is unbounded on $E$.
\endproof

The attentive reader may have noted 
the similarity between $E$ and the spaces $\ve{X}$ of Section \ref{sec:constructions}.
We shall now define
$R=\cup_{\alpha\in\omega_1}R_\alpha$ such that
$\wb{R_\alpha}$ is countable and connected for each
limit $\alpha$. 
Our construction 
is a simple variant on Roy's lattice space (Example 126 in \cite{CEIT}).
Denote by $\Lambda$ the limit ordinals in $\omega_1$.
Let $\{C_n\,:\,n\in\omega\}$ be a disjoint collection of dense subsets of $\Q$. We assume that $0\in C_0$.
For each $\alpha\in\omega_1$, write $\alpha = \beta + \varphi(\alpha)$ with $\varphi(\alpha)\in\omega$ and $\beta\in\Lambda$.
We let then $R$ be $\{\langle r,\alpha\rangle\,:\,\alpha\in\omega_1,r\in C_{\varphi(\alpha)}\}$. 
When $\alpha\in\Lambda$, we denote by $q_\alpha$ the point $\langle 0,\alpha\rangle$.
The topology is the following. 
Given $\langle r,\alpha\rangle\in R$ and $\epsilon >0$, denote by 
$U_{r,\alpha,\epsilon}$
the set $\{\langle s,\alpha\rangle\in R\,:\,|r-s|<\epsilon\}$.
If $\varphi(\alpha)$ is odd, 
then neighborhood of $\langle r,\alpha\rangle$ is given by the $U_{r,\alpha,\epsilon}$ for $\epsilon >0$.
A neighborhood of $\langle r,\alpha\rangle$ when $\varphi(\alpha)$ is even and $>0$
is a stack of $3$ intervals $U_{r,\alpha -1,\epsilon}\cup U_{r,\alpha,\epsilon}\cup U_{r,\alpha +1,\epsilon}$ for $\epsilon>0$
(notice that $\alpha -1$ makes sense when $\varphi(\alpha)>0$).
If $\varphi(\alpha)=0$ (and thus $\alpha\in\Lambda$), a neighborhood of 
a point different from $q_\alpha$ is given by a stack of $2$ intervals $U_{r,\alpha,\epsilon}\cup U_{r,\alpha +1,\epsilon}$ for $\epsilon>0$.
Finally, the neighborhoods of $q_\alpha=\langle 0,\alpha\rangle$
are given by $U_{0,\alpha,\epsilon}\cup U_{0,\alpha +1,\epsilon}\cup V_{\gamma}$,
where $\epsilon>0$ ,$\gamma<\alpha$, and
$V_{\gamma}=\{\langle r,\beta\rangle\,:\,\gamma < \beta < \alpha\}$.
We let $R_\alpha=\{ \langle r,\beta\rangle\,:\,\beta<\alpha\}$.
  
\begin{lemma}\label{lemma:RTypeI}
  $R$ is Hausdorff, Type I and each $\wb{R_\alpha}$ is connected (thus $R$ as well).
\end{lemma}
\proof
  Hausdorffness is immediate since the $C_n$ are disjoint.
  $R=\cup_{\alpha\in\omega_1}R_\alpha$,
  and $\wb{R_\alpha}$ is $R_{\alpha}$ if $\varphi(\alpha)$ is even and $>0$,
  $R_{\alpha+1}$ if $\varphi(\alpha)$ is odd, and
  $R_{\alpha}\cup\{q_\alpha\}$ if $\alpha\in\Lambda$.
  In any case, it is countable and thus Lindel\"of, so $R$ is of Type I.
  A clopen set containing 
  an entire level
  $C_{\varphi(\alpha)}\times\{\alpha\}$
  must contain $C_{\varphi(\beta)}\times\{\beta\}$
  for $\beta$ the immediate predecessor of $\alpha$ if $\alpha$ is successor, or $V_\gamma$ for some $\gamma<\alpha$
  if $\alpha\in\Lambda$.
  (See \cite[Example 126]{CEIT} for details,
  notice that we switched the role of `odd' and `even' in our version.)
  Thus, any clopen set containing $q_\alpha$ for $\alpha\in\Lambda$ contains 
  $\wb{R_\alpha}$, which is then connected.
\endproof

\begin{lemma}\label{lemma:Rnottrue}
  For any $n,m\in\omega$,
  Any $f:\LL_{\ge 0}^n\to R^m$ is constant, and any $g:R^n \to \LL_{\ge 0}^m$ is constant.
\end{lemma}
\proof
  A canonical sequence for $R^n$ is given by $R_\alpha^n$, and $\wb{R_\alpha^n}=\wb{R_\alpha}^n$.
  For a club set of $\alpha$, $f([0,\alpha]^n)\subset \wb{R_\alpha}^m$.
  Thus
  $[0,\alpha]^n=\cup_{x\in\wb{R_\alpha}^m}f^{-1}(\{x\})$ and is therefore
  a countable union of disjoint closed subsets. This is impossible
  except if all but one are empty. The function is thus constant on $[0,\alpha]^n$
  for a club set of $\alpha$, and thus on all of $\LL_{\ge 0}^n$. \\
  Notice that $\wb{R_\alpha}^n$ is connected,
  thus $g$
  must send $\wb{R_\alpha}^n$ to a connected and countable subspace of $\LL_{\ge 0}^m$, 
  and thus be constant on $\wb{R_\alpha}^n$.
  The map is is thus constant on all of $R^n$.  
\endproof

\begin{lemma} \label{lemma:qtrueRdir}
  The set
  $\{q_\alpha\,:\,\alpha\in\Lambda\}$ is a true $R$-direction.
\end{lemma}
\proof
  $\{q_\alpha\,:\,\alpha\in\Lambda\}$ is a closed copy of $\omega_1$, and we apply Lemma \ref{lemma:abcd} (a).
\endproof

\proof[Proof of Theorem \ref{thm:truedir}]
\
\\
(a) By Lemma \ref{lemma:Rnottrue}.\\
(b) Let $X$ be $R^2\times\LL_{\ge 0}^2$. 
    By Lemmas \ref{lemma:RTypeI}--\ref{lemma:qtrueRdir},
    we may apply Lemma \ref{lemma:noZdir} with $Z=R^2$, $Y=\LL_{\ge 0}$
    to find a true $\LL_{\ge 0}$-direction in $R^2\times\LL_{\ge 0}\times\{0\}$ which is not an $R$-direction.
    Similarly, there is a true $R$-direction in $\{q\}\times R\times\LL_{\ge 0}^2$ (for any $q\in R$) which is not an $\LL_{\ge 0}$-direction,
    by applying Lemma \ref{lemma:noZdir} with $Z=\LL_{\ge 0}^2$, $Y=R$.
\endproof

The proof of Theorem \ref{thm:notLdir} is very similar.
Let $T$ be the subset of $\LL_{\ge_0}\times[-1,1]$ given by taking a closed double sided topologist sine 
curve inside each $[\alpha,\alpha+1]$. That is,
$T$ contains $\{\alpha\}\times[-1,1]$ for each $\alpha$,
and curves inside $(\alpha,\alpha+1)\times[-1,1]$ oscillating as $\sin(1/x)$ and accumulating on both
$\{\alpha\}\times[-1,1]$ and $\{\alpha+1\}\times[-1,1]$.
Then $T$ is a closed subset of the $\omega$-bounded normal space $\LL_{\ge_0}\times[-1,1]$ and inherits these properties.
Denote $T\cap[0,\alpha)$ by $T_\alpha$.

\begin{lemma}\label{lemma:notpathconn}
   $T$ is a true $T$- and $\LL_{\ge 0}$-direction in itself. The image of
   any map $\LL_{\ge 0}^2\to T$ is contained in $T_{\alpha+1}-T_\alpha$ for some $\alpha$.
   Given $f:T\to T$, then $f(\wb{T_{\beta+1}}-T_\beta)$ is contained $\wb{T_{\alpha+1}}-T_\alpha$ for some $\alpha$.
\end{lemma}
\proof 
  It is easy to see that $T$ is an $\LL_{\ge 0}$-direction, and thus a $T$-direction by Lemma \ref{lemma:abcd} (c).
  The other claims are immediate.
\endproof

\proof[Proof of Theorem \ref{thm:notLdir}]
  By Lemma \ref{lemma:notpathconn},  $Z=\LL_{\ge 0}$ and $Y=T$ satisfy the assumptions of Lemma \ref{lemma:noZdir}. 
\endproof

\begin{prob}
   Are there regular Type I spaces $X,Y$, with $Y$ path connected, and a club $D\subset X$ such that
   $D$ is an $Y$-direction but not an $\LL_{\ge 0}$-direction in $X$~?
\end{prob}

\bibliographystyle{plain}
\bibliography{../biblio}

{\footnotesize
\vskip .4cm
\noindent
Mathieu Baillif \\
Haute \'ecole du paysage, d'ing\'enierie et d'architecture (HEPIA) \\
Ing\'enierie des technologies de l'information TIC\\
Gen\`eve -- Suisse
}

\end{document}